\documentclass[11pt,twoside]{amsart}
\usepackage{mathrsfs}
\usepackage{txfonts}
\usepackage{bbm}
\usepackage{amsfonts}
\usepackage{amsmath}
\numberwithin{equation}{section}
\usepackage{esint}
\usepackage{amssymb}
\usepackage{wasysym}
\usepackage{appendix}
\usepackage{psfrag}
\usepackage{graphics,epsfig,amsmath}
\usepackage{comment}
\usepackage{enumerate}
\usepackage{times}
\usepackage{fancyhdr}
\usepackage{cite}
\usepackage{mathrsfs}

\newtheorem{Theorem}{Theorem}[section]
\newtheorem{Lemma}[Theorem]{Lemma}
\newtheorem{Proposition}{Proposition}[section]

\newtheorem{remark}[Theorem]{Remark}

\usepackage{metalogo}
\usepackage[cmyk]{xcolor}
\selectfont
\usepackage{mathtools}
\mathtoolsset{showonlyrefs}
\renewcommand{\epsilon}{\varepsilon}
\mathtoolsset{showonlyrefs}
\newcommand{\R}{\mathbb{R}}

\renewcommand{\leq}{\leqslant}
\renewcommand{\le}{\leqslant}
\renewcommand{\geq}{\geqslant}
\renewcommand{\ge}{\geqslant}

\allowdisplaybreaks
\usepackage{geometry}
\title[Uniqueness and Nondegeneracy of Positive Solutions]
{Uniqueness of positive solutions of  mixed Laplacian and fractional Laplacian in $\R$}

\author[X. Su]{Xifeng Su}
\address{School of Mathematical Sciences, Laboratory of Mathematics and Complex Systems (Ministry of Education)\\
	Beijing Normal University,
	No. 19, XinJieKouWai St., HaiDian District, Beijing 100875, P. R. China}
\email{xfsu@bnu.edu.cn, billy3492@gmail.com}

\author[C.Zhang]{Chengxiang Zhang}
\address{School of Mathematical Sciences, Laboratory of Mathematics and Complex Systems (Ministry of Education)\\
	Beijing Normal University,
	No. 19, XinJieKouWai St., HaiDian District, Beijing 100875, P. R. China}
\email{zcx@bnu.edu.cn}

\author[J. Zhang]{Jiwen Zhang}
\address{Academy of Mathematics and Systems Science\\
	Chinese Academy of Sciences,
	No. 55 Zhongguancun East Road, HaiDian District, Beijing 100190, P. R. China}
\email{jwzhang@amss.ac.cn}


\thanks{X. Su is supported by the National Natural Science Foundation of China (Grant No. 12371186,12571207). C. Zhang is supported by the National Natural Science Foundation of China (Grant No. 12371107). J. Zhang is supported by the China Postdoctoral Science Foundation (Grant No. 2025M783149).}

\begin{document}

	\begin{abstract}
    We establish the uniqueness  and nondegeneracy for
positive finite-energy solutions of the following one-dimensional mixed
Schr\"odinger equation
\begin{equation}
    -Q'' + (-\Delta)^s Q + Q = |Q|^p Q
    \quad \text{in } \mathbb{R},
\end{equation}
where $0<s<1$ and $p>0$. 

For every $0<s<1, p>0$, we obtain the corresponding results for positive solutions whose linearized operators have Morse index equal to one. We extend the nodal analysis developed by Frank-Lenzmann-Silvestre  to the mixed operator  establishing the monotonicity of
the quotient of a hypothetical nontrivial even kernel element
by $-Q'$. We then combine this monotonicity with oscillation estimates to prove nondegeneracy, and 
uniqueness follows by continuation.

Moreover, for every $1/2\le s<1,p>0$, both uniqueness and nondegeneracy hold for  positive solutions.
Our novelty is to use the properties of Green kernel of the mixed operator to
derive the upper bound and lower bound on the maximum of each solution. This
peak estimate allows us to prove uniqueness by direct comparison
and to establish nondegeneracy.
	\end{abstract}
	
		\maketitle	
	
	\tableofcontents
	
	\section{Introduction}

	The main purpose of this paper is to derive some uniqueness and nondegeneracy results for positive solutions of the following Schr\"{o}dinger equation driven by  the mixed local/nonlocal operator:
	 \begin{equation}\label{main equation}
	 \begin{cases}
	 		-  Q'' +  (-\Delta)^s Q+Q= |Q|^{p}Q \quad \hbox{in $\mathbb{R}$,}\\
	 		Q\in H^1(\mathbb{R}),  \quad u\geqslant 0.
	 \end{cases}
	  \end{equation}
   Here, the fractional Laplacian $ (-\Delta )^{s} $ is defined as
  \begin{equation}\label{fractional laplacian}
  	(-\Delta )^{s} u(x)=c_{1,s}\,\text{PV}\int_{\mathbb{R}} \frac{u(x)-u(y)}{|x-y|^{1+2s}}\, dy,
  \end{equation}
  where $ c_{1,s} $ is a dimensional constant that depends on $n$ and $s$, precisely given
\begin{equation}\label{c_{1.s}}
    c_{1,s}=\left(\int_{{\mathbb{R}}} \frac{1-\cos(\xi_1)}{|\xi|^{1+2s}}\, d\xi\right)^{-1}=\frac{2^{2s}s\Gamma(s+\tfrac12)}
{\sqrt{\pi}\,\Gamma(1-s)}.
\end{equation}

We consider standing waves of the focusing nonlinear
Schr\"odinger equation with mixed dispersion
\begin{equation}\label{eq:nls}
 i\partial_t\psi
 =-\partial_x^2\psi+(-\partial_x^2)^s\psi-|\psi|^p\psi.
\end{equation}
By seeking solutions of the form
$\psi(t,x)=e^{i\omega t}Q_\omega(x)$, we obtain
\begin{equation}\label{eq:frequency}
 -Q_\omega''+(-\partial_x^2)^sQ_\omega+\omega Q_\omega
 =Q_\omega^{p+1}.
\end{equation}
Note that, we can deduce equation~\eqref{main equation} by taking $\omega=1$ in equation~\eqref{eq:frequency}. 

We are interested in the uniqueness and nondegeneracy of positive
standing-wave profiles. Uniqueness determines the profile up to
translation, while nondegeneracy means that the kernel of the real
linearized operator is spanned by $Q'$. After fixing the center,
we can use nondegeneracy to construct a local family of profiles
as the frequency varies. For related results on solitary waves
of the generalized Benjamin-Ono and Benjamin-Bona-Mahony
equations, we refer to~\cite{MR3070568}.

These properties are also important in the study of stability,
as in the work of Weinstein and Grillakis, Shatah, and Strauss
\cite{MR783974,MR901236,MR1081647}.
To prove orbital stability, we must also establish well-posedness
of the evolution, verify the remaining spectral assumptions,
and determine the sign of
$\frac{d}{d\omega}\|Q_\omega\|_2^2$.
For equations with a single homogeneous dispersion term, we can
compute the frequency dependence of the mass by rescaling a
fixed profile. In the mixed equation, however, a rescaling that
normalizes the frequency changes the relative coefficients
of the local and fractional terms. The mass derivative
therefore also depends on how the normalized profile varies
with the frequency, and its sign requires a separate analysis.

A well-developed theory of uniqueness and nondegeneracy
is available for equations with a single dispersive order. For the classical Laplacian, the uniqueness theory have been established in
 \cite{Coffman72, MR886933, MR969899, MR1201323}, and these  conclusions were extended into the nonlocal case in \cite{AT91, MR3070568, MR3207007, MR3530361} and references therein. In the mixed setting, in~\cite{SZZ2026}, the authors
proved uniqueness and nondegeneracy of least-energy solutions in arbitrary dimension in the Sobolev subcritical
range, when $s$ is sufficiently close to $0$ or $1$. These endpoint results raise the question of whether the same conclusions hold at intermediate fractional orders and for positive solutions beyond variational ground states.

To address these questions, we need an argument that accounts for both kinetic terms. The nonlocal term prevents a direct application of the classical
ODE shooting method. Although the fractional term is of lower differential order, it has a fixed coefficient and governs the algebraic decay of
positive solutions at infinity. Hence, it cannot simply be treated
as a small perturbation, and its interaction with the local
operator must be taken into account.

 Our first main result establishes uniqueness and nondegeneracy
among positive finite-energy solutions of~\eqref{main equation}
whose linearized operators have Morse index one, for every
$0<s<1$ and $p>0$. The proof extends the nodal analysis of Frank-Lenzmann-Silvestre~\cite{MR3070568,MR3530361}
to the mixed operator. Under the spectral assumption, the resulting nodal information yields monotonicity of the
quotient of a hypothetical nontrivial even kernel element
by $-Q'$ on the positive half-line.
Together with a dilation identity, this monotonicity excludes nontrivial even kernel elements and establishes nondegeneracy. Uniqueness within this class then follows by continuation.

When $\frac12\le s<1$, a different argument allows us to remove
the Morse index assumption and establish uniqueness and
nondegeneracy for all positive finite-energy solutions and
every $p>0$. Consequently, each such solution has Morse index
one. We derive a new explicit peak estimate from the properties
of the Stieltjes density associated with the Green kernel
and half-line Hamiltonian identities. With this bound at hand, we compare two positive solutions directly through a Hamiltonian identity and a strict Gram inequality,
without using Sturm oscillation theory.
Uniqueness among all positive solutions for $0<s<1/2$ remains open in general.

Our main methodological contribution is a new comparison argument
that exploits the interaction between the local and fractional
operators. By combining Green kernel estimates, Hamiltonian
identities, and spectral arguments, we account for the different
scaling of the two kinetic terms and establish uniqueness beyond
the range covered by endpoint perturbation methods.

\bigskip

To well formulate the main results in this paper,	let us first introduce the  functional framework.  
		
		Consider  the classical Sobolev space $H^1(\mathbb{R})$  with the following equivalent norm 
	\[ \|u\|_{s}^2:=	\frac12 \int_{\mathbb{R}}\int_{\mathbb{R}} \frac{|u(x)-u(y)|^2}{|x-y|^{1+2s}}\,dx dy+\int_{{\mathbb{R}}} | u'|^2\,dx+\int_{{\mathbb{R}}} u^2\, dx \quad \text{for any } u\in H^1(\mathbb{R}).\]
	We define a weak notion of solution $u\in H^1(\mathbb{R})$ for equation~\eqref{main equation} by the identity
	\begin{equation}\label{weak solution}
		\frac12 \int_{\mathbb{R}}\int_{\mathbb{R}} \frac{(u(x)-u(y))(\varphi(x)-\varphi(y))}{|x-y|^{1+2s}}\,dx dy+\int_{{\mathbb{R}}}  u'\varphi'\,dx+\int_{{\mathbb{R}}} u\,\varphi\, dx=\int_{{\mathbb{R}}}|u|^{p}u\varphi\, dx\quad \text{for any } \varphi\in H^1(\mathbb{R}).
	\end{equation}
	In addition, we define  the scalar product by
	\begin{equation}\label{eq: inner product}
		\left\langle u,v \right\rangle_{s}:=\frac12 \int_{\mathbb{R}}\int_{\mathbb{R}} \frac{(u(x)-u(y))(v(x)-v(y))}{|x-y|^{1+2s}}\,dx dy+\int_{{\mathbb{R}}}  u' v'\,dx+\int_{{\mathbb{R}}} u\,v\, dx,\quad\text{for any } u,v\in H^1(\mathbb{R}).
	\end{equation}
    
We next introduce the corresponding minimization problem of the  energy functional of 
	\[ 
	\mathscr{F}_s(u):=\frac{1}{2}\|u\|_s^2-\frac{1}{p+2}\int_{{\mathbb{R}}}|u|^{p+2}\, dx \qquad \forall u\in H^1(\mathbb{R})
	\]
	over the Nehari manifold
	\[ \mathcal{B}_s:=\left\{u\in H^1(\mathbb{R}):u\not\equiv 0 \text{ and }\, \|u\|_s^2=\|u\|^{p+2}_{L^{p+2}(\R)} \right\}. \]

	One may observe that, by some suitable rescaling, to find a minimizer of the above minimization problem of $\mathscr{F}_s$ over $\mathcal{B}_s $ could be transformed into finding
	a minimizer $u_s\in H^1(\mathbb{R})$ of a minimization problem of the following functional 
	\begin{equation}\label{j}
		\mathscr J_s(v):=\frac{\|v\|_s^2}{\|v\|^{2}_{L^{p+2}(\R)}}  \qquad \text{ for any } v\in H^1(\mathbb{R}) \text{ with } v\not\equiv 0. 
	\end{equation}
	Note that  the functional $\mathscr J_s$ has indeed minimizers  $u_s\in H^1(\mathbb{R})$ due to  the concentration-compactness principle (see e.g. \cite{MR778974}).

From~\cite{DSVZ24,SZZ2026}, the following existence result along with fundamental properties of nonnegative solutions for equation~\eqref{main equation} can be summarized below:   
\begin{Theorem}\label{th: existence and properity} 
Let $s\in(0,1)$, and $p>0$. Then the following holds:
	\begin{itemize}
		\item[(i)]\text{EXISTENCE}: There exists a nonnegative solution  $Q\in H^1(\mathbb{R})$ of equation \eqref{main equation}, which is a minimizer of the above functional  $\mathscr J_s$.
		\item[(ii)] REGULARITY, SYMMETRY AND DECAY: If $Q\in H^1(\mathbb{R})$ with $Q\geqslant 0$ and $Q\not\equiv 0$ solves equation~\eqref{main equation}, then 
		\begin{itemize}
		\item there exists some $x_0\in\mathbb{R}$ such that $Q(\cdot-x_0)$ is even, positive, and non-increasing in $|x-x_0|$. 
		\item Moreover, the function $Q\in C^{2}(\mathbb{R})\cap H^2(\mathbb{R})$ and there exist two suitable constants $C_2\geqslant C_1>0$ such that 
		\[ \frac{C_1}{|x|^{1+2s}}\leqslant Q(x)\leqslant  \frac{C_2}{|x|^{1+2s}} \quad \forall |x|>1. \]
		\end{itemize}
	\end{itemize}
\end{Theorem}

It would then be interesting to introduce  the following notion for  those positive solutions of~\eqref{main equation} obtained in  Theorem~\ref{th: existence and properity}. We always center a positive solution at its symmetry point.
Let us denote the set of positive finite-energy solutions of equation~\eqref{main equation} for $s\in[0,1]$ by 
\begin{equation}\label{m_s}
	\mathcal{M}_s:=
	\left\{\begin{array}{ll}
	Q\in H^1(\mathbb{R}):   Q \text{ solves equation } \eqref{main equation} \text{ with } Q>0, \text{ and } Q(0)=\max_{x\in\
    \mathbb R}Q(x)
	\end{array}
	\right\}.
\end{equation}

Then, for any $Q\in\mathcal{M}_s$, we  denote  the corresponding linearized operator at $Q$ by
\begin{equation}\label{Ls}
	L_{s,Q}= - \Delta  +  (-\Delta)^s +1 - (p+1)Q^{p}
\end{equation}
acting on $L^2(\mathbb{R})$. Note that
 $L_{s,Q}$ is a self-adjoint operator on $L^2(\mathbb{R})$ with quadratic-form domain $H^1(\mathbb{R})$ and
 operator domain $H^2(\mathbb{R})$.

Moreover, for fixed $p>0$, we introduce a subclass of $\mathcal M_s$
consisting of solutions whose linearized operators have Morse index equal to one:
\[
\mathcal G_s
:=
\left\{
Q\in\mathcal M_s:
L_{s,Q}\text{ has Morse index one}
\right\}.
\]
 By~\cite[Proposition~A.1]{SZZ2026}, the class is nonempty, since every global minimizer of $\mathscr J$ belongs to $\mathcal G_s$.

Our first result establishes uniqueness  and nondegeneracy within $\mathcal G_s$ for every $s\in(0,1)$.

\begin{Theorem}
\label{th:nondegeneracy low s}
Let $0<s<1$ and $p>0$.  Then $\mathcal G_s$ consists of
a single solution $Q_s$, and its linearized operator satisfies
\[
\ker L_{s,Q_s}
=\operatorname{span}\{Q_s'\}.
\]
\end{Theorem}

Note that, the spectral condition provides the nodal information needed for this argument, and our conclusions are therefore restricted to
solutions in $\mathcal G_s$. 

To extend these conclusions to
all positive finite-energy solutions, we develop a different strategy based on Green kernel estimates and half-line
Hamiltonian identities, obtaining the following result without a spectral assumption, although the range of
fractional orders is restricted to $\frac12\le s<1$.

\begin{Theorem}\label{th:nondegeneracy s>1/2}
	Let $s\in [1/2,1) $ and $p>0$.  Then $ \mathcal{M}_s$ is a singleton.
\end{Theorem}

 \begin{remark}
 \begin{itemize}
     \item [(i)] From this uniqueness result, we know that the positive solution $Q$ of~\eqref{main equation} has Morse index equal to one, i.e.,  the  corresponding linearized operator $L_{s,Q}$ has exactly one strictly negative eigenvalue. Moreover, in view of Theorem~\ref{th:nondegeneracy low s}, we derive that $L_{s,Q}$ is nondegenerate; i.e., 
	 \[  \ker L_{s,Q}=\operatorname{span}\left\{Q'\right\}. \]
     \item[(ii)] Both uniqueness and nondegeneracy rely essentially on the
peak estimate
\begin{equation}\label{eq:upperpeak}
\frac{p+2}{2}<Q^p(0)<p+1.
\end{equation}
We prove these two results independently, without using
either conclusion to establish the other.
 \end{itemize}
\end{remark}

We point out that the upper bound in~\eqref{eq:upperpeak} is invalid when $s$ is sufficiently small. Indeed, assume that $Q\in \mathcal M_s$, with $p>0$ fixed, then
\[
 Q_s^p(0)\longrightarrow p+2>p+1
 \qquad\text{as }s\downarrow0.
\]
Hence, to extend both conclusions to all positive finite-energy solutions throughout $0<s<1$, we must therefore modify or replace the peak condition, rather than merely sharpen
the kernel estimates.
\medskip

We now turn to the stationary mixed local and nonlocal Schr\"odinger equation~\eqref{eq:frequency} with frequency $\omega>0$.
Note that the scaling  $Q(x)=\omega^{1/p}u(\sqrt{\omega}\,x)$ reduces this equation to the normalized problem with coupling coefficient
$\kappa=\omega^{s-1}$. Hence, as a straightforward consequence of Theorems~\ref{th:nondegeneracy low s}
and~\ref{th:nondegeneracy s>1/2}, we obtain the following result.

\begin{Theorem}\label{thm:corollary}
Let $0<s<1$, $p>0$, and $\omega>0$.
Equation~\eqref{eq:frequency} admits, up to translation, a unique
positive solution $Q_\omega\in H^1(\mathbb R)$ of Morse index one.
Moreover, this solution is nondegenerate:
\[
\ker L_{\omega,s,Q_\omega}=\operatorname{span}\{Q_\omega'\},
\]
where $L_{\omega,s,Q_\omega}
:=-\partial_x^2+(-\partial_x^2)^s+\omega-(p+1)Q_\omega^p.$

If, in addition, $s\ge1/2$ and $\omega\ge1$,  all positive finite-energy solutions  are unique up to translation and nondegenerate.
\end{Theorem}

\begin{remark}
The scaling preserves both the Morse index and the dimension of
the kernel. Since uniqueness and nondegeneracy in the Morse index equal to one class hold for the normalized problem for every $\kappa>0$, the same
conclusions hold for equation~\eqref{eq:frequency} at every
$\omega>0$.

The argument for all positive solutions, however, relies on a Green kernel comparison that, for $s\ge1/2, \omega\ge1$, yields the peak bound \[
\frac{\omega(p+2)}2<Q^p(0)<\omega(p+1).
\]
This bound allows us to prove uniqueness and nondegeneracy for all
positive finite-energy solutions. 
However, when $0<\omega<1$, the corresponding coupling satisfies $\kappa>1$, and the required Green kernel inequality is invalid when $\omega$ is sufficiently small. Thus, the same argument no longer provides the peak bound needed to establish the above both conclusions.
\end{remark}

\medskip


The paper is organized as follows.
In Section~\ref{sec:Green kernels and common analytic identities},
we establish basic estimates and identities for the convolution
kernel of the resolvent
$(-\Delta+(-\Delta)^s+1)^{-1}$.
We also derive a uniform lower bound on the peak values of
positive finite-energy solutions by using a monotonicity formula
for the associated extension problem in the upper half-plane
$\mathbb R^2_+$. Section~\ref{sec: An oscillation estimate} introduces a sharp oscillation estimate for the second eigenfunction of mixed Schr\"odinger operator, which provides us with a bound in agreement with Sturm–Liouville theory for ODE. 

In Section~\ref{sec:uniqueness}, we prove the uniqueness up to translations and nondegeneracy properties for $0<s<1$ and every $p>0$ within the class $\mathcal G_s$.
Section~\ref{sec:Nondegeneracy s ge 1/2} establishes the corresponding results   of all positive finite-energy
solutions  for $1/2\leq s<1$ and every
$p>0$, as stated in
Theorem~\ref{th:nondegeneracy s>1/2}. 
The appendix contains several technical results and their proofs
used in the main arguments.


\section{Green kernels and a lower sharp bound}
\label{sec:Green kernels and common analytic identities}

In this section, we investigate some properties of the Green kernel associated with the mixed operator and derive the lower bound in~\eqref{eq:upperpeak} using the extension method.
\medskip

We begin by introducing the Green kernel. For $s\in(0,1)$, let $G_s$ denote the convolution kernel of $(-\Delta+(-\Delta)^s+1)^{-1}$, given by
\begin{equation}\label{eq:bsdefinition}
G_s(x)=\frac{1}{2\pi}\int_{\mathbb R}
\frac{e^{ix\xi}}{1+\xi^2+|\xi|^{2s}}\,d\xi,
\qquad
G_s(0)=\frac{1}{\pi}\int_0^\infty
\frac{d\xi}{1+\xi^2+\xi^{2s}}.
\end{equation}

\subsection{The properties of Green kernel}\label{sec:The value of Green kernel at the origin}

\begin{Lemma}\label{lemma: G_s definition }
    Let $G_s$ be the convolution kernel of $(- \Delta +  (-\Delta)^s +1)^{-1}$.
    Let $r>0$, define
\begin{equation}\label{rho_s representation}
    \rho_s(r)=\frac{\sin(\pi s)}{\pi}
\frac{r^s}
{\left(1-r+r^s\cos(\pi s)\right)^2
+r^{2s}\sin^2(\pi s)}.
\end{equation}
Then \begin{equation}\label{eq:green mixture}
G_s(x)=\int_0^\infty\frac{\rho_s(r)}{2\sqrt r}e^{-\sqrt r|x|}\, dr
=\int_0^\infty \rho_s(\lambda^2)e^{-\lambda|x|}\, d\lambda.
\end{equation}
Moreover, we have
\begin{equation}\label{eq:green properties}
\int_{\R}G_s=1,\quad G_s'(x)=-\lambda G_s(x) \quad \forall  x>0, \quad
\lim_{x\to 0_+}G_s'(x)=-\frac12,\quad \text{ and }\quad 
\|G'\|_{L^\infty(\R)}\le\frac12.
\end{equation}
In particular,  we define  \begin{equation}\label{eq:K mixture}
    K_s(x):=\int_0^\infty \lambda^2\rho_s(\lambda^2)e^{-\lambda|x|}\, d\lambda
\qquad \text{for every } x\neq0,
\end{equation}  then 
\begin{equation}
    G_s''=K_s-\delta_0.
\end{equation}
Here, the second derivative of $G$ is a distribution, that is, its Dirac mass records the jump of $G_s'$ at the origin.
\end{Lemma}
\begin{proof} From Proposition~\ref{prop: rho_s}, it follows that $\rho_s(r)>0$, and
\begin{equation}\label{rho_s identity}
\frac{1}{1+x+x^s}
=\int_0^\infty\frac{\rho_s(r)}{x+r}\,dr,
\qquad x\ge0.
\end{equation} 
In particular,
\begin{equation}\label{rho_s estimate}
    \int_0^\infty\frac{\rho_s(r)}r\, dr=1,
\qquad
\int_0^\infty\rho_s(r)\, dr=1.
\end{equation}
 Recall that the one-dimensional resolvent of $-\Delta+r$ has kernel
$e^{-\sqrt r|x|}/(2\sqrt r)$, that is,
\[\frac1{2\pi}\int_{\R}\frac{e^{ix\xi}}{|\xi|^2+r}\, d\xi=\frac{e^{-\sqrt r|x|}}{2\sqrt r}.\]
As a result of this and~\eqref{eq:bsdefinition}, we obtain the identity~\eqref{eq:green mixture}.   A direct calculation implies the properties of $G$ in~\eqref{eq:green properties}.
\end{proof}

\begin{Proposition}
\label{prop:green tail}
For $0<s<1$, let $G_s$ be the convolution kernel of $(- \Delta +  (-\Delta)^s +1)^{-1}$. Then
\begin{equation}\label{eq:green tail}
\lim_{|x|\to\infty}
|x|^{1+2s}G_s(x)
=c_{1,s},
\end{equation}
where $c_{1,s}$ be introduced in~\eqref{c_{1.s}}.  Moreover, for some $C_G>0$,
\begin{equation}\label{eq:green bound}
0<G_s(x)
\leq C_G(1+|x|)^{-1-2s}
\qquad\forall x\in\R.
\end{equation}
\end{Proposition}

\begin{proof} By the definition of $\rho_s$ in~\eqref{rho_s representation}, we know that $\rho_s(\lambda^2)\leq C_s\lambda^{2s}$ for some constant $C_s>0$, and 
\[
\rho_s(\lambda^2)
\sim
\frac{\sin(\pi s)}{\pi}\lambda^{2s}
\qquad \lambda\to 0, \quad \text{and}\quad \rho_s(\lambda^2)
\sim
\frac{\sin(\pi s)}{\pi}\lambda^{2s-4}
\qquad \lambda\to\infty. 
\]
Thus $\rho_s(\lambda^2)\in L^1(0,\infty)$.

Moreover, by~\eqref{eq:green mixture} and the change of variables
$v=r\sqrt{t}$, we obtain, for every $r>0$,
\begin{equation}\label{asdfsda}
r^{1+2s}G_s(r)
=
\int_0^\infty e^{-v}r^{2s}\rho_s(v^2/r^2)\,dv
\le C_s\int_0^\infty e^{-v}v^{2s}\,dv
<\infty.
\end{equation}
For every fixed $v>0$,
\[
r^{2s}\rho_s(v^2/r^2)
\longrightarrow
\frac{\sin(\pi s)}{\pi}v^{2s}
\qquad r\to\infty.
\]
By the dominated convergence theorem, we obtain 
\[
\lim_{r\to\infty}r^{1+2s}G_s(r)=
\frac{\sin(\pi s)}{\pi}
\int_0^\infty e^{-v}v^{2s}\,dv=
\frac{\Gamma(1+2s)\sin(\pi s)}{\pi}
=
c_{1,s}.
\]
This ends the proof.
\end{proof}

\begin{Proposition}\label{prop: Q'>0}
  Let $Q\in \mathcal M_s$, then $Q'(r)<0$  for all $r>0$.
\end{Proposition}
\begin{proof} 
    Since $Q\in C^2(\R)$ solves~\eqref{main equation}, we observe that 
   \begin{equation}\label{eq:resolvent profile}
Q=G_s*Q^{p+1},\qquad (1+(-\Delta)^s)Q=K_s*Q^{p+1},\qquad Q''=K_s*Q^{p+1}-Q^{p+1}.
\end{equation} 
In view of Theorem~\ref{th: existence and properity}, we know that $Q$ is even and $Q'(r)\leq 0$ for all $r>0$ and $Q'\not\equiv0$. Since $|x-y|<x+y$ for every $x,y>0$, thus by  differentiating the convolution $Q=G_s*Q^{p+1}$, we see that
\begin{equation}\label{eq:strict decrease}
Q'(x)=(p+1)\int_0^\infty
\left[G_s(|x-y|)-G_s(x+y)\right]Q^p(y)Q'(y)\,dy<0\qquad \forall x>0.
\end{equation}
This ends the proof.
\end{proof}

We notice that for every $\xi\in\R$, $1+\xi^2\le1+\xi^2+|\xi|^{2s}\le 2(1+\xi^2)$, recalling the definition of $G_s(0)$ in~\eqref{eq:bsdefinition}, one thus has that $1/4\le G_s(0)\le1/2$.
\begin{Proposition}\label{prop: property of G_s}
    The function $s\mapsto G_s(0)$  is  continuous and strictly concave in $[0,1]$. In particular,
\begin{equation}\label{eq:G_s values}
G_0(0)=G_1(0)=\frac1{2\sqrt2},\qquad \qquad G_s(0)>\frac1{2\sqrt2}\qquad \forall 0<s<1.
\end{equation}
\end{Proposition}
\begin{proof}
  For fixed   $\xi>0$, a direct
differentiation yields
\[
\frac{\partial^2}{\partial s^2}\frac1{1+\xi^2+|\xi|^{2s}}
=\frac{4|\xi|^{2s}(\ln \xi)^2(|\xi|^{2s}-1-{\xi^2})}{(1+\xi^2+|\xi|^{2s})^3}<0\qquad  \text{for every }\xi\ne1.
\]
This implies that $G_s(0)$ is strictly concave. Since $(1+\xi^2)^{-1}$ is integrable, the dominated convergence theorem proves the continuity at the endpoints. Thus, by the strict concavity of $G_s(0)$, we see that 
\begin{equation}\label{G_s(0)}
    G_s(0)>(1-s) G_{0}(0)+sG_1(0)=\frac1{2\sqrt2}\qquad \forall 0<s<1.
\end{equation}
This ends the proof.
\end{proof}

\subsection{A energy extension and a lower sharp bound}
In this section, we aim to establish a uniform lower bound on the peak values of positive finite-energy solution.

We start by recalling the extension principle in \cite{MR2354493} that expresses the nonlocal operator $(-\Delta)^s$ on $\mathbb{R}$ with $s\in(0,1)$ as a generalized Dirichlet-Neumann map for a suitable local elliptic problem posed on the upper half-space $\mathbb{R}^{2}_+=\left\{(x,t):\ x\in\mathbb{R},\ t>0\right\}$. That is,  for a measurable $f:\mathbb{R}\rightarrow \mathbb{R}$,
	we  define $\mathscr{E}_s f$ as the $s$-Poisson extension of $f$ to   $\mathbb{R}_+^{2}$ by
	\begin{equation}\label{keydf}
		\left(\mathscr{E}_s f\right)(x,t)=\int_{\mathbb{R}}P_s(x-y,t)f(y)\, dy.
	\end{equation}
	Here the generalized Poisson kernel $P_s(\tau,t)$ of order $s$ is given by 
	\begin{equation}\label{eq:poisson-kernel}
		P_{s}(\tau,x)=
\frac{\Gamma(s+1/2)}{\sqrt\pi\,\Gamma(s)}
\frac{\tau^{2s}}{(\tau^2+x^2)^{(1+2s)/2}}.
	\end{equation}	
 Under suitable assumptions on $f$, it is known that $\left(\mathscr{E}_s f\right)(x,t)$
	 solves the degenerate elliptic boundary value problem
	\begin{equation}
		\begin{cases}\label{eq:extension energy}
			\text{div }(t^{1-2s}\nabla \left(\mathscr{E}_s f\right))=0\qquad &\text{ in } \mathbb{R}^{2}_+\\
			\left(\mathscr{E}_s f\right)(x,0)=f &\text{ in } \mathbb{R}\\
			-d_s t^{1-2s}\partial_t \left(\mathscr{E}_s f\right)|_{t\rightarrow 0}=(-\Delta)^sf
			&\text{ on }  \mathbb{R}.
		\end{cases}
	\end{equation}
	Here the positive constant $d_s=2^{2s-1}\frac{\Gamma(s)}{\Gamma(1-s)}$.
From~\cite[Proposition~3.1]{MR3070568}, it follows that
\begin{equation}\label{bdsjhkisdfs}
    d_s\int_{\R}\int_{0}^\infty
\tau^{1-2s}|\nabla (\mathscr{E}_s f)|^2\,\mathrm d\tau\,\mathrm dx
=\|(-\Delta)^{s/2}f\|_2^2=\frac{c_{1,s}}2\int_{\R}\int_{\R}\frac{|f(x)-f(y)|^2}{|x-y|^{1+2s}}\, dxdy, \qquad \forall f\in H^1(\R).
\end{equation}

Let $f,g\in H^1(\R)\cap C_b^1(\R)\cap C^{1,1}_{\mathrm{loc}}(\R)$, we define the following energy and symmetric bilinear form
\begin{equation}\label{A_s}
    \mathcal A_s(f)
=d_s\int_0^\infty\tau^{1-2s}
|(\mathscr{E}_s f)_\tau(0,\tau)|^2\,d\tau,\quad 
\mathcal B_s(f,g)
=d_s\int_0^\infty\tau^{1-2s}
(\mathscr{E}_s f)_\tau(0,\tau) \,(\mathscr{E}_s f)_\tau(0,\tau)\,d\tau.
\end{equation}
From Lemma~\ref{lemma Poisson extension}, it follows that  
\begin{equation}\label{mathcal B_s}
    A_s(f)<\infty\quad \text{and}\quad |\mathcal B_s(f,g)|^2\leqslant A_s(f)\,A_s(g).
\end{equation}

\begin{Lemma}\label{lemma A_s>0}
Let $0<s<1$ and $u\in H^1(\R)\cap C_{\mathrm{loc}}^{1,1}(\R)$. Assume that $u(0)>0$.  
Then
\[
0<\mathcal A_s(u)<\infty.
\]
\end{Lemma}

\begin{proof}
Let $\mathscr{E}_s u$ be its $s$-extension, one has 
\[
(\mathscr{E}_s u)_\tau(0,\tau)
=
\frac12\int_{\R}
\partial_\tau P_s(z,\tau)
\bigl(u(z)+u(-z)-2u(0)\bigr)\,dz.
\]
From Lemma~\ref{lemma Poisson extension}, it follows that
\begin{equation}
    |(\mathscr{E}_s u)_\tau(0,\tau)|\leq \begin{cases}
        C_{s,u}\tau^{2s-1}\qquad &0<\tau\leq 1\\
        \frac{1+2s}{\tau}\|u\|_{L^\infty(\R)} &\tau\geq 1.
    \end{cases}
\end{equation}
This implies that $\mathcal A_s(u)$ is finite.

\medskip
 We now turn to show that $\mathcal A_s(u)>0$.

Suppose, argue by contradiction that
$\mathcal A_s(u)=0.$ 
Since \(d_s>0\) and \(\tau^{1-2s}>0\), the nonnegativity of the
integrand yields
$(\mathscr{E}_s u)_\tau(0,\tau)=0 $ for almost every $
\tau>0.$
We notice that, the function \((\mathscr{E}_s u)_\tau(0,\cdot)\) is continuous on \((0,\infty)\), therefore  $(\mathscr{E}_s u)_\tau(0,\tau)=0, \forall \tau>0.$
Hence there exists a constant \(c\) such that
\begin{equation} \label{eq:constant axis}
(\mathscr{E}_s u)(0,\tau)=c
\qquad\text{for all }\tau>0.
\end{equation}

\medskip

We claim that $c=u(0)>0$. 

Indeed, by  scaling, one has 
\[
(\mathscr{E}_s u)(0,\tau)
=
\int_{\R}P_s(r,1)u(-\tau r)\,dr.
\]
For each fixed \(r\), continuity of \(u\) at the origin implies that $u(-\tau r)\to u(0) $ as $\tau\to0.$

Moreover, we observe that
\[
|P_s(r,1)u(-\tau r)|
\leq
\|u\|_{L^\infty(\R)}P_s(r,1)\in L^1(\R).
\]
By combining dominated convergence theorem and the normalization of \(P_s\), we obtain 
\[
\lim_{\tau\downarrow0}(\mathscr{E}_s u)(0,\tau)
=
u(0)\int_{\R}P_s(r,1)\,dr
=
u(0).
\]
Together with \eqref{eq:constant axis}, this yields
\begin{equation}
\label{eq:positive constant}
c=u(0)>0.
\end{equation}

\medskip

On the other hand,  since $u\in H^1(\R)$,  the H\"{o}lder inequality implies 
\[
|(\mathscr{E}_s u)(0,\tau)|
\leq
\|P_s(\cdot,\tau)\|_{L^2(\R)}
\|u\|_{L^2(\R)} =\tau^{-1/2}\|P_s(\cdot,1)\|_{L^2(\R)}
\|u\|_{L^2(\R)}
\longrightarrow0
\qquad\text{as }\tau\to\infty.
\]
From this and~\eqref{eq:constant axis}, it follows that $c=0$, which  
contradicts~\eqref{eq:positive constant}.
Therefore, $\mathcal A_s(u)>0$.
\end{proof}

\begin{Proposition}\label{prop:lower bound}
   Let $s\in (0,1), p>0$ and  $Q\in \mathcal M_s$. We define 
   \begin{equation}\label{eq:Hamiltonian}
      \mathcal H_Q(x)
=\frac12Q'(x)^2-\frac12Q^2(x)+\frac{Q^{p+2}(x)}{p+2}+\frac{d_s}{2}\int_0^\infty
\tau^{1-2s}((\mathscr{E}_s Q)_x^2-(\mathscr{E}_s Q)_\tau^2)(x,\tau)\,d\tau.
\end{equation}
Then $\mathcal H_Q(x)\equiv 0$. 

Moreover, we have
\begin{equation}\label{vsdfvdsfcds}
    \mathcal A_s(Q)=
\frac{2Q^{p+2}(0)}{p+2}-Q^2(0).
\end{equation}
In particular, 
  \begin{equation}\label{eq:axis profile}
{Q^{p}(0)}>\frac{p+2}2.
\end{equation} 
\end{Proposition}
\begin{proof} Since $Q\in C^2(\R)\cap C_b^1(\R)$,  in view of Lemma~\ref{lemma Poisson extension}, we know that for every compact interval $I\subset\R$,
\begin{equation}\label{eq:weighted majorant}
\tau^{1-2s}\bigl(|(\mathscr{E}_s Q)_x|^2+|(\mathscr{E}_s Q)_\tau|^2\bigr)(x,\tau)
\leq
\begin{cases}
C_I\bigl(\tau^{1-2s}+\tau^{2s-1}\bigr),
&0<\tau\leq1,\\[2mm]
C_I\tau^{-1-2s},
&\tau\geq1,
\end{cases}
\qquad x\in I,
\end{equation}
 and
\begin{equation}\label{vsdknvmlds12}
    \sup_{x\in I}
\int_0^\infty
\tau^{1-2s}
\left(
|(\mathscr{E}_s Q)_x (x,\tau)|^2
+
|(\mathscr{E}_s Q)_\tau(x,\tau)|^2
\right)\,d\tau
<\infty.
\end{equation}
Define
\[
J(x)=\frac{d_s}{2}\int_0^\infty
\tau^{1-2s}\bigl[(\mathscr{E}_s Q)_x^2-(\mathscr{E}_s Q)_\tau^2\bigr](x,\tau)\,d\tau
\]
and, for \(0<\epsilon<1<T\),
\[
J_{\epsilon,T}(x)=
\frac{d_s}{2}\int_\epsilon^T
\tau^{1-2s}\bigl[(\mathscr{E}_s Q)_x^2-(\mathscr{E}_s Q)_\tau^2\bigr](x,\tau)\,d\tau.
\]
We notice that the integral defining $J(x)$ is absolutely convergent.
Moreover, by~\eqref{eq:weighted majorant}, it follows that
\begin{equation}\label{eq:uniform truncation}
\sup_{x\in I}|J(x)-J_{\epsilon,T}(x)|
\leq
C_I\left(
\frac{\epsilon^{2-2s}}{2-2s}+
\frac{\epsilon^{2s}}{2s}+
\frac{T^{-2s}}{2s}
\right).
\end{equation}
Thus,  $J_{\epsilon,T}$ converges to $J$ locally uniformly in $x$.

Owing to~\eqref{eq:weighted majorant}, together with continuity of the
kernel derivatives for $\tau>0$, we have  that $J$ is continuous.

\medskip

We claim that $J'(x)=-(-\Delta)^sQ(x)Q'(x)$ is continuous.

Indeed, we first prove that
\begin{equation}\label{eq:horizontal trace}
(\mathscr{E}_s Q)_x(\cdot,\epsilon)\longrightarrow Q'
\qquad\text{locally uniformly as }\epsilon\to0.
\end{equation}
We choose a compact interval $I$, and set
$I_1=I+[-1,1]$.
For $0<\delta<1$, define
\[
\omega_I(\delta)=
\sup\left\{
|Q'(a)-Q'(b)|:
a,b\in I_1,\ |a-b|\leq\delta
\right\}, \qquad \omega_I(\delta)\to0\text{ as } \delta\to 0,
\]
thanks to the continuity of $Q'$ on the compact interval $I_1$.

By combining positivity with normalization of $P_s$, we obtain
\begin{equation}
    \begin{split}
        |(\mathscr{E}_s Q)_x(x,\epsilon)-Q'(x)|
&\leq
\int_{|z|\leq\delta}
P_s(z,\epsilon)|Q'(x-z)-Q'(x)|\,dz+
\int_{|z|>\delta}
P_s(z,\epsilon)|Q'(x-z)-Q'(x)|\,dz\\
&\leq
\omega_I(\delta)
+
2\|Q'\|_{L^\infty(\R)}
\int_{|z|>\delta}P_s(z,\epsilon)\,dz\\
&\leq \omega_I(\delta)
+2\|Q'\|_{L^\infty(\R)}\frac{c_s}{s}\left(\frac{\epsilon}{\delta}\right)^{2s}.
    \end{split}
\end{equation}
Taking $\epsilon\to0$, with $\delta$ fixed, and then letting $\delta\to 0 $, it follows~\eqref{eq:horizontal trace}.

\medskip

We next prove
\begin{equation}\label{eq:flux limit}
-d_s\epsilon^{1-2s}(\mathscr{E}_s Q)_\tau(\cdot,\epsilon)
\longrightarrow (-\Delta)^sQ
\qquad\text{locally uniformly as }\epsilon\to0.
\end{equation}
Since $ \int_{\R}\partial_\epsilon
P_s(z,\epsilon)\,dz=0$, 
we find
\begin{equation}\label{eq:symmetric flux}
\epsilon^{1-2s}(\mathscr{E}_s Q)_\tau(x,\epsilon)
=
\frac12\int_{\R}\epsilon^{1-2s}\partial_\tau P_s(z,\epsilon).
\bigl(Q(x+z)+Q(x-z)-2Q(x)\bigr)\,dz.
\end{equation}
Let
$M_I=\sup_{x\in I_1}|Q''(x)|.$
Taylor's formula implies that for each $x\in I,$
\begin{equation}
   |Q(x+z)+Q(x-z)-2Q(x)|
\leq m_I(z)=\begin{cases}
        M_I|z|^2,
\qquad \quad &|z|\leq1\\
4\|Q\|_{L^\infty(\R)} &|z|>1.
    \end{cases}
\end{equation}
Combining these estimates with
\eqref{eq:symmetric flux}, using the dominated convergence theorem, we have 
\begin{equation}
    \begin{split}
        &\sup_{x\in I}
\left|
\epsilon^{1-2s}(\mathscr{E}_s Q)_\tau(x,\epsilon)
-
sc_s\int_{\R}
\frac{Q(x+z)+Q(x-z)-2Q(x)}
{|z|^{1+2s}}\,dz
\right|\\
\leq&
\frac12\int_{\R}
\left|c_s
\frac{2sz^2-\epsilon^2}
{(z^2+\epsilon^2)^{(3+2s)/2}}-2sc_s|z|^{-1-2s}\right|m_I(z)
\,dz
\to0
    \end{split}
\end{equation}
Recalling~\eqref{fractional laplacian}, we see that $c_{1,s}=2sd_sc_s.$
This proves~\eqref{eq:flux limit}. In addition, it is immediate that $(-\Delta)^sQ$ is continuous.

Combining~\eqref{eq:horizontal trace} and
\eqref{eq:flux limit}, we obtain
\begin{equation}
\label{eq:lower product}
d_s\epsilon^{1-2s}
(\mathscr{E}_s Q)_\tau(\cdot,\epsilon)(\mathscr{E}_s Q)_x(\cdot,\epsilon)
\longrightarrow
-(-\Delta)^sQ\,Q'
\qquad\text{locally uniformly}.
\end{equation}

\medskip

Moreover, fix $0<\epsilon<T<\infty$,
we compute
\[
J_{\epsilon,T}'(x)=
d_s\int_\epsilon^T
\tau^{1-2s}
\bigl((\mathscr{E}_s Q)_x(\mathscr{E}_s Q)_{xx}-(\mathscr{E}_s Q)_\tau (\mathscr{E}_s Q)_{x\tau}\bigr)\,d\tau.
\]
By the extension equation~\eqref{eq:extension energy},
it follows that $\tau^{1-2s}(\mathscr{E}_s Q)_{xx}
=
-\partial_\tau\bigl(\tau^{1-2s}(\mathscr{E}_s Q)_\tau\bigr),$ and hence
\[
J_{\epsilon,T}'(x)
=
-d_sT^{1-2s}(\mathscr{E}_s Q)_\tau(x,T)(\mathscr{E}_s Q)_x(x,T)+
d_s\epsilon^{1-2s}
(\mathscr{E}_s Q)_\tau(x,\epsilon)(\mathscr{E}_s Q)_x(x,\epsilon).
\]
Integrating this identity over a fixed interval $[a,b]$, one has
\begin{equation}\label{eq:integrated truncated}
J_{\epsilon,T}(b)-J_{\epsilon,T}(a)
=
-d_s\int_a^b
T^{1-2s}(\mathscr{E}_s Q)_\tau(x,T)(\mathscr{E}_s Q)_x(x,T)\,dx+
d_s\int_a^b
\epsilon^{1-2s}
(\mathscr{E}_s Q)_\tau(x,\epsilon)(\mathscr{E}_s Q)_x(x,\epsilon)\,dx.
\end{equation}

Furthermore, by Lemma~\ref{lemma Poisson extension}, for $T\geq1$,
\[
T^{1-2s}|(\mathscr{E}_s Q)_\tau(x,T)(\mathscr{E}_s Q)_x(x,T)|
\leq
C_s\|Q\|_{L^\infty(\R)}^2T^{-1-2s}.
\]
 Hence,  for every $x\in [a,b]$
\[|J_{\epsilon,T}'(x)|\leqslant C_s \|Q\|_{L^\infty(\R)}^2+C_{I}\|Q'\|_{L^\infty(\R)}.\]
Here $C_I$ depends on $s,\|Q\|_{L^\infty(\R)}$ and the Lipschitz constant of $Q'$ on $[a-1,b+1]$.

By combining~\eqref{eq:uniform truncation} with~\eqref{eq:lower product}, passing to the limit, we conclude that $J$ is locally Lipschitz continuous, and
\begin{equation}
\label{eq:J integral identity}
J(b)-J(a)
=
-\int_a^b(-\Delta)^sQ(x)Q'(x)\,dx.
\end{equation}
 Since $a$ and $b$ are arbitrary,
the fundamental theorem of calculus yields
\[
J'(x)=-(-\Delta)^sQ(x)Q'(x).
\]

\medskip

We now turn to consider the functional $\mathcal H_Q$. Using the identity for $J'$, we find
\[
\mathcal H_Q'(x)
=
Q'(x)Q''(x)-Q(x)Q'(x)
+Q^{p+1}(x)Q'(x)
-(-\Delta)^sQ(x)Q'(x)=0
\]
Therefore, $\mathcal H_Q$ is constant on $\R$.

\medskip

We claim that $\mathcal H_Q\equiv 0$.

Let us define a nonnegative function
\[
E(x)=Q^2(x)+(Q'(x))^2
+
\int_0^\infty
\tau^{1-2s}\bigl(|(\mathscr{E}_s Q)_x|^2+|(\mathscr{E}_s Q)_\tau|^2\bigr)(x,\tau)\,d\tau.
\]
From the extension energy identity~\eqref{bdsjhkisdfs}
and Tonelli's theorem, it follows that
$\int_{\R}E(x)\,dx<\infty.$
We thus can choose an increasing sequence $R_j\to+\infty$
such that $E(R_j)<\frac1j.$

Hence, up to this sequence, one has $Q(R_j)\to0$, $ Q'(R_j)\to0,$
and $|J(R_j)|\to 0$. 
It follows that $\mathcal H_Q(R_j)\longrightarrow0.$
Since $\mathcal H_Q$ is constant, we conclude that
\begin{equation}\label{eq:zero Hamiltonian}
\mathcal H_Q(x)=0
\qquad\text{for every }x\in\R.
\end{equation}

\medskip

Finally, we prove~\eqref{eq:axis profile}.

Since $Q$ and $P_s(\cdot,\tau)$ are even, their convolution
$(\mathscr E_sQ)(\cdot,\tau)$ is also even. Consequently, $Q'(0)=0$ and $(\mathscr E_sQ)_x(0,\tau)=0$ for every $\tau>0.$ This implies 
\[
0=\mathcal H_Q(0)=
-\frac12Q^2(0)
+\frac{Q^{p+2}(0)}{p+2}
-\frac{d_s}{2}\int_0^\infty
\tau^{1-2s}|(\mathscr{E}_s Q)_\tau(0,\tau)|^2\,d\tau.
\]
As a consequence of this and Lemma~\ref{lemma A_s>0}, we obtain
\begin{equation}
    0<\mathcal A_s(Q)
=
d_s\int_0^\infty
\tau^{1-2s}|(\mathscr{E}_s Q)_\tau(0,\tau)|^2\,d\tau
=
\frac{2Q^{p+2}(0)}{p+2}-Q^2(0).
\end{equation}
We complete the proof.
\end{proof}
\bigskip

\section{An oscillation estimate for $H_*=-\Delta+(-\Delta)^s+V$}\label{sec: An oscillation estimate}

This section serves as a preliminary discussion for Section~\ref{sec:uniqueness}, where we prove the nondegeneracy and uniqueness  results of Theorem~\ref{th:nondegeneracy low s}.  Let  $s\in(0,1)$,  we are concerned with a general class of mixed local/nonlocal Schr\"{o}dinger operators 
	\begin{equation}\label{keydvv}
		H_*=-\Delta+(-\Delta)^s+V,
	\end{equation}
	where $V:\mathbb R\to\mathbb R$ belongs to the Kato class
$K_s$ defined in~\cite[Definition~B.1]{SZZ2026}.

The Kato condition implies that, for every $\varepsilon>0$,
there exists $C_\varepsilon>0$ such that
\[
\int_{\mathbb R}|V||f|^2\,dx
\leq
\varepsilon\left(
\|f'\|_2^2+\|(-\Delta)^{s/2}f\|_2^2
\right)
+C_\varepsilon\|f\|_2^2,
\qquad f\in H^1(\mathbb R).
\]
Hence, the associated quadratic form is closed
and bounded from below, with domain $H^1(\mathbb R)$,
and defines a self-adjoint realization of $H_*$.
Its $L^2$ eigenfunctions are bounded and continuous;
see~\cite[Remark~B.1]{SZZ2026}.

Under the assumptions stated in the respective subsections,
we establish the following nodal properties:
\begin{itemize}
\item
In Section~\ref{sec:The second radial eigenfunction of linear operator},
an even eigenfunction associated with the second eigenvalue
of $H_*$ changes sign at most once on $(0,\infty)$.

\item
In Section~\ref{sec:A unique nodal point of the even kernel function},
for $Q\in\mathcal G_s$, any hypothetical nontrivial even element of $\ker L_{s,Q}$ has exactly one nodal point on $(0,\infty)$ and has opposite strict signs on the two intervals separated by that point.
\end{itemize}

\subsection{The second eigenfunction of linear operator}\label{sec:The second radial eigenfunction of linear operator}
This section is devoted to proving that the second eigenfunction of $ H_*$ changes its sign at most twice on $\R$, as follows.
\begin{Theorem}\label{theorem second eigenvalue } Let $s\in (0,1)$. 
	Assume that  $\lambda_2<\inf\sigma_{\rm ess}(H_*)$ be the second eigenvalue of $H_*$ and let $\phi_2\in H^1(\mathbb{R})$ be  a corresponding real-value eigenfunction. Then $\phi_2 = \phi_2(x) $ changes its sign at most twice on $\R$.
    
    In particular, if $\phi_2 = \phi_2(|x|) $ is an even eigenfunction, then $\phi_2$ changes its sign at exactly once on the positive axis $\{x>0\}$.
\end{Theorem}
    \begin{remark}\label{remark changes sign}
        By Perron-Frobenius arguments (see e.g. \cite{SZZ2026} Appendix B) the first eigenfunction
$\phi_1>0$ of $H_*$ is always strictly positive. Hence, by the self-adjointness
of $H_*$, we easily see that $\phi_2$ changes its sign at least once in order to satisfy the
orthogonality condition.
    \end{remark}

	\subsubsection{Bounds on the number of nodal domain}\label{sec:Bounds on the number of nodal domian}
	 Our objective is to derive oscillation bounds for radial eigenfunctions of operator $ H_*=-\Delta+(-\Delta)^s+V $ under a mild condition on the potential $V$ (i.e. $V$ belongs to a ``Kato class" ). The  core of our approach adapted from~\cite{MR3070568, MR3530361}, lies in studying an associated variational problem formulated in the upper half-space $\mathbb{R}^{2}_+.$  The presentation here follows the lines of~\cite{MR3070568, MR3530361}, which treated the fractional case, and will therefore be kept concise, with details provided only where necessary.

However, a decisive difference to ~\cite{MR3070568, MR3530361} will be that the energy functional is defined on a closed subspace of a classical weighted Sobolev space.	This restriction arises from the nature of the mixed-order operator. Consequently, special attention must be paid during the proof to the effects caused by the term $-\Delta$. 
	\bigskip
    
\noindent\par\textbf{Variational formulation on $\mathbb{R}^{2}_+$.}

	By  following the approach in~\cite{MR3070568, MR3530361},  we leverage the localization technique (see~\eqref{eq:extension energy}) to present a variational  characterization of the eigenvalues of the mixed operator $H_*$ via a local energy functional.

To this end,we first recall the homogeneous Sobolev space $\dot{H}^1(\mathbb{R}^{2}_+; \,t^{1-2s})$. From~\cite[Theorem 18.57]{MR3726909}, we know that there exists a linear bounded trace
	operator Tr:$ \dot{H}^1(\mathbb{R}^{2}_+;\, t^{1-2s})\rightarrow \dot{H}^s(\mathbb{R})$ such that $ \operatorname{Tr}U=U(x,0)$, and  the following inequality holds:
	\begin{equation}\label{keysdvseew}
		\int_{\R}\int_{0}^\infty t^{1-2s}|\nabla U|^2\,dt dx\geqslant \frac{1}{d_s}\int_{\mathbb{R}}|(-\Delta)^{s/2} \operatorname{Tr}U|^2\, dx.
	\end{equation}
In particular,	 equality holds in~\eqref{keysdvseew} if and only if $U$ 
is given as the $\mathscr{E}_s$-extension of its trace. 

  Based on~\eqref{eq:extension energy},  we define $\mathscr{H}^1(\mathbb{R}^{2}_+;\, t^{1-2s})$ as the completion of $C_c^\infty(\overline{\mathbb R^2_+})$ with respect to the norm
\[
\|U\|^2_{\mathscr{H}^1}=d_s\int_{\R}\int_{0}^\infty t^{1-2s}|\nabla U|^2\,dt dx+\int_{\mathbb{R}}\left(| U(x,0)|^2+|\nabla_x U(x,0)|^2\right)\, dx,
\]
where  $U(x,0)$ denotes the trace on $\partial \mathbb{R}^{2}_+$, and $C_c^\infty(\overline{\mathbb R^2_+})$ denotes smooth functions with compact support that may meet the boundary.
	
	 
	Moreover,  we introduce the functional
	\begin{equation}\label{keyvfdv}
		\mathcal I (U):=d_s\int_{\R}\int_{0}^\infty t^{1-2s}|\nabla U|^2dtdx+\int_{\mathbb{R}}V(x)|U(x,0)|^2\, dx+ \int_{\mathbb{R}}|\nabla_x U(x,0)|^2\, dx,
	\end{equation}
defined for $U\in \mathscr{H}^1(\mathbb{R}^{2}_+;\, t^{1-2s})$.
	 
	If the  potential $V\in K_s$ is even, we also define  a corresponding closed subspace.
	\begin{equation}\label{working space1}
		\mathscr{H}_{\rm even}^1(\mathbb{R}^{2}_+;\, t^{1-2s}):=\left\{U\in	\mathscr{H}^1(\mathbb{R}^{2}_+;\, t^{1-2s}):
U(-x,t)=U(x,t)
\text{ for almost every }(x,t)\in\mathbb R^2_+ \right\}.
	\end{equation}
	We are now in a position to characterize the discrete eigenvalues of the operator $H_*=-\Delta+(-\Delta)^s+V$ in terms of the local energy functional $\mathcal{I}(U)$ introduced above. See also~\cite{MR3070568, MR3530361} for a similar
	result for $(-\Delta)^s$.
	
	\begin{Lemma}\label{lemma eigenvalue linear}
		Let  $s\in(0,1)$, and $V \in K_s$. Suppose that $m \geqslant 1$ is an integer and assume that $H_* = -\Delta+(-\Delta)^s + V$ acting on $L^2(\mathbb{R})$ has at least $m$ eigenvalues
		\[
		E_1 \leq E_2 \leq \cdots \leq E_m < \inf \sigma_{\mathrm{ess}}(H_*).
		\]
		Furthermore, let $\mathscr{N}$ be an $(m-1)$-dimensional subspace of $L^2(\mathbb{R})$ spanned by eigenfunctions corresponding to the eigenvalues $E_1, \ldots, E_{m-1}$. Then we have
		\begin{equation}\label{hjjdhs}
		    E_m = \inf\limits_{U\in \mathscr{M}}\mathcal{I}(U), 
		\end{equation}
	where $ \mathcal{I}(U)$  is defined in~\eqref{keyvfdv}, and
		\[ \mathscr{M}:=\left\{ U\in \mathscr{H}^1(\mathbb{R}^{2}_+;\, t^{1-2s}), \, \int_{\mathbb{R}} |U(x, 0)|^2 \, dx = 1, \, U(\cdot, 0) \perp_{2} \mathscr{N} \right\} \]
		 Moreover, the infimum is attained if and only if $U = \mathscr{E}_s f$ with $f \in H^1(\mathbb{R})$, where $\|f\|_{L^2(\mathbb{R})} = 1$ and $f \in \mathscr{N}^{\perp}$ is a linear combination of eigenfunctions of $H_*$ corresponding to the eigenvalue $E_m$.
		
		Moreover, if $V \in K_s$ is even, then the same result holds true with $L^2(\mathbb{R})$ and $\mathscr{H}^1(\mathbb{R}^{2}_+;\, t^{1-2s})$ replaced by $L_{\mathrm{even}}^2(\mathbb{R})$ and $\mathscr{H}_{\rm even}^1(\mathbb{R}^{2}_+;\, t^{1-2s})$, respectively, and the eigenvalues counted in the corresponding subspaces.
	\end{Lemma}
\begin{proof}

We define
\[
 B[f]
=
\int_{\mathbb R}|f'|^2\,dx
+
\int_{\mathbb R}|(-\Delta)^{s/2}f|^2\,dx
+
\int_{\mathbb R}V|f|^2\,dx,
\qquad \forall f\in H^1(\mathbb R).
\]
Let $U\in\mathscr M$ and set $f=\operatorname{Tr}U$.
By the definition of $\mathscr M$, we observe that $\|f\|_2=1$
and $f\perp_2\mathscr N$. Owing to the extension inequality~\eqref{keysdvseew} and the min-max principle, one derives 
\[
\mathcal I(U)\geq B[f]\geq E_m.
\]

We now choose a normalized eigenfunction $f_m\in\ker(H_*-E_m)\cap\mathscr N^\perp$, then $\mathscr E f_m\in \mathscr M$, and $\mathcal I(\mathscr E_sf_m)=B[f_m]=E_m.$
This implies that
\[
\inf_{U\in\mathscr M}\mathcal I(U)
=
\inf_{{
f\in H^1(\mathbb R),\ \|f\|_2=1,
f\perp_2\mathscr N}}B[f]=E_m,
\]
 and hence $\mathscr E f_m$ attains the infimum.

 On the other hand, suppose that $U\in\mathscr M$
attains the infimum and set $f=\operatorname{Tr}U$.
Hence, $E_m\leq B[f]\leq\mathcal I(U)=E_m.$  Recalling the identity~\eqref{bdsjhkisdfs}, we see that 
$U=\mathscr E_sf$, while equality in the spectral variational principle implies $f\in\ker(H_*-E_m)\cap\mathscr N^\perp.$

If $V$ is even, $H_*$ leaves the even subspace invariant, while the Poisson extension preserves symmetry in $x$.
The conclusion follows by restricting the preceding argument to the corresponding even subspaces.
\end{proof}

\begin{remark}
	Since $H_* = -\Delta+(-\Delta)^s + V$ is a real operator, all of its eigenfunctions can be chosen to be real-valued. Furthermore, we note that for $V \in K_s$, every eigenfunction $\psi$ of $H_* $ is continuous. This yields that its extension $\mathscr{E}_s \psi\in C(\overline{\mathbb{R}^{2}_+})$. Hence we can consider its nodal domains, which are defined as the connected components of the open set $\{(x, t) \in \mathbb{R}^{2}_+ : (\mathscr{E}_s \psi)(x, t) \neq 0\}$ in the upper half-space $\mathbb{R}^{2}_+$.
\end{remark}

With Lemma~\ref{lemma eigenvalue linear} at hand, we now derive bounds on the number of nodal domains for the extension of eigenfunctions of $H_*$ to the upper half-space $\mathbb{R}^{2}_+$,  see~\cite{MR3070568, MR3530361} for related results concerning the operator $(-\Delta)^s$.

\begin{Lemma}\label{Lemma number of nodal domains}
		Let  $0 < s < 1$, and $V \in K_s$. Suppose that $m\geqslant 1$ is an integer and assume that $H_* = -\Delta+(-\Delta)^s + V$ acting on $L^2(\mathbb{R})$ has at least $m$ eigenvalues
		\[
		E_1 \leqslant E_2 \leqslant \cdots \leqslant E_m < \inf \sigma_{\mathrm{ess}}(H_*).
		\]
		If $\psi_m \in H^1(\mathbb{R})$ is a real eigenfunction of $H_*$ with eigenvalue $E_m$, then its extension $\mathscr{E}_s \psi_m$ has at most $m$ nodal domains on $\mathbb{R}^{2}_+$.
	\end{Lemma}

To begin with, we introduce some notation.

Let $\psi\in H^1(\R)$ be a real weak eigenfunction, $H_*\psi=E\psi.$
A nodal domain of $W:=\mathscr E_s\psi$ is a connected component of $\{(x,t)\in\R^2_+: W(x,t)\ne0\}.$
For every $U\in\mathscr H^1$, with trace $u$, the weak
extension equation is
\begin{equation}\label{eq:weak extension}
d_s
\int_{\R}\int_{0}^\infty 
t^{1-2s}\nabla W\cdot\nabla U\,dt\,dx+
\int_{\R}
\bigl(\psi'u'+V\psi u\bigr)\,dx
=E\int_{\R}\psi u\,dx.
\end{equation}

\begin{Lemma}
\label{lem:nodal restrictions}
Let $\Omega$ be a nodal domain of $W:=\mathscr E_s\psi$. Then the zero-extended
restriction
$W_\Omega=W\mathbf1_\Omega $
belongs to $\mathscr H^1$. Define the trace by
$\psi_\Omega=\operatorname{Tr}W_\Omega$. Then $\psi_\Omega$
is nontrivial and satisfies
\begin{equation}\label{eq:restriction energy}
\mathcal I(W_\Omega)
=
E\|\psi_\Omega\|_{L^2(\R)}^2.
\end{equation}
In particular, distinct nodal restrictions have disjoint
traces and no mixed terms in the extended energy.
\end{Lemma}

\begin{proof}
 Since $\psi\in H^1(\R)\subset C^{1/2}(\R)$, $\psi$ vanishes at infinity, and
its Poisson extension is continuous up to the boundary. For  $t>0$ large enough, by H\"{o}lder inequality, we see that
\[
\|W(\cdot,t)\|_{L^\infty(\R)}
\leq
\|P_s(\cdot,t)\|_{L^2(\R)}
\|\psi\|_{L^2(\R)}
=
C_s t^{-1/2}\|\psi\|_{L^2(\R)}
\to0 \qquad \text{as } t\to \infty.
\]
For $0<t\leq T$, we observe that
\[
|W(x,t)|
\leq
\sup_{|y|\geq|x|/2}|\psi(y)|+
C_s\|\psi\|_{L^\infty(\R)}
T^{2s}|x|^{-2s}\to 0\qquad \text{as } |x|\to \infty.
\]
This implies that $|W|\to 0$ as $|(x,t)|\to \infty$.

Suppose first that $\Omega$ is a positive nodal domain of $W$ on $\R_+^2$. 
Let $\widetilde{\Omega}$ denote the corresponding connected
component of $\{W>0\}$ in the relative topology of
$\overline{\R^2_+}$, and hence $\widetilde\Omega\cap\mathbb R^2_+=\Omega.$ For $\epsilon>0$, define
\[
W_{\Omega,\epsilon}
=
(W-\epsilon)_+\mathbf1_{\widetilde{\Omega}}.
\]
Since $W$ vanishes at infinity,
$W_{\Omega,\epsilon}$ has compact support in
$\overline{\R^2_+}$.

We claim that $W_{\Omega,\epsilon}\in\mathscr H^1$.

Indeed, $W$ vanishes on the relative boundary of
$\widetilde{\Omega}$. By continuity, each point of this boundary has a relative neighborhood on which $W<\epsilon$. Hence $W_{\Omega,\epsilon}$ is trivial in that neighborhood, this ensures that extending the truncated function by zero outside the component introduces
no jump across the nodal boundary. Therefore, we have the global
weak gradient identity
\[
\nabla W_{\Omega,\epsilon}
=
\mathbf1_{\Omega\cap\{W>\epsilon\}}\nabla W
\qquad\text{a.e. in }\R^2_+.
\]
In particular,
\begin{equation}\label{gsdagsdg}
   \int_{\R}\int_{0}^\infty
t^{1-2s}|\nabla W_{\Omega,\epsilon}|^2\,dt\,dx
\leq
\int_{\R}\int_{0}^\infty
t^{1-2s}|\nabla W|^2\,dt\,dx.
\end{equation}

We next identify the boundary trace of $W_{\Omega,\epsilon}$
and show that it belongs to $H^1(\mathbb R)$. Set $I_\Omega
=\{x\in\mathbb R:(x,0)\in\widetilde\Omega\}.$
Since $W(x,0)=\psi(x)$, each interval component of
$\{\psi>0\}$ lies in a single relative positive component
of $W$. Thus $I_\Omega$ is a union of disjoint open
intervals $J_k$, each of which is a connected component
of $\{\psi>0\}$. By continuity, $\psi$ vanishes at every
finite endpoint of these intervals.

The boundary trace of $W_{\Omega,\epsilon}$ is therefore $f_{\Omega,\epsilon}
=(\psi-\epsilon)_+\mathbf1_{I_\Omega}.$
We claim that $f_{\Omega,\epsilon}\in H^1(\R)$. 

Indeed,  define $g_k=(\psi-\epsilon)_+\mathbf1_{J_k}. $  
We observe that $g_k$ vanishes at every finite
endpoint of $J_k$, its extension by zero belongs to $H^1(\mathbb R)$, with $g_k'
=\mathbf1_{J_k\cap\{\psi>\epsilon\}}\psi' $ a.e. in $\mathbb R.$
Since the intervals are disjoint, one deduces that 
\[
\sum_k\|g_k\|_{H^1(\mathbb R)}^2=
\int_{I_\Omega}(\psi-\epsilon)_+^2\,dx
+
\int_{I_\Omega\cap\{\psi>\epsilon\}}|\psi'|^2\,dx\leq \|\psi\|_{H^1(\mathbb R)}^2.
\]
From this bound and the disjoint supports, it follows that $\sum_k g_k$ converges to $f_{\Omega,\epsilon}$ in
$H^1(\mathbb R)$, and hence
\[
f_{\Omega,\epsilon}\in H^1(\mathbb R),
\qquad
f_{\Omega,\epsilon}'
=
\mathbf1_{I_\Omega\cap\{\psi>\epsilon\}}\psi'
\quad\text{a.e. in }\mathbb R.
\]

Moreover, applying the  dominated convergence theorem, we obtain 
\[
\nabla W_{\Omega,\epsilon}
\longrightarrow
\mathbf1_\Omega\nabla W
\quad\text{in }
L^2(\R^2_+;t^{1-2s}),\quad \text{and}\quad f_{\Omega,\epsilon}
\longrightarrow
\psi\mathbf1_{I_\Omega}
\quad\text{in }H^1(\R)\quad \text{as } \epsilon\to0.
\]
Hence 
$W_{\Omega,\epsilon}\to W_\Omega$  in $\mathscr H^1$. In particular,
\[
W_\Omega\in\mathscr H^1,\qquad
\nabla W_\Omega=\mathbf1_\Omega\nabla W,
\]
and its trace satisfies
\[
\psi_\Omega=\psi\mathbf1_{I_\Omega},
\qquad
\psi_\Omega'=\psi'\mathbf1_{I_\Omega}
\quad\text{a.e. in } \R.
\]
For a negative nodal domain, the same conclusions follow
by applying the argument to $-W$ and then restoring the sign.

We may now use $W_\Omega$ as a test function in
\eqref{eq:weak extension}. The restriction identities
above yield
\[
d_s\int_{\R}\int_{0}^\infty
t^{1-2s}|\nabla W_\Omega|^2\,dt\,dx
+
\int_{\R}
\bigl(|\psi_\Omega'|^2+V\psi_\Omega^2\bigr)\,dx
=
E\|\psi_\Omega\|_{L^2(\R)}^2,
\]
which is \eqref{eq:restriction energy}.
This identity also guarantees that $\psi_\Omega$ is nontrivial.

\medskip

Finally, we show that restrictions associated with distinct nodal
domains have no cross terms. The estimates for the weighted gradient integral above are supported in disjoint domains, while their traces and
trace derivatives are supported in disjoint selected
boundary intervals. Therefore, if $\Omega_i\ne\Omega_j$,
\[
d_s\int_{\R}\int_{0}^\infty
t^{1-2s}
\nabla W_{\Omega_i}\cdot\nabla W_{\Omega_j}\,dt\,dx+
\int_{\R}
\left(
\psi_{\Omega_i}'\psi_{\Omega_j}'
+
V\psi_{\Omega_i}\psi_{\Omega_j}
\right)\,dx
=0.
\]
This ends the proof.
\end{proof}

\begin{proof}[Proof of Lemma~\ref{Lemma number of nodal domains}]
Set $W=\mathscr E_s\psi_m$.
Suppose, argue by contradiction, that $W$ has at least $m+1$ nodal domains. Choose
$\Omega_1,\ldots,\Omega_{m+1}$.
Define $W_j=W\mathbf1_{\Omega_j}
$ and $
\psi_j=\operatorname{Tr}W_j.$
From Lemma~\ref{lem:nodal restrictions},  it follows that these functions are
admissible, their traces are nontrivial and disjoint, and
\[
\mathcal I(W_j)
=
E_m\|\psi_j\|_{L^2(\R)}^2.
\]

Let $\mathscr N=
\operatorname{span}
\{\phi_1,\ldots,\phi_{m-1}\}$
be as in Lemma~\ref{lemma eigenvalue linear}. Since these traces are nontrivial and disjoint, 
there exist real coefficients $c_1,\ldots,c_m$, with $c_i\neq0$ for some $i$,
such that
\[
z=\sum_{j=1}^m c_j\psi_j
\perp_2\mathscr N,\qquad \|z\|_{L^2(\R)}=1.
\]
Set  $Z=\sum_{j=1}^m c_jW_j.$
Recall that all terms involving two different nodal restrictions
are zero, we obtain
\begin{equation}\label{eq:trial-energy}
\mathcal I(Z)=
\sum_{j=1}^m c_j^2\mathcal I(W_j)
=E_m\sum_{j=1}^m
c_j^2\|\psi_j\|_{L^2(\R)}^2
=E_m\|z\|_{L^2(\R)}^2
=E_m.
\end{equation}
Since $z\perp_2\mathscr N$,  by combining Lemma~\ref{lemma eigenvalue linear} with~\eqref{keysdvseew}, we derive
\[
E_m\leq \int_{\mathbb{R}}|(-\Delta)^{s/2}z|^2+\int_{\mathbb{R}}V(x)|z|^2\, dx+\int_{\mathbb{R}}|\nabla z|^2\, dx\leq \mathcal I(Z).
\]
From this and 
\eqref{eq:trial-energy}, we obtain $Z=\mathscr E_s z.$
Thus $Z$ is the Poisson extension of its trace and satisfies
\begin{equation}\label{sdgkisdjmgbo}
    Z_{xx}+Z_{tt}+\frac{1-2s}{t}Z_t=0
\qquad\text{in }\R^2_+.
\end{equation}

On the other hand, $Z$ was constructed using only the
restrictions of $W$ to $\Omega_1,\ldots,\Omega_m$. Hence $Z=0$ in $\Omega_{m+1}.$ For $t>0$, the equation~\eqref{sdgkisdjmgbo} is locally uniformly elliptic
and its coefficients are real analytic. Interior analytic
regularity arguments guarantee that $Z$ is real analytic
in $\R^2_+$. Since the upper half-plane is connected
and $\Omega_{m+1}$ is a nonempty open subset, the identity
principle for real analytic functions implies
\[
Z\equiv0\qquad\text{in }\R^2_+.
\]
Applying the Sobolev trace theorem, we conclude
\[
z=\operatorname{Tr}Z=\operatorname{Tr}0=0.
\]
This contradicts the fact that
$\|z\|_{L^2(\R)}=1$. We thus complete the proof.
\end{proof}

\subsubsection{The second eignfunction of $H_*$}\label{sec:The proof of theorem second eigenvalue}

  \begin{proof}[The proof of Theorem~\ref{theorem second eigenvalue }]
Let $\Phi_2=\mathscr E_s\phi_2$
be the $s$-Poisson extension of $\phi_2$. Since $\phi_2\in H^1(\R)$, $\phi_2$ is continuous and bounded in $\R$. Hence,
$\Phi_2$ is continuous on $\overline{\R^2_+}$ and
satisfies $\Phi_2(x,0)=\phi_2(x).$

 By Remark~\ref{remark changes sign},    $\phi_2$ takes both positive and negative values,
continuity up to the boundary ensures that $\Phi_2$ is
positive at some points in $\R_+^2$ and negative at others. Thus $\Phi_2$ has at least one positive
nodal domain and one negative nodal domain.

On the other hand, in view of
Lemma~\ref{Lemma number of nodal domains}, applied with
$m=2$,  we see that $\Phi_2$ has exactly two  nodal domains,  which we denote by
\[
\{\Phi_2>0\}=\Omega_+,
\qquad
\{\Phi_2<0\}=\Omega_-.
\]
Suppose, argue by contradiction, that $\phi_2$
changes sign at least three times.
Then there  exist $x_1<x_2<x_3<x_4$
such that, after replacing $\phi_2$ by $-\phi_2$
if necessary,
\[
\phi_2(x_1)>0,\qquad
\phi_2(x_2)<0,\qquad
\phi_2(x_3)>0,\qquad
\phi_2(x_4)<0.
\]
By  continuity of $\Phi_2$ again, we choose  $\epsilon_0>0$
such that
\[
(x_i,t)\in\Omega_+
\quad\text{for }i=1,3,\quad \text{and}
\quad (x_i,t)\in\Omega_-
\quad\text{for }i=2,4,\]
for every $0<t\leq\epsilon_0$. Thus, there exist short vertical segments above
$(x_1,0)$ and $(x_3,0)$ on $\Omega_+$.

Note that    the
connected open set $\Omega_+\subset \R_+^2$  must be arcwise connected. 
Thus we conclude that
there exists a simple continuous curve 
$\gamma_+:[0,1]\to
\overline{\R^2_+}$
such that
\[
\gamma_+(0)=(x_1,0),\qquad
\gamma_+(1)=(x_3,0),\quad \text{and}\quad \gamma_+((0,1))\subset\Omega_+.
\]
The union $J=
\gamma_+([0,1])
\cup
\bigl([x_1,x_3]\times\{0\}\bigr)$
is a simple closed curve.
By the Jordan curve theorem, its complement has a
bounded component $D$ and an unbounded component,
with $\partial D=J.$
Since $J$ lies in $\overline {\R^2_+}$, its
bounded component $D$ lies in the open upper half-plane.

It is immediate that for sufficiently small $\epsilon>0$, we have
\[
(x_2,\epsilon)\in D,
\qquad
(x_4,\epsilon)\notin\overline D, \qquad (x_2,\epsilon)\in\Omega_-,\qquad (x_4,\epsilon)\in\Omega_-.
\]
Since $\Omega_-$ is open and connected, these two points
can be joined by a continuous path lying entirely
in $\Omega_-$. This path starts inside $D$ and ends outside
$\overline D$, so it must  intersect the arc $\gamma_+$.  This is impossible, since $\Phi_2$ is negative along the path
and positive along the interior of $\gamma_+$.

Thus four boundary points with alternating signs
cannot exist, and $\phi_2$ changes sign at most twice
on $\R$.

Finally, suppose that $\phi_2$ is even.  then clearly $\phi_2$ can change its sign on $\{x>0\}$
at most once, since otherwise $\phi_2$ would change  sign at least four times on $\R$,
contradicting the result just derived. By Remark~\ref{remark changes sign}, we deduce that $\phi_2$ changes sign exactly once on
$(0,\infty)$.
\end{proof}

\subsection{A unique nodal point of the even kernel function on $\R_+$}\label{sec:A unique nodal point of the even kernel function}

\begin{Lemma}\label{lemma boundary zero}
Let $s\in (0,1)$ and $V:\R\to\R$ be bounded and continuous.
Suppose that  $ h\in H^1(\R)\cap L^1(\R)$ is  real-valued 
satisfying
\begin{equation}\label{fdbhgdfshgf}
    h\neq 0\qquad \text{and}\quad 
- h''+ (-\Delta)^sh+h -Vh=0
\qquad\text{in }\R.
\end{equation}
Let $W=\mathscr E_s h$ be the $s$-Possion extension. 
If $h(r)=0$, then every relative neighborhood of $(r,0)$
in $\overline{\R_2^+}$ contains both a strictly
positive and a strictly negative point of $W$ in $\R_2^+$.

Moreover, $h$ cannot vanish identically on any nonempty
open interval.
\end{Lemma}

\begin{proof}
To this end, we first establish the boundary regularity needed in the
argument, and then construct an explicit barrier near a boundary zero.

We recall that, from~\cite[Proposition A.1]{SZZ2026} and the regularity arguments developed in~\cite{DSVZ24},  we know that $h\in H^2(\R)\cap C^2({\R})$. This implies that $h$ is locally  uniformly continuous, hence
 $W$ is continuous up to the boundary and
$W(x,0)=h(x).$

Moreover, from~\eqref{eq:extension energy}, we know that
\begin{equation}\label{sdfvs}
    \begin{cases}
    W_{xx}+W_{tt}+\frac{1-2s}{t}W_t=0\qquad &\text{in } \R_+^2\\
    -d_s\lim_{t\to0}t^{1-2s}W_t(x,t)
= (-\Delta)^s h(x) &\text{in } \R.
\end{cases}
\end{equation}
We observe that, this equation is uniformly
elliptic and has real analytic coefficients near every point with $t>0$.
Interior analytic regularity therefore implies that
$W$ is real analytic in $\R_+^2$.

In particular, if $W$ vanished on a nonempty interior
open set, the identity principle would imply
\[
W\equiv0
\qquad\text{in }\R^2_+.
\]
Furthermore, its boundary trace would then be trivial, which is in contradiction with
$h\ne0$.

\medskip

Suppose that $h(r)=0$.

\textbf{First of all, assume that  $W\ge0$ in a relative neighborhood $\mho$ of $(r,0)$.}

Choose $\rho>0$ and $T>0$ sufficiently small that $[r-2\rho,r+2\rho]\times[0,2T]\subset \mho$. On the connected open rectangle $(r-2\rho,r+2\rho)\times(0,2T),$ since  $\operatorname{Tr}W=h\neq0$,
the strong maximum principle guarantees that \[
W(x,t)>0 \qquad \text{if }|x-r|<2\rho,\quad 0<t<2T.
\]

We next construct a barrier $B$ in the smaller rectangle $\mathcal R=(r-\rho,r+\rho)\times(0,T).$
By comparing $W$ with a suitable positive multiple of $B$, we shall obtain a lower bound for $W$.

Set $k=\frac{\pi}{2\rho}$
and define
\begin{equation}\label{ikjgfrdgd}
    d_0=1,
\qquad
d_n
=\frac{(k^2/4)^n}
{n!\prod_{j=1}^n(s+j)}
\quad\text{for }n\ge1.
\end{equation}
Let
\[
S(t)=t^{2s}\sum_{n=0}^{\infty}d_n t^{2n}\qquad \text{and}\qquad B(x,t)=\cos(k(x-r))S(t).
\]
Since
$\frac{d_{n+1}t^{2n+2}}{d_n t^{2n}}
=
\frac{k^2t^2}
{4(n+1)(n+1+s)}\to 0 $ uniformly when $t$ remains bounded.
Thus the series converges uniformly on every bounded interval
of $t\ge0$. It can also be differentiated term by term for
$t>0$, with the differentiated series converging uniformly
on compact subsets of $(0,\infty)$.

We claim that \begin{equation}\label{nfsdfsdf}
    B_{xx}+B_{tt}+\frac{1-2s}{t}B_t=0
\qquad\text{in }\mathcal R.
\end{equation}
A direct calculation yields
\[
S''(t)+\frac{1-2s}{t}S'(t)=
\sum_{n=1}^{\infty}
4n(n+s)d_n t^{2s+2n-2}=
\sum_{n=0}^{\infty}
4(n+1)(n+1+s)d_{n+1}t^{2s+2n}=
k^2S(t).
\]
On the other hand, we have $B_{xx}=-k^2B.$ This proves~\eqref{nfsdfsdf}. 

Since  $d_n>0$ and $\cos(k(x-r))>0$ for every $|x-r|<\rho$, we see that
$B(x,t)>0 $ in $
\mathcal R.$ Moreover, we also have
\begin{equation}\label{fgdfghbfd}
    B(x,0)=0,
\qquad
B(r-\rho,t)=B(r+\rho,t)=0.
\end{equation}
In particular,  when $x=r$, one has 
\begin{equation}\label{fbgfdsbh}
    B(r,t)=t^{2s}\bigl(1+O(t^2)\bigr)
\qquad\text{as } \, t\to 0.
\end{equation}
\medskip

We claim that there exists $\eta>0$, such that  $W\ge\eta B $ in $ \mathcal R.$

Indeed, $W(x,T)>0$ for every $x\in [r-\rho,r+\rho]$.
By the continuity of $W$ and $B$, we see that
\[
m=\min_{|x-r|\le\rho}W(x,T)>0,\quad \text{and}
\quad M=\max_{|x-r|\le\rho}B(x,T)>0.\]
Choose $\eta=\frac{m}{2M}>0.$
Then
\[
W(x,T)-\eta B(x,T)\ge\frac m2>0.
\]
 On the vertical sides and the bottom, since $W\ge0$, recalling~\eqref{fgdfghbfd},  we see that $W-\eta B\ge0.$

Let $Z=W-\eta B.$ 
For $0<\epsilon<T$, consider
$\mathcal R_\epsilon
=
(r-\rho,r+\rho)\times(\epsilon,T)$
and set
$m_\epsilon
=
\eta\max_{|x-r|\le\rho}B(x,\epsilon).$
Moreover, for every $x\in [r-\rho,r+\rho]$
\[
Z(x,\epsilon)+m_\epsilon
=
W(x,\epsilon)-\eta B(x,\epsilon)+m_\epsilon
\ge0.
\]
Thus $Z+m_\epsilon\ge 0$ on $\partial \mathcal R_\epsilon.$

By combining~\eqref{sdfvs} with~\eqref{nfsdfsdf}, $Z+m_\epsilon$ satisfies
\begin{equation}\label{dfvdag}
    \operatorname{div}
\bigl(t^{1-2s}\nabla(Z+m_\epsilon)\bigr)=0\qquad \text{in }\mathcal R_\epsilon.
\end{equation}
Since $\overline{\mathcal R_\epsilon}\subset\R_+^2$,
both $W$ and $B$ are smooth near this rectangle, and hence
$Z+m_\epsilon\in H^1(\mathcal R_\epsilon)$.
Moreover, $\epsilon\le t\le T$ gives positive upper and
lower bounds for the weight $t^{1-2s}$. Its negative part has zero boundary trace and is an
admissible test function. By testing~\eqref{dfvdag} against $(Z+m_\epsilon)_-$, one has
\[
0 =-\int_{\mathcal R_\epsilon}
t^{1-2s}
\left|
\nabla(Z+m_\epsilon)_-
\right|^2\,dx\,dt.
\]
Thus $\nabla (Z+m_\epsilon)_-=0$ a.e. in $\mathcal R_\epsilon$. Its zero boundary trace guarantees that $(Z+m_\epsilon)_-=0$ in $\mathcal R_\epsilon$.
Therefore $Z\ge-m_\epsilon $ in $
\mathcal R_\epsilon.$
The barrier expansion~\eqref{ikjgfrdgd} implies that $m_\epsilon
\le C\epsilon^{2s}
\to 0$ as $\epsilon\to 0$.
For any fixed point of $\mathcal R$, we may now let
$\epsilon\to 0$ and conclude that
\begin{equation}\label{cxhgbdfhg}
    W\ge\eta B
\qquad\text{in }\mathcal R.
\end{equation}

\medskip
Finally, we complete the proof of the first point.

Recall that $W\ge0$ near $(r,0)$, thus, $h(x)\ge0 $ for $
x\text{ near }r.$
Since $h(r)=0$ and $h\in C^2$, we see that
$r$ is a local minimum and $h''(r)\ge0.$
Moreover, 
owing to~\eqref{sdfvs}, we find that
\begin{equation}\label{sdfgfdgd}
    \ell
=\lim_{t\to0}t^{1-2s}W_t(r,t)=-d_s^{-1}(-\Delta)^s h(r)=d_s^{-1}(-h''(r)+(1-V)h(r))=-d_s^{-1}h''(r)\leq 0.
\end{equation}

However, by the comparison~\eqref{cxhgbdfhg}, it follows   that $\ell>0$.  Indeed, since
 $t^{1-2s}W_t(r,t)\to\ell$, one has $t^{1-2s}W_t(r,\cdot)$ is bounded near $t=0$.
We observe that, for $0<\epsilon<t$,
\[
W(r,t)-W(r,\epsilon)
=
\int_\epsilon^t
W_t(r,\cdot)\,d\tau=\int_\epsilon^t\tau^{2s-1}
[\tau^{1-2s}W_t(r,\cdot)]\,d\tau.
\]
Since  $\int_0^t\tau^{2s-1}\,d\tau
<\infty$, the dominated convergence theorem
implies that 
\[
W(r,t)
=\int_0^t\tau^{2s-1}[\tau^{1-2s}W_t(r,\cdot)]\,d\tau.
\]
where we use the fact that $W(r,0)=h(r)=0$.
Moreover, we find that
\[
\left|
W(r,t)-\frac{\ell}{2s}t^{2s}
\right|\le
\sup_{0<\tau<t}|\tau^{1-2s}W_t(r,\cdot)-\ell|
\int_0^t\tau^{2s-1}\,d\tau=
\frac{t^{2s}}{2s}
\sup_{0<\tau<t}|\tau^{1-2s}W_t(r,\cdot)-\ell|=
o(t^{2s}).
\]
As a consequence of this, 
\[
W(r,t)
=\frac{\ell}{2s}t^{2s}+o(t^{2s}).
\]
Combining~\eqref{fbgfdsbh} with~\eqref{cxhgbdfhg}, we obtain that
\[
W(r,t)
\ge
\eta t^{2s}\bigl(1+O(t^2)\bigr).
\]
Consequently, $\ell\ge2s\eta>0.$ This is in contradict with~\eqref{sdfgfdgd}.

\medskip

We have already ruled out a nonnegative relative neighborhood.
The same argument excludes a nonpositive relative neighborhood.

Consequently, $W$ takes both positive and negative values
in every relative neighborhood of a boundary zero.
If either sign occurs at a boundary point, continuity
ensures that the same sign also occurs at nearby interior
points. Hence there are interior points of both signs
arbitrarily close to the boundary zero.

\medskip

\textbf{We now turn to show that the trace cannot vanish on an open interval.}

Suppose that $h=0$ on a nonempty open interval $I$.
Then $h''=0$ on $I$, and the equation~\eqref{fdbhgdfshgf} implies
$(-\Delta)^sh=0$ on $I$.
By exterior unique continuation property of the fractional Laplacian~\cite[Theorem~1.2]{Ghosh}, we conclude that
$h\equiv0$, contradicting the assumption $h\ne0$.

\end{proof}

\begin{Lemma}\label{lemma unique zero}
Let $s\in (0,1)$ and $p>0$. Let $Q_s\in \mathcal{G}_s$.
Suppose that $h\in H^1(\R)\cap L^1_{\rm even}(\R)$ with $h\neq 0$, satisfies
\[
L_{s,Q_s}h=0.
\]
Then, after replacing $h$ by $-h$ if necessary, there
exists a unique $r>0$ 
\[
h(r)=0,\qquad
h(x)>0\quad\text{for }|x|<r,\qquad
h(x)<0\quad\text{for }|x|>r.
\]
With this choice of sign,
\begin{equation}\label{nod positive moments}
\int_{\R}h(x)\,dx<0,
\qquad
\int_{\R}Q_s(x)h(x)\,dx<0,
\qquad
\int_{\R}Q_s^{p+1}(x)h(x)\,dx=0.
\end{equation}
\end{Lemma}

\begin{proof}
Recall that $Q_s$ solves equation~\eqref{main equation}, we know that $L_{s,Q_s}Q_s=-pQ_s^{p+1}$.  Since $L_{s,Q_s} h=0$, together with the fact that  $Q_s,h\in H^2(\R)$,  integration by parts for the local term
and symmetry of the fractional term yield
\begin{equation}\label{njdnfkdms}
    \int_{\R}Q_s^{p+1}h\,dx=0.
\end{equation}
Note that $Q_s^{p+1}>0$ is even, and  $h\neq 0$ is continuous and even, the identity~\eqref{njdnfkdms}
forces $h$ to take both signs on $\R_+$.

By assumption, $L_{s,Q_s}$ has exactly one strictly negative
eigenvalue. Using the previously established fact that
$0<\inf\sigma_{\mathrm{ess}}(L_{s,Q_s})=1$,  it follows  that $0$ is the indeed second eigenvalue, 
and $h/\|h\|_2$ can be chosen as the second eigenfunction.

We notive that the potential $V=1-(p+1)Q_s^p$ is bounded and even,
and belongs to $K_s$.
Applying Theorem~\ref{theorem second eigenvalue } to $h$,
and replacing $h$ by $-h$ if necessary, we obtain
$r>0$ with $h(r)=0$ such that
\begin{equation}\label{fbvkdfbk}
h(x)\ge0\quad\text{for }0<x<r,
\qquad
h(x)\le0\quad\text{for }x>r.
\end{equation}
Let $W=\mathscr E_s h$ be the $s$-Poisson extension of $h$.
Lemma~\ref{Lemma number of nodal domains} implies that $W$ has exactly one positive nodal domain
and exactly one negative nodal domain.
In particular, both sets
$\{W>0\} $ and $
\{W<0\}$
are connected. 

Moreover, Lemma~\ref{lemma boundary zero} applies to $h$
with the bounded continuous potential $V=1-(p+1)Q_s^p,$ we know that  
$h$ cannot vanish on a nonempty open interval,
and $W$ cannot have a one-sign relative neighborhood
at a boundary zero of $h$.
\medskip

\textbf{We claim that $r$ is the unique root of $h$ on $\R_+$.}

Suppose, for a contradiction, that $h(\zeta)=0 $ for some $
\zeta\notin\{-r,r\}.$

Choose an open interval $I_\zeta$ containing $\zeta$
and disjoint from $\{-r,r\}$.
By \eqref{fbvkdfbk}, there exists
$\sigma\in\{-1,1\}$ such that $\sigma h\ge0$ on $ I_\zeta.$

Set $U=\sigma W.$
Thus $U(x,0)=\sigma h(x)\ge0$ for every $x\in I_\zeta.$
Since $h$ has no nonempty open interval of zeros,
we can choose $a<\zeta<b $ such that $
[a,b]\subset I_\zeta,$ and
\[
\sigma h(a)>0,
\qquad
\sigma h(b)>0.
\]
Since $U(a,0)>0$ and $U(b,0)>0$,  continuity guarantees short vertical segments above these points on which
$U>0$. Their upper endpoints lie in the connected open set $\{U>0\}$ and can therefore be joined by
a polygonal path contained in this set.
Attaching the two vertical segments  and removing loops, we obtain a simple arc $\Gamma$ with no
self-intersections, joining $(a,0)$ to $(b,0)$ and satisfying
\[
U>0
\qquad\text{on }\Gamma\setminus\{(a,0),(b,0)\}.
\]
The arc and the horizontal segment form a Jordan curve $\mathcal C=
\Gamma\cup\bigl([a,b]\times\{0\}\bigr).$
Let $D$ be its bounded complementary component. Then $\partial D=\mathcal C,
$ with $
D\subset\R^2_+.$
By construction, we observe that $U\ge0$ on $\partial D.$ Since $D$ is connect, combining continuity of $U$ with  strong maximum principle, we derive that $U\ge0 $ throughout $D.$
\medskip



Finally, choose $\epsilon>0$ sufficiently small
that the ball centered at $(\zeta,0)$ with radius
$\epsilon$ avoids $\Gamma$ and the endpoints
$(a,0),(b,0)$.
Thus the Jordan curve $\mathcal C$ consists only
of the horizontal segment inside this ball.

The lower half-ball belongs to the unbounded
complementary component.
The upper half-ball is connected and disjoint from the Jordan curve.
Since $(\zeta,0)\in\partial D$ and $D$ lies above the real axis, the upper half-ball contains points
of $D$ and therefore lies entirely in $D$.

We have  already proved that $U\ge0$ in this upper half-ball.
 We also have $U(x,0)=\sigma h(x)\ge0$
on its bottom segment.
Thus $U\ge0$ in a relative neighborhood of
$(\zeta,0)$, while $U(\zeta,0)=\sigma h(\zeta)=0.$
This yields that $W=\sigma^{-1}U$ has one sign in a relative
neighborhood of $(\zeta,0)$.
This contradicts  Lemma~\ref{lemma boundary zero}.

Consequently, $h$ has no root other than $\pm r$.
Combining this fact with the weak inequalities
in \eqref{fbvkdfbk}, we obtain
\begin{equation}\label{rgregrdfd}
    h(x)>0\quad\text{for }|x|<r,
\qquad
h(x)<0\quad\text{for }|x|>r.
\end{equation}
In particular,
\[
h(0)>0,
\]
and $r$ is the unique nodal point on $(0,\infty)$.
\bigskip

We now show the estimates in~\eqref{nod positive moments}.

For $\theta\in\{0,1\}$,  together with the facts that $Q_s(r)$ is positive and strict decreasing in $|x|$, using~\eqref{njdnfkdms} and~\eqref{rgregrdfd}, it follows that
\[
\int_{\R}Q_s^\theta h\,dx=
\int_{\R}
Q_s^\theta
\left[
1-
\left(\frac{Q_s}{Q_s(r)}\right)^{p+1-\theta}
\right]h\,dx<0.
\]
Taking $\theta=0$ and $\theta=1$, we conclude that~\eqref{nod positive moments}. We thus prove all
the asserted identities and inequalities.
\end{proof}

\section{Uniqueness for positive solutions with Morse index equal to one}\label{sec:uniqueness}

In this section, we establish uniqueness and nondegeneracy
for solutions in $\mathcal G_s$ for every $s\in(0,1)$.

We use Lemma~\ref{lemma unique zero} to prove monotonicity
of the quotient of a hypothetical nontrivial even kernel
element by $-Q_s'$ on the positive half-line.
This property determines the sign of a pairing involving
the local part of the operator, while a dilation identity
forces the opposite sign. The resulting contradiction
excludes nontrivial even kernel elements.
On the odd subspace, positive semi-definiteness of the
quadratic form and simplicity of the lowest eigenvalue identify the kernel as $\operatorname{span}\{Q_s'\}$.
Together, these conclusions establish nondegeneracy.

We then prove uniqueness by continuation in $s$, following the strategy of~\cite{MR3070568,MR3530361}. Nondegeneracy ensures that solution branches are locally unique, while compactness of $\mathcal G_s$ on closed subintervals of $(0,1]$
allows us to continue and count all solutions in this class.

\subsection{Nondegeneracy}
This section is devoted  to establish the nondegeneracy result of solutions
to~\eqref{main equation} whose linearized operators have Morse index equal to one.
	The proof will be divided  into two main steps: nondegeneracy both in the odd sector and even sector.

\subsubsection{Nondegeneracy in the even sector} 
 We now show the triviality of the kernel of linearized operator in the space of even function.

\begin{Lemma}
\label{th:nondegeneracy s<1/2}
Let $0<s<1$ and $p>0$. If $Q_s\in\mathcal G_s$, then
\[
\ker L_{s,Q_s}\cap L^2_{even}(\R)
=\{0\}.
\]
\end{Lemma}

To begin with, we introduce some preliminary results.

\begin{Lemma}\label{lemma: u tail at infinity}
Let $s\in (0,1)$ and $Q_s\in \mathcal M_s$. Suppose that  $u\in H^1(\R)$ solves
\begin{equation}\label{eq:linearized equation}
-u''+(-\Delta)^s u+u
=(p+1)Q_s^p u
\qquad\text{in }\R.
\end{equation}
Then $u\in L^1(\R)$, and
\begin{equation}\label{eq:u asymptotic}
u(x)
=c_{1,s}
\left(\int_{\R}u(y)\,dy\right)
|x|^{-1-2s}
+
o\bigl(|x|^{-1-2s}\bigr)
\qquad\text{as }|x|\to\infty,
\end{equation}
where $c_{1,s}$ is defined in~\eqref{c_{1.s}}.
In particular, if $M:=\int_{\R}u(y)\,dy\ne0,$
there exist $R_0>1$ and $c_u\in(0,1)$ such that
\begin{equation}\label{eq:u lower}
|u(x)|
\geq c_u|x|^{-1-2s}
\qquad\text{for every }|x|\geq R_0.
\end{equation}
\end{Lemma}

\begin{proof}
Recall that  $Q_s\in C^{2,\alpha}(\R)\cap L^\infty(\R)$ and $Q_s(x)\leqslant C|x|^{-1-2s}$ for all $|x|>1$, by~\cite[Proposition A.1]{SZZ2026},  employing the regularity arguments developed in~\cite{DSVZ24},  we know that $u\in H^2(\R)\cap C^2({\R})$. 

   Moreover, in the light of ~\cite[Lemma C.1]{SZZ2026},  there exists $C>0$ such that  \[|u(x)| \leq {C_{s,u}}{|x|^{-1-2s}} \qquad \forall |x|>1.\]  Hence $u\in L^1(\R) $ and $  (p+1)Q_s^p u\in L^2(\R)\cap L^\infty(\R)$.
   
\smallskip
We claim that \begin{equation}\label{bvsjdjnk}
\lim_{|x|\to\infty}|x|^{1+2s} u(x)
=
c_{1,s}\int_{\R}u(y)\,dy.
\end{equation}
Indeed, we observe that $u=G_s\ast (p+1)Q_s^pu, $ since $u\in L^1(\R)$,
Fubini's theorem yields that
\begin{equation}\label{eq:mass identity}
\int_{\R}u(x)\,dx=
(p+1)\left(\int_{\R}G_s\right)
\left(\int_{\R}Q_s^pu\right)= (p+1)
\int_{\R}Q_s^p(y)u(y)\,dy.
\end{equation}
For $|x|\geq1$, we have
\[
u(x)=(p+1)\int_{|y|\leq |x|/2}
G_s(x-y)Q_s^p(y)u(y)\,dy+(p+1)\int_{|y|>|x|/2}
G_s(x-y)Q_s^p(y)u(y)\,dy=:I_1(x)+I_2(x).
\]
For every fixed $y$, we find that $\frac{|x-y|}{|x|}\to1$ as $|x|\to\infty$.
 From Proposition~\ref{prop:green tail}, we know that $|x|^{1+2s} G_s(x-y)
\to c_{1,s}.$
As for $|y|\leq |x|/2$, one also has
\[
|x|^{1+2s} G_s(x-y)
\leq
C_G\frac{|x|^{1+2s}}{(1+|x|/2)^{1+2s}}
\leq2^{1+2s} C_G.
\]
Applying dominated convergence theorem to $(p+1)|x|^{1+2s} G_s(x-y)Q_s^p(y)u(y)
\mathbf1_{\{|y|\leq |x|/2\}}$, we obtain
\begin{equation}\label{eq:inner limit}
|x|^{1+2s} I_1(x)=  (p+1)|x|^{1+2s}\int_{|y|\leq |x|/2}
G_s(x-y)Q_s^p(y)u(y)\,dy
\to
c_{1,s}(p+1)\int_{\R}Q_s^p(y)u(y)\,dy.
\end{equation}

We now turn to compute $I_2$, by the decay of $u$ and $Q_s$ at infinity, one has
\[
|I_2(x)|\leq
C|x|^{-(p+1)(1+2s)}
\int_{|y|>|x|/2}G_s(x-y)\,dy\leq C|x|^{-(p+1)(1+2s)}.
\]
This implies that $|x|^{1+2s}|I_2(x)|
\leq C|x|^{-p(1+2s)}
\to 0$ as $|x|\to \infty$.

As a result of this and~\eqref{eq:inner limit},
we prove~\eqref{bvsjdjnk} which is equivalent to~\eqref{eq:u asymptotic}.

\medskip

Finally, if $M\ne0$, the limit just proved supplies $R_0>1$ such that
\[
\left|
|x|^{1+2s} u(x)-c_{1,s}M
\right|
\leq
\frac{c_{1,s}|M|}{2}
\qquad \forall|x|\geq R_0.
\]
This implies that 
\[
|x|^{1+2s}|u(x)|
\geq
\frac{c_{1,s}|M|}{2}.
\]
Thus the lower estimate~\eqref{eq:u lower} follows with
\[
c_u=
\min\left\{
\frac12,\frac{c_{1,s}|M|}{2}
\right\}
\in(0,1).
\]
This ends the proof.
\end{proof}

Assume that $Q_s\in \mathcal G_s$, and $h\in\ker L_{s, Q_s}\cap L^2_{even}(\R)$ with $h\neq 0$, 
in view of Lemma~\ref{lemma unique zero}, we know that  $h$ has a unique positive nodal point $r$ 
such that
\[
h(0)>0,\qquad
h>0\text{ on }(0,r),\qquad
h<0\text{ on }(r,\infty).
\]
Moreover, by~\eqref{nod positive moments}, we have $\int_{\R} h<0$. Thus we are in a position to apply the lower bound in~\eqref{eq:u lower} to $h$, and derive the following the monotonicity of $h/(-Q_s')$.

\begin{Lemma}\label{lemma comparision monotonicity}
    Let $0<s<1$ and $p>0$. Let $Q_s\in \mathcal G_s$. Suppose  that $h\in\ker L_{s, Q_s}\cap L^2_{even}(\R)$ with $h(0)>0$.  Denote 
\[
w(x)=\frac{h(x)}{-Q_s'(x)}\quad \text{for every } x>0.
\]
Then $w$ is non-increasing on $(0,\infty)$ and
\begin{equation}\label{localpairing}
\langle - Q_s'',h\rangle_{L^2}\ge0.
\end{equation}
\end{Lemma}
\begin{proof}
From the previously established regularity in Lemma~\ref{lemma: u tail at infinity}, we know that  $h,Q_s\in H^2(\R)\cap L^1(\R)\cap C^2(\R)
$, and hence $L_{s,Q_s} Q_s'=0.$ 
Since $Q_s'(0)=0 $ and $
Q_s'(x)<0$ for every $x>0$, the quotient $w(x)$ is well defined and locally $C^1$ on $(0,\infty)$.
By continuity, $h(x)\to h(0)>0 $ and $
-Q_s'(x)\to 0 $ as $
x\to 0.$ Therefore
\begin{equation}\label{rig left ratio}
\lim_{x\to 0}w(x)=+\infty.
\end{equation}

\medskip

We claim that
\begin{equation}\label{rig tail ratio}
\inf_{x>b}w(x)=-\infty
\qquad\text{for every } b>0.
\end{equation}
First, 
since $h$  satisfies~\eqref{eq:linearized equation}, by~\eqref{nod positive moments}, we see that $\int_{\R}h\,dx<0.$
Thus Lemma~\ref{lemma: u tail at infinity} implies that there exists $R_0>1$ such that
\[|h(x)|\geq c_{h}|x|^{-1-2s}\qquad \forall |x|>R_0.\]
As a consequence of this, for every $s\in (0,1)$, we derive 
\begin{equation}\label{eq:infinite first moment}
\int_{\R}|x|^{2s}\,|h(x)|\,dx\geq \int_{|x|>R_0}c_{h}|x|^{-1}\,dx =+ \infty.
\end{equation}

Furthermore,  by the decay of $Q_s$, we obtain $x^{2s}Q_s(x)\to 0$ as $|x|\to \infty$. For each $a>0$, integration by parts yields
\begin{equation}\label{eq:translation first moment}
-\int_a^\infty x^{2s}Q_s'(x)\,dx
=a^{2s}Q_s(a)+2s\int_a^\infty x^{2s-1}Q_s(x)\,dx
<\infty.
\end{equation}
Suppose now that \eqref{rig tail ratio} fails.
Then there exist $b>0$ and $C>0$ such that
$w(x)\ge -C $ for all $
x>b.$ 
Choose $\bar R>\max\{b,r\}.$
Since $h<0$ and $-Q_s'>0$ on $(\bar R,\infty)$, 
using \eqref{eq:translation first moment}, we find
\[
\int_{\bar R}^\infty x^{2s}|h(x)|\,dx
\le-
C\int_{\bar R}^\infty x^{2s}Q_s'(x)\,dx
<\infty.
\]
By continuity and boundedness of $h$, we see that $\int_0^{\bar R} x^{2s}|h(x)|\,dx<\infty.$
Finally, since  $h$ is even, one has 
\[
\int_{\R}|x|^{2s}\,|h(x)|\,dx
=
2\int_0^\infty x^{2s}|h(x)|\,dx
<\infty,
\]
which  contradicts~\eqref{eq:infinite first moment}.

\medskip
We now show that the quotient is non-increasing.

Suppose, argue by contradiction, that there exist $0<x_1<x_2$
such that $w(x_1)<w(x_2).$
Choose a  number $c$ satisfying
$ w(x_1)<c<w(x_2). $
By~\eqref{rig left ratio}, we can choose
$x_0\in(0,x_1)$ such that
$w(x_0)>c.$
By~\eqref{rig tail ratio}, we can also choose
$x_3>x_2$ such that $w(x_3)<c.$

Consider the function $u=h+cQ_s'.$ We observe that  $L_{s,Q_s}u=0.$
Moreover, $u(0)=h(0)+cQ_s'(0)=h(0)>0$, hence $u$ is nontrivial.

For every $x>0$, one has $u(x)=-Q_s'(x)\bigl(w(x)-c\bigr).$
Since $-Q_s'(x)>0$, we deduce that
\begin{equation}\label{sdgbfrdebgfde}
    u(x_0)>0,\qquad
u(x_1)<0,\qquad
u(x_2)>0,\qquad
u(x_3)<0.
\end{equation}
Since  $0$ is the second eigenvalue of $L_{s,Q_s}$. Theorem~\ref{theorem second eigenvalue } guarantees that $u$ changes its sign at most twice on $\R$,
which contradicts~\eqref{sdgbfrdebgfde}. 
We thus conclude that $w$ is non-increasing.
Since $w$ is locally $C^1$, we have 
\[
w'(x)\le0
\qquad\text{for every }x>0.
\]

\medskip
Finally, we prove~\eqref{localpairing}.

For $0<\epsilon<T<\infty$, since $w$ is $C^1$ away from the origin, one has 
\begin{align}
\int_\epsilon^T(-Q_s'')(x)h(x)\,dx
=-\frac12Q_s'(x)h(x)|_\epsilon^T
-\frac12\int_\epsilon^T [-Q_s']^2w'(x)\,dx.
\label{rig truncated pairing}
\end{align}
We notice that,   by continuity, $Q_s'(\epsilon)h(\epsilon)\to 0,$ as $\epsilon \to 0$. Also, owing to the decay estimate~\eqref{eq:u asymptotic}, we derive  $Q_s'(T)h(T)\to 0$ as $T\to\infty$.
Furthermore, we find that
\begin{equation}\label{sdvbgdfb}
\int_0^\infty |Q_s''h|\,dx
\le
\|Q_s''\|_{L^2(0,\infty)}
\|h\|_{L^2(0,\infty)}
<\infty.
\end{equation}
Thus by rewriting~\eqref{rig truncated pairing}, we see that  
\[
\int_{\epsilon}^{T}(-Q_s')^2(-w')\,dx
=Q_s'(x)h(x)|_{\epsilon}^{T}
+2\int_{\epsilon}^{T}(-Q_s'')h\,dx\to 2\int_{0}^{\infty}(-Q_s'')h\,dx, \quad \text{as } \epsilon\to 0,\,  T\to\infty,
\]
thanks to the dominated convergence theorem.
On the left-hand side, the integrand is nonnegative
and the intervals increase to $(0,\infty)$.
Thus monotone convergence theorem and the estimate~\eqref{sdvbgdfb} yield
\[
\int_{0}^{\infty}(-Q_s')^2w'\,dx
= -2\int_{0}^{\infty}(-Q_s'')h\,dx
<\infty.
\]
In particular,
\[
\int_0^\infty(-Q_s'')h\,dx
=
-\frac12\int_0^\infty (-Q_s')^2w'\,dx
\ge0.
\]

Finally, $Q_s''$ and $h$ are even, we conclude that
\[
\langle -Q_s'',h\rangle_{L^2}
=
-\int_0^\infty (-Q_s')^2(x)w'(x)\,dx
\ge0.
\]
This ends the proof.
\end{proof}

With the preliminary work at hand, we now complete the proof of Lemma~\ref{th:nondegeneracy s<1/2}.

\begin{proof}[Proof of Lemma~\ref{th:nondegeneracy s<1/2}]
    Argue by contradiction. Suppose that $h\in \ker L_{s,Q_s}\cap L^2_{even}(\R)$ with $h\neq 0$, then Lemmata~\ref{lemma unique zero} and~\ref{lemma: u tail at infinity} imply that there exists a unique $r>0$ such that $h(r)=0$, 
    \begin{equation}\label{fdfdsgds}
        h(0)>0, \qquad h(x)>0 \quad \forall x\in (0,r)\qquad\text{and}\qquad  h(x)<0 \quad \forall x\in (r,\infty).
    \end{equation}
    
On the other hand,  by~\cite[Lemma~3.2]{MR5122449}, we know that $|Q'(x)|\leqslant C(1+|x|)^{-2-2s}$ for every $x\in\R$. Thus, for every $x\in\R$, 
\[|xQ_s'(x)|\leqslant C(1+|x|)^{-1-2s}.\]
Hence $xQ_s'\in L^2(\R)$.
\medskip

We claim that $xQ_s'\in H^2(\R)$
and $L_{s,Q_s}(xQ_s')=-2Q''+2s(-\Delta)^sQ.$

Indeed,   for simplicity, set $z=xQ_s'$, $H_*=-\Delta+(-\Delta)^s+1$ and $V=(p+1)Q_s^p$.
For $\lambda>0$, $\lambda\ne1$, define
\[
Q_{s,\lambda}(x)=Q_s(\lambda x),
\qquad
v_\lambda=\frac{Q_{s,\lambda}-Q_s}{\lambda-1}\in H^2(\R).
\]
The fundamental theorem of calculus yields
\[
Q_{s,\lambda}(x)-Q_s(x)
=\int_1^\lambda \theta^{-1}z(\theta x)\,d\theta.
\]
Since dilation is strongly continuous on $L^2$ and
$z\in L^2$, it follows that
\begin{equation}\label{dsgvklfdsvg}
    v_\lambda\longrightarrow z
\qquad\text{in }L^2
\quad\text{as }\lambda\to1.
\end{equation}
Moreover, the scaling identities for $H_*$ give 
\[
H_*Q_{s,\lambda}(x)
=
Q_{s,\lambda}^{p+1}(x)
+(\lambda^2-1)(-Q_s'')(\lambda x)
+(\lambda^{2s}-1)[(-\Delta)^sQ_s](\lambda x).
\]
Subtracting $H_*Q_{s}=Q_s^{p+1}$ and dividing by $\lambda-1$, we obtain
\[
H_*v_\lambda(x)
= c_\lambda v_\lambda
+(\lambda+1)(-Q_s'')(\lambda x)
+\frac{\lambda^{2s}-1}{\lambda-1}
[(-\Delta)^sQ_s](\lambda x),
\]
where
\[
c_\lambda(x)=(p+1)\int_0^1
\bigl[(1-t)Q_s(x)+tQ_{s,\lambda}(x)\bigr]^p\,dt.
\]
These coefficients are uniformly bounded and converge
pointwise to $V$. Hence, by dominated convergence theorem and
the $L^2$ convergence of $v_\lambda$, we see that $c_\lambda v_\lambda\to Vz $ in $
L^2.$
Since $Q_s'',(-\Delta)^sQ_s\in L^2$, strong continuity of dilation
also implies 
\[
H_*v_\lambda\longrightarrow
g:=Vz+2(-Q_s'')+2s(-\Delta)^sQ_s
\qquad\text{in }L^2.
\]

The inverse $H_*^{-1}:L^2\mapsto H^2$ is a bounded map.  Hence
\[
v_\lambda=H_*^{-1}(H_*v_\lambda)
\longrightarrow H_*^{-1}g
\qquad\text{in }H^2.
\]
Comparing this with~\eqref{dsgvklfdsvg}, we see that  $z=H_*^{-1}g\in H^2.$
Therefore
\[
L_{s,Q_s}z=H_*z-Vz=-2Q_s''+2s(-\Delta)^sQ_s.
\]
This proves the claim.
\medskip

In particular,  $h\in H^2(\R)$ satisfies $L_{s,Q_s}h=0$, we obtain 
\[
0=\langle L_{s,Q_s}h,xQ_s'\rangle_{L^2} =
\langle h,L_{s,Q_s}(xQ_s')\rangle_{L^2}
=
2\langle -Q_s'',h\rangle_{L^2}
+
2s\langle (-\Delta)^sQ_s,h\rangle_{L^2}.
\]
Combining this with~\eqref{nod positive moments}, owing to~\eqref{localpairing}, one has 
\[
0\leq\langle -Q_s'',h\rangle_{L^2}=\frac s{1-s}\langle Q_s,h\rangle_{L^2}<0.\]
This is impossible. We thus complete the proof.
\end{proof}

\subsubsection{Nondegeneracy in the odd sector}
We now  identify the real-valued odd kernel of this linearized operator.

\begin{Lemma}\label{lemma:nondegeneracy2}
Let  $0<s<1$, $p>0$, and  $Q\in \mathcal M_s$. Then
\begin{equation}\label{sdbvgdfbgf}
    \ker L_{s,Q}\cap L^2_{\mathrm{odd}}(\R)
=\operatorname{span}\{Q'\}.
\end{equation}
\end{Lemma}
\begin{proof}
It is  known that $Q'\in\ker L_{s,Q}\cap L^2_{\mathrm{odd}}(\R).$
Suppose, for contradiction, that the odd kernel has dimension at least two. Then we can choose a nontrivial function
$\phi\in\ker L_{s,Q}\cap L^2_{\mathrm{odd}}(\mathbb R)$
such that
\begin{equation}\label{3.12}
    \langle \phi,Q'\rangle_{L^2(\R)}=0,
\end{equation}
where $Q',\phi\in H^2(\R)$. Indeed, choose $\psi\in \ker L_{s,Q}\cap L^2_{\mathrm{odd}}(\R)$ and $\psi\not\in\operatorname{span}Q' $. Set $$\phi =\psi-\frac{\langle \psi,Q'\rangle_{L^2(\R)}}{\|Q'\|_2^2}Q',$$ then $\phi\not\equiv 0$, $L_{s,Q}\phi=0$ and the identity~\eqref{3.12} holds. 

We now denote the restriction of $\phi$ to the positive half-line by
\[
w=\phi|_{(0,\infty)}\in H^1(0,+\infty).
\]
Since $H^1(\R)\subset C^{0,1/2}(\R)$ and
$\phi$ is odd, it follows that $\phi(0)=0,$
and hence $w(0)=0.$
This means that $w\in H^1_0(0,+\infty)$ in the trace sense.

\medskip

We claim that $w$ must change sign on $(0,\infty)$.

Indeed, if $w\not\equiv0$ has one sign, we could replace $\phi$ by $-\phi$
and assume that $w\ge0$. Since $Q'(r)<0$ for every $r>0$, one has 
\[
\langle \phi,Q'\rangle_{L^2(\R)}
=
2\int_0^\infty w(r)Q'(r)\,dr
<0,
\]
which is impossible, thanks to~\eqref{3.12}. 

Therefore, $w$ must take both positive and negative values.
Since $\phi\in H^1(\R)\subset C^{0,1/2}(\R)$,
the positive and negative sets both have positive measure.

Define the
odd function $\phi^*$ by
\[ 
\phi^*(x)=| w(x)| \quad \text{if } x>0, \qquad \phi^*(0)=0,\qquad  \phi^*(x)= -| w(-x)| \quad \text{if } x<0.
\]
It follows that $\phi^*\in H^1(\R)\cap L^2_{\mathrm{odd}}(\R)$. By the strict inequality~\eqref{ksdjsdf}, one has
\[
\langle\phi^*,L_{s,Q}\phi^*\rangle
=\int_{\mathbb R}|(\phi^*)'|^2\,dx
+\int_{\mathbb R}|(-\Delta)^{s/2}\phi^*|^2\,dx
+\int_{\mathbb R}(1-(p+1)Q^p)|\phi^*|^2\,dx
<\langle \phi,L_{s,Q}\phi\rangle=0,
\]
since $L_{s,Q}\phi=0$.
This contradicts the nonnegativity~\eqref{sdvkidsjmf}  of the odd quadratic form.

Thus, the odd kernel has dimension at most one. Since it contains the nontrivial function $Q'$, we obtain the desired result~\eqref{sdbvgdfbgf}.
\end{proof}

\begin{remark}
By Lemma~\ref{lemma:nondegeneracy2}, the odd kernel of
$L_{s,Q}$ is given by
\[
\ker L_{s,Q}\cap L^2_{\mathrm{odd}}(\mathbb R)
=\operatorname{span}\{Q'\}.
\]
This identity holds without the Morse index assumption, which is needed in our nodal argument to exclude nontrivial
even kernel elements.
\end{remark}

\subsection{Uniqueness}\label{sec: uniqueness s<1/2}

Our goal this section is to prove Theorem~\ref{th:nondegeneracy low s}, we can now closely follow the approach provided in \cite{MR3070568,MR3530361}. That is, beginning with a given $Q_{s_0}$ with $s_0 \in (0, 1)$, we construct a branch $s \mapsto Q_s$ satisfying
	\begin{equation}
	-\Delta Q_s+(-\Delta)^s Q_s + Q_s = |Q_s|^{p}Q_s \quad \text{in } \mathbb{R}.
	\end{equation}
	
	First, based on the invertibility of the linearized operator $L_{s,Q_s}$ around $Q_{s_0}$ in the even sector, we are in a position  applying an implicit function Theorem to prove 
	the local existence and uniqueness of $Q_s$ for $s$ close to $s_0$ (see Section~\ref{sec:Construction of the Local Branch}).
	
	With this local branch at hand, we derive a priori bounds for $Q_s$ that allows us to continue the solution to $s = 1$, linking it to the classical setting where $u_{s=1}$ is known to be unique and nondegenerate \cite{MR969899,MR886933}. This strategy culminates in a proof that the branch $Q_s$ itself is unique, which in turn establishes the uniqueness of the initial ground state $Q_{s_0}$ (see Sections~\ref{sec:A-Priori Bounds and Global Continuation} and~\ref{sec:The proof of Theorem uniqueness}).

	To begin with, we need to introduce the following fundamental estimate.   Let us  denote by $\bot_s$ the orthogonality relation in $H^{1}(\mathbb{R})$.
	
	\subsubsection{Construction of the Local Branch}\label{sec:Construction of the Local Branch}
	
	Let  $s \in (0, 1)$. We define the real Banach space $X_p$ as the set of real-valued, even functions
	\begin{equation}
		X_p :=\left\{f\in L^2(\mathbb{R}) \cap L^{p+2}(\mathbb{R}): f \text{ is even } \right\},
	\end{equation}
	equipped with the norm $\|{f}\|_{X_p} \coloneqq \|{f}\|_{L^2} + \|{f}\|_{L^{p+2}}.$

	For $s \in [s_0, 1)$,  let $Q \in X_p$ be a solution to the following equation
	\begin{equation}\label{bvjad}
		-\Delta Q+ (-\Delta)^s Q + Q - |Q|^{p}Q  = 0 \quad \text{in } \mathbb{R},
	\end{equation}
in the sense of distributions. Then,  indeed $Q \in H^{2}(\mathbb{R})\cap C^{2,\alpha}(\mathbb{R})$, thanks to  Theorems~1.3 and~1.4 in~\cite{DSVZ24}.
We nevertheless formulate the problem in $X_p$. This choice is motivated by the fact that $X_p$ is independent of $s$ and gradient of $ Q$, providing a consistent framework for the continuation argument that constructs a local branch $Q_s$ with $s \in [s_0, 1)$. At this stage, $Q$ is not required to be positive.

\begin{Lemma}\label{Lemma implicit theorem}
	Let  $0 < s_0 < 1$, and $p>0$. Suppose that $Q_{s_0} \in X_p$ solves equation~\eqref{bvjad} with $s = s_0$. Moreover, assume that the linearization
	\[
	L_{s_0,Q_{s_0}} = -\Delta+ (-\Delta)^{s_0} + 1 - (p + 1)|Q_{s_0}|^p
	\]
	has trivial kernel on $L^2_{\mathrm{even}}(\mathbb{R})$.
	
	Then, for some $\delta > 0$, there exists a map $Q \in C^1(I; X_p)$ defined on the interval $I = [s_0, s_0 + \delta)$ such that the following holds, where we denote $Q_s = Q(s)$ in the sequel.
	\begin{itemize}
		\item[(i)] $Q_s$ satisfies~\eqref{bvjad} for all $s \in I$.
		\item[(ii)] There exists $\epsilon > 0$ such that $Q_s$ is the unique solution of~\eqref{bvjad} for $s \in I$ in the neighborhood $\{ Q \in X_p : \|{Q - Q_{s_0}}\|_{X_p} < \epsilon \}$. In particular, we have that $Q_{s=s_0} = Q_{s_0}$ holds.
	\end{itemize}
\end{Lemma}
	\begin{proof}
As in~\cite{MR3070568,MR3530361}, we apply an implicit function argument to the map
\begin{equation}
	F : X_p \times [s_0, s_0 + \delta) \to X_p, \quad F(Q, s) = Q - \frac{1}{-\Delta+(-\Delta)^s + 1} |Q|^p Q.
\end{equation}
By construction,  $F(Q, s) = 0$ is equivalent to $Q \in X_p$ being a solution of~\eqref{bvjad}, and we have $F(Q_{s_0}, s_0) = 0$ by hypothesis. Moreover, by combining Young's with  H\"older's inequalities, it is easy to check that $F$ is a well-defined $C^1$ mapping, whence the Fr\'echet derivative with respect to $Q$ at $(Q_{s_0}, s_0)$ is given by
\begin{equation}
	\partial_Q F(Q_{s_0}, s_0) = 1  -\frac{1}{-\Delta+(-\Delta)^{s_0} + 1} (p + 1)|Q_{s_0}|^p:=1+K.
\end{equation}		

We claim that  $\partial_Q F$ has a bounded inverse on $X_p$ at $(Q_{s_0}, s_0)$.

Indeed,  we notice that the  operator $K$ maps $X_p$ to $X_p$ and  is compact on $X_p$ by combining a cutoff argument and the fact that $Q_{s_0}$ vanishes at infinity. Furthermore, by assumption, we know  that  $\ker L_{s_0,Q_{s_0}}\cap L_{\rm even}^2(\R)=\{0\} $, which implies that $-1\not\in \sigma(K)$.  We thus deduce that  the bounded inverse $(1+K)^{-1}$ exists on $X_p$ at $(Q_{s_0}, s_0)$.

 As a result, the implicit function theorem yields that the assertions (i) and (ii)  hold for some $\delta > 0$ and $\epsilon > 0$ sufficiently small.
	\end{proof}

	\subsubsection{A-Priori Bounds and Global Continuation}\label{sec:A-Priori Bounds and Global Continuation}
	
	We aim to  globally extend the locally unique branch $Q_s$, with $s \in [s_0, s_0 + \delta)$, constructed in Lemma~\ref{Lemma implicit theorem}. To this end, we  adapt the general strategy derived in~\cite{MR3070568,MR3530361} for the fractional Laplacian case. 
	
	In what follows, let  $0 < s_0 < 1$, and $p>0$ be fixed. Furthermore,  suppose that $Q_{s_0} \in X_p$ is fixed satisfying the hypotheses of  Lemma~\ref{Lemma implicit theorem}. Correspondingly, let  $Q_s \in C^1(I;X_p)$ with $I = [s_0, s_0 + \delta)$ denote the unique local branch provided by Lemma~\ref{Lemma implicit theorem}. We consider the maximal extension of the branch $Q_s$ with $s \in [s_0, s_*)$, where $s_*$ is given by
\[
s_* := \sup \left\{ s_0 < \tilde{s} < 1 : 
\begin{array}{l}
	\qquad Q_s \in C^1([s_0, \tilde{s}); X_p)  \text{ given  by  Lemma~\ref{Lemma implicit theorem} }\\ \text{and } Q_s 
	\text{ satisfies the assumptions of } \text{Lemma~\ref{Lemma implicit theorem} for } s \in [s_0, \tilde{s})
\end{array}
\right\}.
\]
\begin{remark}\label{remark property of Q} Let $Q_s$ is given by Lemma~\ref{Lemma implicit theorem}.
    Since $Q_s\in X_p$ satisfies~\eqref{bvjad}, using the decay of the kernel $G_s$ at infinity and standard regularity argument, it is immediate that $Q_s\in C^2(\R)$, 
    \[|Q_s(x)|\leqslant C(1+|x|)^{-1-2s},\qquad |Q_s'(x)|\leqslant C(1+|x|)^{-2-2s}.\]
    Hence $Q_s\in L^1(\R)$ and $xQ_s'\in L^2(\R)$. Here $Q_s$ is not required to be positive.
\end{remark}

Our  goal is to prove that $s_*=1$ when $Q_{s_0}\in \mathcal G_s$. For this, we first establish a priori bounds along the maximal branch $Q_s$ for $s\in[s_0,s_*)$. 
	
	\begin{Lemma}\label{lemma priori bound}
		We have the a-priori bounds
		\[
		\int_{\mathbb{R}} | Q_s'|^2 + \int_{\mathbb{R}} |(-\Delta)^{s/2}Q_s|^2 \sim \int_{\mathbb{R}} |Q_s|^{p+2} \sim 1
		\]
		for all $s \in [s_0, s_*)$. Here, we use $X\sim Y $ to denote that $C^{-1} Y\leqslant X\leqslant C Y$ for some constant $C>0$.
	\end{Lemma}
\begin{proof}
	We divide the proof into the following  steps.
	
	\textbf{Step 1. Lower bounds.} We denote the following notations
	\begin{equation}\label{sbvjsdnv}
	M_s = \int_{\mathbb{R}} |Q_s|^2, \quad 	N_s = \int_{\mathbb{R}} | Q_s'|^2, \quad  T_s= \int_{\mathbb{R}} |(-\Delta)^{s/2}Q_s|^2,\quad 
		V_s = \int_{\mathbb{R}} |Q_s|^{p+2},
	\end{equation}
	for $s \in [s_0, s_*)$. We first claim that
	\begin{equation}\label{sdvs}
		 T_s+N_s \sim V_s.
	\end{equation}
Indeed, integrating equation~\eqref{bvjad} against $Q_s$ and $x \cdot  Q_s'$, respectively, one derives that
\begin{equation}\label{bvjsdd}
	N_s+T_s + M_s = V_s, \quad\text{and}\quad \frac{-N_s}{2}+\frac{1 - 2s}{2}T_s + \frac{1}{2}M_s = \frac{1}{p + 2}V_s.
\end{equation}
As a straightforward consequence of~\eqref{bvjsdd}, we see that \begin{equation}\label{keybdks}
	(1-s)T_s+\left(\frac{1}{2}-\frac{1}{p+2}\right)V_s=N_s+T_s .
\end{equation}
Since $0<\frac{1}{2}-\frac{1}{p+2}<1 $ and $T_s\geqslant 0$, one concludes that $ T_s+N_s \sim V_s.$

On the other hand, we notice that the classical Gagliardo-Nirenberg-Sobolev inequality implies that 
	$N_s^{\frac{p}{4}} M_s^{\frac{p+2}{2}-\frac{p}{4}} \geqslant C_{opt}V_s,$
where the optimal constant $C_{opt}$ is independent of $s $. From this, one can infer that
\begin{equation}\label{sdbcksd}
	N_s+ T_s \sim V_s \geqslant  C_0>0,
\end{equation}
that is, the quantities are uniformly bounded away from zero. It thus  remains to find an upper bound for one of the quantities  $N_s+T_s$ or $V_s$.

\textbf{Step 2. Upper bounds.}
We first claim that
\begin{equation}\label{svbajs}
	\left\| (-\Delta)^{1 - \frac{p}{4(p+2)}} Q_s \right\|_{L^2} \leqslant C V_s^{\frac{p+1}{p+2}},
\end{equation}
for some constant $C>0$ independent
of $s$. (This is a regularity statement since $1 - \frac{p}{4(p+2)}> 1/2$.) 

Indeed, recall that $Q_s$ solves equation~~\eqref{bvjad}, employing Plancherel's identity, one
deduces that
\begin{equation}\label{snvkds}
	\| (-\Delta)^{1 - \frac{p}{4(p+2)}} Q_s \|_{L^2} = \left\| \frac{(-\Delta)^{1 - \frac{p}{4(p+2)}}}{-\Delta+(-\Delta)^s + 1} |Q_s|^p Q_s \right\|_{L^2}
	\leqslant \| (-\Delta)^{- \frac{p}{4(p+2)}} (|Q_s|^p Q_s) \|_{L^2}
\end{equation}
Since $ - \frac{p}{4(p+2)}>-1/2$, the operator thus $(-\Delta)^{- \frac{p}{4(p+2)}}$ on $\mathbb{R}$ is given by convolution with
$|x|^{-\left(1-\frac{p}{2(p+2)}\right)}$ up to a multiplicative constant depending only on $p$. Hence,
the Hardy-Littlewood-Sobolev inequality yields  that
\begin{equation}\label{keysv}
	\| (-\Delta)^{1 - \frac{p}{4(p+2)}} Q_s \|^2_{L^2} \leqslant C_{p} \left\| |x|^{-\left(1-\frac{p}{2(p+2)}\right)} \ast (|Q_s|^p Q_s) \right\|^2_{L^2} 
	\leqslant C_{p} \left\| |Q_s|^{p+1} \right\|^2_{L^{\frac{p+2}{p+1}}} =C_{p} V_s^{\frac{2(p+1)}{p+2}}.
\end{equation}
Hence, the claim~~\eqref{svbajs} is proved.

With the estimate~\eqref{svbajs}  at hand, we shall provide a uniform upper bound for $V_s$ as follows.  Let us differentiate
the equation satisfied by $Q_s$ with respect to $s$, one has that 
$L_{s,Q_s} \dot{Q_s} = -(-\Delta)^s \ln (-\Delta)Q_s$
with $\dot{Q_s} = {dQ_s}/{ds}$. From 
$L_{s,Q_s}Q_s = -p |Q_s|^{p}Q_s$, it follows  that
\begin{equation}\label{dsvsdj}
	\frac{d}{ds} V_s = (p+2) \left\langle |Q_s|^{p}Q_s, \dot{Q_s}  \right\rangle_{L^2} = \frac{p+2}{p} \left\langle Q_s, (-\Delta)^s \ln(-\Delta)Q_s  \right\rangle_{L^2}.
\end{equation}

Owing to the formulas~\eqref{sdvs} and~\eqref{svbajs},  one has that
\begin{equation}
	\begin{split}\label{kesdvsy}
		\left\langle Q_s, (-\Delta)^s \ln(-\Delta)Q_s  \right\rangle_{L^2}
		&\leqslant  2\int_{|\xi|\leq R} |\xi|^{2s} (\ln |\xi|) |\hat{Q}_s(\xi)|^2  d\xi + 2 \int_{|\xi|>R} |\xi|^{2} (\ln |\xi|) |\hat{Q}_s(\xi)|^2  d\xi \\
		&\leqslant 2(\ln R+s^{-1}) \int_{\mathbb{R}} |\xi|^{2s}|\hat{Q}_s(\xi)|^2  d\xi + 2R^{\frac{p}{p+2}-2}(\ln R) \int_{\mathbb{R}} |\xi|^{4-\frac{p}{p+2}}|\hat{Q}_s(\xi)|^2  d\xi \\
		&\leqslant 2(\ln R+s_0^{-1})T_s + 2R^{\frac{p}{p+2}-2}(\ln R)V_s^{\frac{2(p+1)}{p+2}}\\
		&\leqslant 2(\ln R+s_0^{-1})V_s + 2R^{\frac{p}{p+2}-2}(\ln R)V_s^{\frac{2(p+1)}{p+2}},
	\end{split}
\end{equation}
provided that we take $R > e^{\frac{1}{2-\frac{p}{p+2}}}>1$.

Moreover, we observe that, there exists  positive constant $c_0$ independent of $s$ such that 
\[e<R^{2-\frac{p}{p+2}}=c_0V_s^{\frac{p}{p+2}}, \]
thanks to the estimate~\eqref{sdbcksd}. By combining~\eqref{dsvsdj} and~\eqref{kesdvsy}, one derives that
\begin{equation}\label{keysdds}
		\frac{d}{ds} V_s \leqslant C_{s_0,p,c_0} \left(1+\ln\left( c_0^{1+\frac{2}{p}} V_s\right)\right)V_s.
\end{equation}
Utilizing integration yields the uniform upper bound $V_s\leqslant C$ for $s\in[s_0,s_*)$, where the positive constant $C$ depends on $s_0, p,V_{s_0}$.  The
proof of Lemma~\ref{lemma priori bound} thus is now complete.
\end{proof}
	
	As a first application of the a-priori bounds from Lemma~\ref{lemma priori bound}, we establish the continuous dependence of the linearized operators in the norm-resolvent sense with respect to $s$.
	
	We consider two key self-adjoint and lower semi-bounded operators on $L^2(\mathbb{R})$ arising from the linearization: 
	\[
	L_{s,Q_s} = -\Delta+ (-\Delta)^s + 1 - (p + 1)|Q_s|^p\qquad \text{and}\qquad L_{-,s,Q_s} =-\Delta+ (-\Delta)^s + 1 - |Q_s|^p.
	\]

	 Recall that convergence
	of a sequence $\{A_n\}_{n=1}^{\infty}$ of self-adjoint operators to an operator $A$ in norm-resolvent sense
	means that 
	\[
	\left\| \frac{1}{A_n + z} - \frac{1}{A + z} \right\|_{L^2 \to L^2} \to 0
	\]
	as $n \to \infty$, where $z \in \mathbb{C}$ with $\operatorname{Im} z \neq 0$.
	This notion of convergence, which is stronger than strong resolvent convergence, ensures the continuity of spectral properties.
	
	\begin{Lemma}\label{Lemma operator continunity}
		Let $\{s_k\}_{k=1}^{\infty} \subset [s_0, s_*)$ be a sequence with $s_k \to \tilde{s} < s_*$. Then  $L_{s_k,Q_{s_k}} \to L_{\tilde{s},Q_{\tilde{s}}}$ and $L_{-,s_k,Q_{s_k}} \to L_{-,\tilde{s},Q_{\tilde{s}}}$
		in norm-resolvent sense.
	\end{Lemma}
\begin{proof}
	A straightforward calculation yields that norm-resolvent convergence  holds, provided we can show that $V_s = |Q_s|^p$ satisfies
	\begin{equation}\label{dsvsv}
		\|{V_s - V_{\tilde{s}}}\|_q \to 0 \quad \text{as } s \to \tilde{s}
	\end{equation}
	for some $q > 1$.
	
	For this, we notice that, from  Lemma~\ref{lemma priori bound},
	it follows that  the $H^1$ norm of ${Q_s}$ is uniformly bounded with respect to $s$. Recall that $Q_s \to Q_{\tilde{s}}$ in $X_p$, we thus deduce  that
	\[
	\|{Q_s - Q_{\tilde{s}}}\|_{H^{1}} \to 0 \quad \text{as } s \to \tilde{s}.
	\]
	Hence, by Sobolev inequalities
	and H\"older's inequality,  we can take $q=1+2/p$ for $p>1$, and $q={2}/{p} $ for $0<p\leqslant 1$.
\end{proof}

As a result of the norm-resolvent convergence,  it is immediate that the Morse index of $L_{s,Q_s}$  remain constant along the solution branch, as follows.

\begin{Lemma}\label{lemma morse index continutity}
	We denote the Morse index of $L_{s,Q_s}$ in the sector of even functions by 
 \[\mathcal{N}_{\rm even}(L_{s,Q_s}):=\sharp\left\{e<0:\, e \text{ is a eigenvalue of } L_{s,Q_s} \text{ restricted on } L^2_{\rm even}(\R)\right\}\]
 where multiplicities of eigenvalues are taken into account.
 
	 Then,  we have  $\mathcal{N}_{\rm even}(L_{s,Q_s})=\mathcal{N}_{\rm even}(L_{s_0,Q_{s_0}})$ for all $s \in [s_0, s_*)$.
\end{Lemma}

	Based on Lemmata~\ref{lemma priori bound} and~\ref{Lemma operator continunity}, we next prove  a uniform decay bound with respect to $s$.
	
	\begin{Lemma}\label{Lemma uniform decay}
	Suppose that $Q_{s_0}(|x|) > 0$ is positive. Then, for all $s \in [s_0, s_*)$, we have
	\begin{equation}\label{ndvsdi}
	    Q_s(|x|) > 0 \quad \text{for } x \in \mathbb{R}, \qquad Q_s(|x|) \leqslant  C|x|^{-1} \quad \text{for } |x| \geqslant 1,
	\end{equation} for some positive constant $C$ independent of $s$.
	\end{Lemma}
	\begin{proof}
With the Lemma~\ref{Lemma operator continunity} at hand,	from the proof of Corollary~5.10 in~\cite{MR3070568}, in virtue of the  fact that the operator $L_{-,s,Q_s}$ enjoys a Perron–Frobenius property (see~\cite[Lemma~B.4]{SZZ2026}),  the standard spectral arguments implies  the positivity of $Q_s(|x|) > 0$ for $ s\in[s_0, s_*)$, provided that $Q_{s_0}(|x|) >0 $ initially
holds.
\medskip

We now turn to derive the uniform decay bound in~\eqref{ndvsdi}.

Recall that $Q_s \in L^1(\mathbb{R})$. By integrating equation~\eqref{main equation} over $\mathbb{R}$, and employing the fact that $\int_{\mathbb{R}} (-\Delta)^s Q_s = 0$,  we obtain the identity
\begin{equation}\label{dsvjsdbvj}
\int_{\mathbb{R}} - Q_s''  +	\int_{\mathbb{R}} Q_s = \int_{\mathbb{R}} Q_s^{p+1}.
\end{equation}
Furthermore, owing to~\cite[Lemma~C.1]{SZZ2026}, we see that $| Q_s'(x)|\leqslant C|x|^{-1-2s}$ as $|x|>1$ for some positive constant $C$.  This implies that  $\int_{\mathbb{R}} -Q_s'' = 0$ holds. Hence, 
\begin{equation}\label{ascdsvjsdbvj}
		\int_{\mathbb{R}} Q_s = \int_{\mathbb{R}} Q_s^{p+1}.
\end{equation}
By H\"older's inequality, we see that
\[
\int_{\mathbb{R}} Q_s^{p+1} \leqslant \left( \int_{\mathbb{R}} Q_s \right)^{\frac{1}{1+p}} \left( \int_{\mathbb{R}} Q_s^{p+2} \right)^{\frac{p}{(1+p)}}.
\]
As a consequence of this and the a-priori bound from Lemma~\ref{lemma priori bound}, one infers the uniform bound
\begin{equation}\label{sdvnsd}
	\int_{\mathbb{R}} Q_s \leqslant C.
\end{equation}
for some constant $C>0$ independent of $s$.

From Theorem~\ref{th: existence and properity}, we know that the function $Q_s(x) > 0$ is decreasing in $|x|$. Hence, for any $R > 0$,
\begin{equation}\label{sdvasd}
	\int_{\mathbb{R}} Q_s \geqslant \int_{|x|\leqslant R} Q_s(R) \geqslant C R Q_s(R).
\end{equation}
From this and~\eqref{sdvnsd}, we conclude that the estimate~\eqref{ndvsdi} holds, completing the proof of Lemma~\ref{Lemma uniform decay}.	\end{proof}
	
Gathering these facts above,	we conclude the following convergence fact.
	
	\begin{Lemma}\label{Lemma convergence}
		Let $\left\{s_k\right\}_{k=1}^{\infty} \subset [s_0, s_*)$ be a sequence such that $s_k \to s_*$ and suppose that
		$Q_{s_k}(|x|) > 0$ for $k \in \mathbb{N}$. Then, after possibly passing to a subsequence, we have $Q_{s_k} \to u_*$
		in $L^2(\mathbb{R}) \cap L^{p+2}(\mathbb{R})$ as $k \to +\infty$. Moreover, the function $Q_*(|x|) > 0$ is positive and
		satisfies
		\begin{equation}\label{sdvsdvd}
		-\Delta Q_*+	(-\Delta)^{s_*} Q_* + Q_* - Q_*^{p+1} = 0 \quad \text{in } \mathbb{R}.
		\end{equation}
	\end{Lemma}
	\begin{proof}
		Let us denote $Q_k = Q_{s_k}$ in the following.  From  Lemma~\ref{lemma priori bound}, it follows that 
		 $$\|{Q_k}\|^2_{H^{1}}= \|Q_k\|^{p+2}_{L^{p+2}} \leqslant C_{p,s_0,Q_{s_0}}.$$
		 
		This implies that,  passing to a subsequence if necessary, there exists $Q_*\in H^1(\mathbb{R})$ such that  $Q_k \rightharpoonup Q_*$
		weakly in $H^{1}$. Furthermore, owing to  local Rellich compactness, we see that $Q_k \to Q_*$ in
		$L^2_{\mathrm{loc}}(\mathbb{R})$ and pointwise a.e. in $\mathbb{R}$. 
		Recall that by Lemma~\ref{Lemma uniform decay}, the sequence $\{Q_k\}$ has a uniform decay at infinity. 
		Combined with the local compactness already established, this allows us to upgrade the convergence 
		to strong convergence in $L^2(\mathbb{R})$, i.e., $Q_k \to Q_*$ in $L^2(\mathbb{R})$.
		
		On the other hand,  the Sobolev embedding
		ensure the uniform bound $\|Q_k\|_{L^q{(\mathbb{R})}} \leqslant C$. 
		for every $q >2$, we may interpolate with the strong $L^2$-convergence: 
		by H\"older's inequality, we deduce that $Q_k \to Q_*$ in $L^{p+2}(\mathbb{R})$ as well.

	We now claim  that  $Q_*(|x|) $  is a solution of equation~\eqref{sdvsdvd} and is a positive function in $L^2(\mathbb{R}) \cap L^{p+2}(\mathbb{R})$.
		
		 Indeed, since $Q_k(|x|) > 0$ and $Q_k \to Q_*$ pointwise a.e. on
		$\mathbb{R}$, one  thus has that $Q_*(|x|)\geqslant 0$. Furthermore, the uniform lower bounds in
		Lemma~\ref{lemma priori bound} combined with the convergence $Q_k \to Q_*$ in $L^{p+2}(\mathbb{R})$ yield that $Q_* \not\equiv 0$.
		By passing to the limit in the equation satisfied by $Q_k$, we deduce that~\eqref{sdvsdvd} holds. 
		This implies the fixed-point relation:
		$Q_* = (-\Delta+(-\Delta)^{s_*} + 1)^{-1} Q_*^{p+1}$. 
		
		Since $Q_* \geqslant 0$ and $Q_* \not\equiv 0$, using the fact that the kernel of the resolvent operator $(-\Delta + (-\Delta)^{s_*} + 1)^{-1}$ is positive, 
		it follows from the positivity property of the kernel that $Q_*(|x|) > 0$ everywhere.
	\end{proof}
	
We now show that every  branch starting at $Q_{s_0}\in\mathcal G_{s_0}$ can be continued up to $s=1$.
	
	\begin{Lemma}\label{lemma maximal estension}
		Let $Q_{s_0} \in \mathcal G_{s_0}$ be a solution of~\eqref{main equation}. Then its maximal
		branch $Q_s$ with $s \in [s_0, s_*)$ extends to $s_* =1$.
	\end{Lemma}
	\begin{proof}
	From Lemma~\ref{Lemma uniform decay}, it follows that $Q_s(|x|) > 0$ for all $s \in [s_0, s_*)$. 
	Furthermore, 
 Lemma~\ref{lemma morse index continutity} implies that 
 $\mathcal{N}_{\rm even}(L_{s,Q_s})= 1 $ for $ s \in [s_0, s_*),$
where the linearized operators $L_{s,Q_s}$ along the branch defined by 
$L_{s,Q_s} = -\Delta +(-\Delta)^s + 1 - (p + 1)Q_s^p.$

We now pick a sequence $s_k \to s_*$ with $s_k < s_*$ and denote $Q_k = Q_{s_k}$. In view of Lemma~\ref{Lemma convergence},  one has that there exists $Q_*\in H^1(\mathbb{R})$, such that $Q_k \to Q_*$ in $L^2(\mathbb{R}) \cap L^{p+2}(\mathbb{R})$ and $Q_*(|x|) > 0$ satisfies equation~\eqref{sdvsdvd}.

We argue by contradiction and suppose that $s_* < 1$. 
We claim that $Q_*\in \mathcal G_{s_*}$ is a solution to~\eqref{sdvsdvd}. 
If this claim holds, then by Lemma~\ref{th:nondegeneracy s<1/2}, the corresponding linearized operator  has a trivial kernel on $L^2_{\mathrm{even}}(\mathbb{R})$. 
Lemma~\ref{Lemma implicit theorem} would then imply that the branch $Q_s$ can be extended beyond $s_*$, contradicting the maximality of $s_*$.

It remains to show the claim that $Q_*\in \mathcal G_{s_*}$. Let $L_{*}$ denote its
corresponding linearization
\begin{equation}\label{sdvnksd}
	L_{*} = -\Delta+(-\Delta)^{s_*} + 1 - (p + 1)Q_*^p.
\end{equation}
We then have that $L_{s_k,Q_k} \to L_{*}$ in the norm resolvent sense.  Indeed,  we notice that $ Q_k \to Q_*$ in $L^2(\mathbb{R}) \cap L^{p+2}(\mathbb{R})$ and $Q_*$ is a solution to equation~\eqref{sdvsdvd}, this yields that $Q_k \to Q_* $ in $ H^1(\mathbb{R}),$
Thus, by adapting the proof of
Lemma~\ref{Lemma operator continunity}, we deduce that $L_{s_k,Q_k} \to L_{*}$ in the norm resolvent sense.

Since the Morse index is lower
semi-continuous with respect to the norm-resolvent topology, we conclude that
\begin{equation}
	\liminf_{k\to\infty} \mathcal{N}_{\rm even}(L_{s_k,Q_k}) \geqslant \mathcal{N}_{\rm even}(L_{*}).
\end{equation}
This shows that $\mathcal{N}_{\mathrm{even}}(L_{*}) \leqslant 1$.

On the other hand, we observe that, $Q_*$ solves equation~\eqref{sdvsdvd}, thus
\[
\langle Q_*, L_{*} Q_* \rangle_{L^2} = -p \int_{\mathbb{R}} |Q_*|^{p+2}  dx < 0.
\]
As a result of this and the min-max principle, we derive that $L_{*}$ acting on
$L^2_{\mathrm{even}}(\mathbb{R})$ has Morse index equal to $1$. 
Hence $Q_*(|x|) \in \mathcal G_{s_*}$.  This implies that the proof of Lemma~\ref{lemma maximal estension} is  complete.
	\end{proof}

	\subsubsection{The proof of uniqueness}\label{sec:The proof of Theorem uniqueness}
	
	Let $0 < s_0 < 1$, and $p>0$ be fixed.
	\begin{proof}
	   We proceed by contradiction. Suppose there exist two distinct even  positive solutions $Q_{s_0}, \tilde{Q}_{s_0} \in \mathcal G_{s_0}$ for equation~\eqref{main equation} with $s = s_0$, satisfying $Q_{s_0} \not\equiv \widetilde{Q}_{s_0}$.

	By Lemma~\ref{th:nondegeneracy s<1/2}, both $Q_{s_0}$ and $\tilde{Q}_{s_0}$ belong to $X_p$ and satisfy the hypotheses of Lemma~\ref{Lemma implicit theorem}. Combining this with Lemma~\ref{lemma maximal estension}, we obtain two distinct $C^1$ branches $Q_s, \tilde{Q}_s \in C^1([s_0, 1); X_p)$ emanating from $Q_{s_0}$ and $\tilde{Q}_{s_0}$, respectively. The local uniqueness guaranteed by Lemma~\ref{Lemma implicit theorem} ensures that $Q_s \not\equiv \tilde{Q}_s$ for all $s \in [s_0, 1)$.
	
	Now, by Lemma~\ref{Lemma convergence}, as $s \to (1)^-$, we have the convergences
	\[
	Q_s \to Q_* \quad \text{and} \quad \tilde{Q}_s \to \tilde{Q}_* \quad \text{in } L^2(\mathbb{R}) \cap L^{p+2}(\mathbb{R}),
	\]
	where both $Q_*$ and $\tilde{Q}_*$ satisfy the limiting equation~\eqref{sdvsdvd}.
	
	However, it is known that the positive  solution $Q_1(|x|) > 0$ in $L^2(\mathbb{R}) \cap L^{p+2}(\mathbb{R})$ to the local problem
	\begin{equation}
		-2 Q_{1}''  + Q_{1} - Q_{1}^{p+1} = 0 \quad \text{in } \mathbb{R}
	\end{equation}
	is unique. Moreover, its linearized operator $L_1 = -2\Delta + 1 - (p+1)Q_1^p$ is nondegenerate and invertible on $L^2_{\mathrm{even}}(\mathbb{R})$. (see, e.g., Kwong~\cite{MR969899}).
	By uniqueness of the centered positive solution,	it follows that $Q_* \equiv \tilde{Q}_* \equiv Q_{1}$, and hence
	\[
	\|{Q_s - \tilde{Q}_s}\|_{L^2 \cap L^{p+2}} \to 0 \quad \text{as } s \to {1}^-.
	\]
	
	Finally,  the operator $L_*=-2\Delta+ 1 - (p+1)Q_*^p$ satisfies the assumptions of Lemma~\ref{Lemma implicit theorem}, an implicit function argument (similar to that in the proof of Lemma~\ref{Lemma implicit theorem}) implies the existence of a unique branch $Q_s \in C^1((1-\epsilon, 1]; X_p)$ solving~\eqref{bvjad} near $Q_*$, for some small $\epsilon > 0$. This contradicts the earlier conclusion that $Q_s \not\equiv \tilde{Q}_s$ on $[s_0, 1)$.
	This completes the proof of Theorem~\ref{th:nondegeneracy low s}. 
	\end{proof}

	\section{Uniqueness for all positive solutions }\label{sec:Nondegeneracy s ge 1/2}

This section is devoted  to establishing the uniqueness and nondegeneracy results of $Q\in\mathcal{M}_s  $ as stated in Theorem~\ref{th:nondegeneracy s>1/2}.
The proof  will be divided  into three main sectors:
\begin{itemize}
    \item[Sec:\ref{sec:A distribution inequality for the Green kernel}] 
Consider a distribution inequality for the Green kernel:
\begin{equation}\label{eq:green doubling}
F_s(y)=2\int_0^yG_s(x)\,dx,\qquad y>0.
\end{equation}
We prove that \begin{equation}\label{dfbfdgdfvdf}
    F_s(2y)>\frac{2F_s(y)}{1+F_s(y)}.
\end{equation}
    \item[Sec:\ref{sec:sharp proof}] An upper estimate is established by a energy extension, that is, \begin{equation}\label{sharp estimate 1}
Q^p(0)<p+1, \qquad Q\in \mathcal{M}_s.
\end{equation}
    \item[Sec:\ref{sec:Nondegeneracy and uniqueness}] Both conclusions  are obtained by using an energy extension based on a monotonicity formula.
\end{itemize}

\subsection{A distribution inequality for the Green kernel}\label{sec:A distribution inequality for the Green kernel}
For simplify, we denote $\alpha=2s$ and
\begin{equation}\label{mu_s}
    D_s(\lambda)
:=(1-\lambda^2+\lambda^\alpha\cos(\pi s))^2
+\lambda^{2\alpha}\sin^2(\pi s),\qquad
\mu_s(\lambda):=\frac{2\sin(\pi s)}{\pi}
\frac{\lambda^{\alpha-1}}{D_s(\lambda)}=\frac{2\rho_s(\lambda^2)}{\lambda}.
\end{equation}
Combining~\eqref{eq:green mixture} with~\eqref{rho_s estimate}, by Tonelli's theorem, it follows that
\begin{equation}\label{green tail}
\int_0^\infty\mu_s=1,\qquad  T_s(y):=1-F_s(y)
=\int_0^\infty e^{-\lambda y}\mu_s(\lambda)\, d\lambda.
\end{equation}
Moreover, we observe that
\begin{equation}\label{mu_s limit}
    \mu_s(\lambda)\sim
\frac{2\sin(\pi s)}{\pi}\lambda^{\alpha-1}\quad\text{as } \lambda\to0,
\qquad
\mu_s(\lambda)\sim
\frac{2\sin(\pi s)}{\pi}\lambda^{\alpha-5}\quad\text{as } \lambda\to\infty.
\end{equation}
Thus $\int_0^\infty\lambda^2\mu_s(\lambda)\, d\lambda<\infty$. In particular,   
\begin{equation} \label{eq:density moments}
m_1:=\int_0^\infty\lambda\mu_s(\lambda)\, d\lambda=2G_s(0),
\qquad
m_2:=\int_0^\infty\lambda^2\mu_s(\lambda)\, d\lambda=1.
\end{equation}
From~\eqref{G_s(0)}, we know that $m_1^2>m_2/2$.
\medskip

Our purpose in introducing the probability density $\mu_s$
is to obtain the distribusion inequality~\eqref{dfbfdgdfvdf} for $F_s$.
Indeed, the representation
\[
F_s(y)
=
1-\int_0^\infty e^{-\lambda y}\mu_s(\lambda)\,d\lambda
\]
reduces the problem to properties of $\mu_s$.
We shall establish the stronger differential inequality
\begin{equation}\label{dfbdfgdgdd}
    yF_s'(y)>F_s(y)\bigl(1-F_s(y)\bigr),
\qquad y>0.
\end{equation}
The estimate~\eqref{dfbfdgdfvdf} then follows by integration.

We first prove a second moment comparison with an exponential density and  use it to establish sufficient conditions for the
differential inequality~\eqref{dfbdfgdgdd}. These conditions are then verified for $\mu_s$.

\begin{Lemma}\label{lemma:three crossings}
Let $\nu$ be a probability density on $(0,\infty)$ with $m=\int_0^\infty\lambda\nu(\lambda)\, d\lambda>0$. 
Suppose that $y>0$, $ym<1$, and
\[
\int_0^\infty e^{\lambda y}\nu(\lambda)\, d\lambda=\frac1{1-ym}.
\]
If $\nu(\lambda)-m^{-1}e^{-\lambda/m}$ has exactly three sign changes,
in the order $-,+,-,+$, then
\[
\int_0^\infty\lambda^2\nu(\lambda)\, d\lambda<2m^2.
\]
\end{Lemma}
\begin{proof}
Since $e^{\lambda y}\geq \frac{(y\lambda)^2}2$, we observe that
\[\int_0^\infty\lambda^2\nu(\lambda)\, d\lambda\leq  2\int_0^\infty e^{\lambda y}y^{-2} \nu(\lambda)\, d\lambda=\frac{2}{(1-ym)y^2}<\infty. \]
Moreover, set $g_m(\lambda):=m^{-1}e^{-\lambda/m}$, then
\begin{equation}\label{sdvsdvd bdf}
    \int_0^\infty g_m\, d\lambda=1,\quad \int_0^\infty \lambda g_m\, d\lambda=m, \quad  \int_0^\infty \lambda^2g_m\, d\lambda=2m^2\, \quad\text{and}\quad \int_0^\infty e^{\lambda y} g_m\, d\lambda=\frac{1}{1-ym}.
\end{equation}
We denote $r(\lambda):=\nu(\lambda)-g_m(\lambda)$, it follows that
\begin{equation}\label{sdv bsjdvnk}
    \int_0^\infty r(\lambda)\, d\lambda=\int_0^\infty \lambda r(\lambda)\, d\lambda= \int_0^\infty e^{\lambda y} r(\lambda)\, d\lambda=0.
\end{equation}

Let $0<a_1<a_2<a_3$ be the sign changing points of $r$ and set 
\begin{equation}\label{sdsdjo}
    P(\lambda):= \sum_{i=1}^3e^{a_i y}\prod_{1\leq j\leq3,j\neq i}\frac{\lambda-a_j}{a_i-a_j}:=A\lambda^2+B\lambda+C.
\end{equation}

We claim that $A>0$ and $f(\lambda):=e^{\lambda y}-P(\lambda)$  changes its sign exactly three times on the half-line, in the order $-,+,-,+$.

Indeed, since $f(a_i)=0$ with $i=1,2,3$, by Rolle's theorem, there exists $\eta\in (a_1,a_3)$ such that $f''(\eta)=0$. This implies that $A=\frac{y^2 e^{\eta y}}2>0$.

On the other hand, set $w(t):=(t-a_1)(t-a_2)(t-a_3)$, and $c_\lambda:=\frac{f(\lambda)}{w(\lambda)}$ with $\lambda\neq a_i$ for each $i=1,2,3$.
Denote $h(t):=f(t)-c_\lambda w(t)$. Then $h(a_i)=h(\lambda)=0$, using Rolle's theorem, we can find $\eta_\lambda$ satisfying $h'''(\eta_\lambda)=0$. Hence $c_\lambda=\frac{y^3e^{\eta_\lambda y}}{6}$ and $f(\lambda)=\frac{y^3e^{\eta_\lambda y}}{6}(\lambda-a_1)(\lambda-a_2)(\lambda-a_3)$. 

\medskip

Consequently, we obtain $r(\lambda)(e^{y\lambda}-P(\lambda))\ge0,$ their product is strictly positive
on a set of positive measure.  Hence,
by combining~\eqref{sdv bsjdvnk} with~\eqref{sdsdjo}, we conclude that 
\[0< \int_0^\infty r(\lambda)(e^{y\lambda}-P(\lambda))\, d\lambda= \int_0^\infty e^{y\lambda} r(\lambda)\, d\lambda -\int_0^\infty (A\lambda^2+B\lambda+C)r(\lambda)\, d\lambda=-A \int_0^\infty \lambda^2r(\lambda)\, d\lambda. \]
As a consequence of this and~\eqref{sdvsdvd bdf}, taking into account the fact that $A>0$, we complete the proof.
\end{proof}

We now use this moment comparison to derive the differential
inequality for the Laplace transform of a probability density.

\begin{Lemma}\label{lemma:score}
Let $\mu\in C^2(0,\infty)$ be a probability measure and $\int_0^\infty \lambda^2\mu(\lambda)\, d\lambda<\infty$. Assume that:
\begin{itemize}
\item[(i)] $m_1^2>m_2/2$, where $m_1:=\int_0^\infty \lambda\mu(\lambda)\, d\lambda$ and $m_2:=\int_0^\infty \lambda^2\mu(\lambda)\, d\lambda$;
\item[(ii)] $\mu'(a)=0$ for some $a>0$,  $\mu'>0$ on $(0,a)$ and $\mu'<0$
on $(a,\infty)$;
\item[(iii)] $q=-(\ln\mu)'$, restricted to $(a,\infty)$, increases
strictly to one maximum and then decreases strictly vanishing at infinity;
\item[(iv)] $e^{\gamma\lambda}\mu(\lambda)\to\infty$ as
$\lambda\to\infty$, for every $\gamma>0$.
\end{itemize}
\medskip

Set $T(y)=\int_0^\infty e^{-\lambda y}\mu(\lambda)\, d\lambda$
and $F(y)=1-T(y)$. Then
\begin{equation}\label{eq:logistic}
yF'(y)>F(y)(1-F(y))\qquad \text{for every } y>0.
\end{equation}
\end{Lemma}

\begin{proof}
Let $l(y):=F'(y)=\int_0^\infty\lambda e^{-y\lambda}\mu(\lambda)\, d\lambda$ and
\[Z(y)=yl(y)-T(y)+T^2(y).\]
It suffices to show that $Z(y)>0$ for all $y>0$. 

To this end, by Taylor formula, we observe that  
\[
T(y)=1-m_1y+\tfrac12m_2y^2+o(y^2),\qquad
l(y)=m_1-m_2y+o(y).
\]
Consequently, we see that $Z(y)=(m_1^2-m_2/2)y^2+o(y^2)>0$ for $y>0$ sufficiently small.
\medskip

Argue by contradiction. Suppose $Z(y)=0$ at some $y>0$. 

We define a  probability measure $\nu$ by
\[
\nu(\lambda)=\frac{e^{-y\lambda}\mu(\lambda)}{T(y)},
\qquad m=\int_0^\infty\lambda\nu(\lambda)\, d\lambda=\frac{l(y)}{T(y)}.
\]
Since $Z(y)=0$, one has $ ym=1-T(y)\in(0,1)$ and 
\[
\int e^{y\lambda}\nu(\lambda)\,d\lambda
=\frac1{T(y)}=\frac1{1-ym}.
\]
Let $\gamma=1/m-y>0$. It is immediate that the sign of $r(\lambda):=\nu(\lambda)-m^{-1}e^{-\lambda/m}$ is the sign of $e^{\gamma\lambda}\mu(\lambda)-T(y)/m.$

We claim that $e^{\gamma\lambda}\mu(\lambda)-T(y)/m$ changes its sign at most three times, and the final sign is positive. 

Indeed,  by conditions (ii) and (iii), one has \[\left(\mu(\cdot) e^{\gamma \cdot}\right)'(\lambda)=\left(\gamma \mu(\lambda)+\mu'(\lambda)\right)e^{\gamma\lambda}>0\quad \forall \lambda\in(0,a),\]
 and for every $\lambda>a$, $\left(\ln[\mu(\cdot) e^{\gamma \cdot}]\right)'(\lambda)=\gamma-q.$
The hypothesis on $q$ therefore permits at most three 
monotonicity intervals, increasing, decreasing, and increasing. The assertion that the final sign is positive follows from assumption (iv).

Therefore, $r(\lambda):=\nu(\lambda)-m^{-1}e^{-\lambda/m}$ changes its sign at most three times.
\medskip

Since $\mu'(r)>0$ near the origin, one has that the limit $\mu_0:=\lim_{\lambda\downarrow 0}\mu(\lambda)$ exists.

We now prove $\mu_0<T(y)/m$. 

If $\mu_0\ge T(y)/m $, then  the initial signs of  both $r(\lambda)$ and $e^{\gamma\lambda}\mu(\lambda)-T(y)/m$ are positive. This implies that both $r(\lambda)$ and $e^{\gamma\lambda}\mu(\lambda)-T(y)/m$ have at most two null points. 

We observe that $\int_0^\infty r(\lambda)\, d\lambda=0$, thus  $r(\lambda) $ changes its sign at least once on $(0,+\infty)$. 

Recall that the final signs of  both $r(\lambda)$ and $e^{\gamma\lambda}\mu(\lambda)-T(y)/m$ are  positive,  there thus exist $a_1,a_2$ with $a_1\neq a_2$ such that $r(a_1)=r(a_2)=0$, $r(\lambda)>0$ on $(0,a_1)\cup(a_2,\infty)$, and $r(\lambda)<0$ on $(a_1,a_2)$. We choose a line $L(\lambda):=k\lambda+b$   such that $e^{a_i y}=L(a_i)$ with $i=1,2$. Set $\bar c_\lambda=\frac{e^{\lambda y}-L(\lambda)}{(\lambda-a_1)(\lambda-a_2)}$, then $\bar h(t)=e^{t y}-L(t)-\bar c_\lambda (t-a_1)(t-a_2)$ has three null points. By Rolle's theorem, it follows that $e^{\lambda y}-L(\lambda)=\frac{y^2 e^{\bar\eta y}}{2}(\lambda-a_1)(\lambda-a_2)$ for some $\bar\eta >0$. 

Based on this, we see that
\begin{equation}\label{sdvsdvsdvsd}
    \int_0^\infty r(\lambda)(e^{\lambda y}-L(\lambda))\, d\lambda>0.
\end{equation}
Recall that
\begin{equation}\label{sdv bsjdvnkf}
    \int_0^\infty r(\lambda)\, d\lambda=\int_0^\infty \lambda r(\lambda)\, d\lambda= \int_0^\infty e^{\lambda y} r(\lambda)\, d\lambda=0.
\end{equation}
From this and~\eqref{sdvsdvsdvsd}, we conclude 
\[0<\int_0^\infty r(\lambda)(e^{\lambda y}-L(\lambda))\, d\lambda=0.\]
This is impossible. Hence $\mu_0<T(y)/m$. The initial sign of $r$ is negative.

\medskip

Moreover, if $r$ changes its sign exactly once, then $\int_0^\infty r(\lambda)(\lambda-d_1)\, d\lambda>0$ where $r(d_1)=0$. This is contradict to the identity~\eqref{sdv bsjdvnkf}. 

As a result, $r$ has exactly three sign changes, in the order $-,+,-,+$. In view of Lemma~\ref{lemma:three crossings}, we obtain 
\[M_2:=\int_0^\infty\lambda^2\nu(\lambda)\, d\lambda<2m^2.\]
Taking the derivative of $Z$ at the null point $y$, we see that
\begin{equation}\label{derivative Z}
    Z'(y)=yl'(y)+2(1-T(y))l(y)=yT(y)(2m^2-M_2)>0.
\end{equation}

 Since $Z>0$ on an interval $(0,\delta)$ for some $\delta>0$, its first point of
vanishing $y_0>0$ is well defined.
This means that $Z(y)>0$ for $0<y<y_0$, and therefore
\[
Z'(y_0)
=
\lim_{h\downarrow0}
\frac{Z(y_0)-Z(y_0-h)}{h}
\leq0.
\]
This contradicts the strict positivity of $Z'$ at
every point where $Z$ vanishes.  This yields~\eqref{eq:logistic}.
\end{proof}

We next show that  $\mu_s$ defined in~\eqref{mu_s} satisfies the all assumptions in Lemma~\ref{lemma:score}, and deduce the desired estimate for $F_s$.

\begin{Lemma}\label{lemma F_s}
    Let $1/2\leq s<1$ and $F_s$ be defined in~\eqref{eq:green doubling}, then \begin{equation}\label{F_s}
    F_s(2y)>\frac{2F_s(y)}{1+F_s(y)}.
\end{equation}
\end{Lemma}
\begin{proof}
By~\eqref{eq:density moments}, we know that $m_1^2>\frac{m_2}{2}$. Also,
the asymptotic behavior in~\eqref{mu_s limit}  implies  that $e^{\gamma\lambda}\mu_s(\lambda)\to\infty
$ as $\lambda\to\infty$ for every $\gamma>0$.
It remains to verify the required properties of $\mu_s$ and its logarithmic derivative.

    Let $1\le\alpha=2s<2$. Set
\[
b=-\cos(\pi\alpha/2)\in[0,1),\quad
t=\lambda^2,\quad u=\lambda^\alpha,\quad
D=(t-1)^2+2bu(t-1)+u^2>0.
\]
Define the differential operator $E=2t\partial_t+\alpha u\partial_u$
and the polynomials
\begin{equation}\label{eq:polynomial definitions}
    H=(\alpha-1)D-ED,\quad
N=(ED)^2-D(E^2D)+D\,(ED)-(\alpha-1)D^2,\quad 
M=EN-5N.
\end{equation}
We notice that, along the curve $t=\lambda^2$, $u=\lambda^\alpha$, the operator $E$
acts as $\lambda\frac d{d\lambda}$. Thus, taking the derivative of $\ln\mu_s$
 with respect to $\lambda$, we obtain
\[
\frac{d}{d\lambda}(\ln\mu_s)
=\frac{H}{\lambda D},
\qquad
\frac{d^2}{d\lambda^2}(\ln\mu_s)
=\frac{N}{\lambda^2D^2}.
\]
By  the definition of $H$, a direct calculation yields
\begin{equation}\label{eq:Hmonotone}
    \frac H{\lambda^2}
=(\alpha-1)\lambda^{-2}+2b\lambda^{\alpha-2}+(6-2\alpha)-(\alpha+1)\lambda^{2\alpha-2}-6b\lambda^\alpha-(5-\alpha)\lambda^2.
\end{equation}
For every  $1<\alpha<2$, we see that $H/{\lambda^2}$ is strictly decreasing,  $ H/{\lambda^2}\to +\infty$ as $\lambda\to0$, and $ H/{\lambda^2}\to -\infty$ as $\lambda\to \infty$. 
For $\alpha=1$, one has $b=0$ and $H/{\lambda^2}=2-4\lambda^2$.

Hence, for every $1\le\alpha<2$, $H$ has exactly one null point $a>0$, with $H'(a)<0$.  It follows  that 
\begin{equation}\label{xbfbdxbf}
   \mu_s \text{ is strictly increasing on } (0,a) \text{ and decreasing on } (a,\infty) 
\end{equation}
From Lemma~\ref{lem:polynomial}, since $H<0$ on $(a,\infty)$,  we know 
that $M(\lambda)>0$ for every $\lambda>a$. 

Furthermore, we observe that
\[
N=\lambda \frac{dH}{d\lambda}D-H(D+ED),\qquad N(a)=aH'(a)D(a)<0,
\]
and, for $\lambda>a$,
\[
\frac{d}{ d\lambda}\left(\frac{N(\lambda)}{\lambda^5}\right)
=\frac{M(\lambda)}{\lambda^6}>0.
\]
Since $D\sim\lambda^4$, $N/\lambda^5\sim(5-\alpha)\lambda^3$ as $\lambda\to\infty$.
Hence there exists $\bar a>a$ such that  $N(\bar a)=0$, $N<0$ on $(a,\bar a)$, and $N>0$ on $(\bar a,\infty)$. Recall that  $q=-(\ln\mu_s)'$, then $q'=-N/(\lambda^2D^2)$. 
It follows that $q$ first increases
strictly to one maximum and then decreases strictly 
on $(a,\infty)$, and
\[
q>0,\qquad q(\lambda)\sim(5-\alpha)/\lambda\to0 \quad \text{as } \lambda\to \infty.
\]
As a result of this and~\eqref{xbfbdxbf},  $\mu_s$ satisfies the hypotheses of Lemma~\ref{lemma:score}, 
we thus obtain
\[
yF_s'(y)>F_s(y)(1-F_s(y)).
\]
Dividing by $F_s(1-F_s)>0$ and integrating from $y$ to $2y$,  we conclude that
\[
\ln\frac{F_s(2y)}{1-F_s(2y)}
-\ln\frac{F_s(y)}{1-F_s(y)}>\ln2,
\]
which proves~\eqref{F_s}.
\end{proof}

\subsection{An upper  estimate}\label{sec:sharp proof}

In this section,	we prove the sharp estimate~\eqref{sharp estimate 1} by a key  monotonicity formula that effectively handles the interplay between operators of different orders.  

Fix $s\in(0,1)$. For $y>0$, we denote
\begin{equation}\label{v_y}
    w_y(x)=[G_s\ast\mathbf1_{[-y,y]}](x)=\int_{x-y}^{x+y}G_s(t)\, dt,\qquad v_y=w_y/F_s(y).
\end{equation}
It follows that
\[
0<v_y\le1,\quad v_y(0)=1,\quad
w_y'(x)=G_s(x+y)-G_s(x-y)<0\quad\forall x>0,
\]
and  $w_y(y)=F_s(2y)/2$. In particular,
\begin{equation}
\label{eq:resolvent}
-w_y''+(-\Delta)^sw_y+w_y=
\mathbf 1_{[-y,y]}
\end{equation}
in the sense of distributions.
Since $\mathbf1_{[-y,y]}\in L^2(\R)$, one has $w_y\in H^2(\R)$. The Lipschitz bound on $G_s$ also implies $w_y\in C^{1,1}(\R)$.

\begin{Lemma}\label{lemma: s-extension w_y}
Fix $y>0$. Let $\mathscr{E}_s w_y$ be $s$-extension of $w_y$. Define 
\[
\mathcal H_y(x)=\frac12w_y'^2(x)-\frac12w_y^2(x)+
\frac{d_s}{2}\int_0^\infty
\tau^{1-2s}((\mathscr{E}_s w_y)_x^2-(\mathscr{E}_s w_y)_\tau^2)\,\mathrm d\tau.
\]
    Then $\mathcal H_y$ is locally absolutely continuous, \[
\mathcal H_y'(x)=
-\mathbf1_{(0,y)}w_y'(x)\qquad\text{ a.e. on }(0,\infty), \quad \text{and}\quad \mathcal H_y(x)=0
\qquad\text{for every }x\geq y. 
\]
\end{Lemma}
\begin{proof}
    Since $w_y\in C^{1,1}(\R)\cap C_b^1(\R)$ solves~\eqref{eq:resolvent}, the argument for the extension integral remains valid. Thus $\mathcal H_y$ is continuous, locally Lipschitz and
\[
\mathcal H_y'(x)=w_y'(x)w_y''(x)-w_y(x)w_y'(x)-(-\Delta)^sw_y(x)w_y'(x)=
-\mathbf1_{(0,y)}w_y'(x)\qquad\text{ a.e. on }(0,\infty).
\]
We observe that, since the identity holds almost everywhere, replacing
$\mathbf1_{(0,y)}$ by $\mathbf1_{(0,y]}$ does not change the conclusion.

For any $0\leq a<b$, the precise integrated identity is
\begin{equation} \label{eq:H integrated}
\mathcal H_y(b)-\mathcal H_y(a)
=
-\int_a^b\mathbf 1_{(0,y)}(x)w_y'(x)\,dx.
\end{equation}
This identity remains valid when the interval crosses $y$.

\medskip

We now show that $\mathcal H_y(x)=0$ for every $x\geq y$.

Indeed, by \eqref{eq:H integrated}, $\mathcal H_y$ is constant on
$[y,\infty)$. We now show that this constant is zero.

Since $w_y\in H^2(\R)\subset H^1(\R)$,
the extension energy identity gives
\[
d_s\int_{\R}\int_0^\infty
\tau^{1-2s}
\left(
|(\mathscr{E}_s w_y)_x|^2+|(\mathscr{E}_s w_y)_\tau|^2
\right)\,d\tau\,dx
=
\|(-\Delta)^{s/2}w_y\|_{L^2(\R)}^2
<\infty.
\]
Therefore, the nonnegative function
\[
E_y(x)=w_y^2(x)+w_y'(x)^2+
\int_0^\infty\tau^{1-2s}
\left(
|(\mathscr{E}_s w_y)_x(x,\tau)|^2
+
|(\mathscr{E}_s w_y)_\tau(x,\tau)|^2
\right)\,d\tau\in L^1(\R).
\]
There exists a sequence $R_j\to\infty$, with $R_j>y$, such that $E_y(R_j)\longrightarrow0.$
This implies that  $|\mathcal H_y(R_j)|\to 0$.
Since $\mathcal H_y$ is constant on $[y,\infty)$, this proves
\begin{equation}
\label{eq:H-outside}
\mathcal H_y(x)=0
\qquad\text{for every }x\geq y.
\end{equation}
This ends the proof.
\end{proof}

\begin{remark}
    The preceding argument implies the explicit expression
\[
\mathcal H_y(x)
=
\begin{cases}
w_y(y)-w_y(x),&0\leq x\leq y,\\
0,&x\geq y.
\end{cases}
\]
We observe that its
one-sided derivatives are
\[
\mathcal H_y'(y_-)=-w_y'(y),
\qquad
\mathcal H_y'(y_+)=0.
\]
In fact, since $G_s$ is strictly decreasing , $w_y'(y)=G_s(2y)-G_s(0)<0,$
thus these one-sided derivatives differ.
Accordingly, the derivative identity is stated almost everywhere,
rather than at the endpoint \(x=y\).
\end{remark}

\begin{Lemma}\label{lemma: A_s (v_y)}
Let $\mathcal A_s$ be defined in~\eqref{A_s}. Then 
    \begin{equation}\label{eq:interval axis}
\mathcal A_s(v_y)=\frac{2F_s(y)-F_s(2y)}{F_s^2(y)}-1,
\end{equation}
where $v_y$ is introduced in~\eqref{v_y}.
\end{Lemma}
\begin{proof}
  Taking $a=0$ and $b=y$ in~\eqref{eq:H integrated}, since $\mathcal H_y(y)=0$, we obtain
\begin{equation}\label{eq:H origin}
\mathcal H_y(0)=w_y(y)-w_y(0).
\end{equation}

On the other hand, since $w_y$ and   $\mathscr{E}_s w_y$ are even in $x$, one has $w_y'(0)=0$ and $(\mathscr{E}_s w_y)_x(0,\tau)=0 $ for every  $ \tau>0.$
By the definition of $\mathcal H_y$, recalling~\eqref{eq:H origin}, it follows that
\begin{equation}\label{eq:H origin energy}
w_y(y)-w_y(0)=\mathcal H_y(0)
=
-\frac12w_y^2(0)
-\frac12\mathcal A_s(w_y).
\end{equation}
By the definition of $w_y$, we conclude that
\begin{equation}
\label{eq:axis energy w}
\mathcal A_s(w_y)
=2F_s(y)-F_s(2y)-F_s^2(y).
\end{equation}

Since  $v_y=\frac{w_y}{F_s(y)}$, one has 
\[\mathcal A_s(v_y)
=
d_s\int_0^\infty
\tau^{1-2s}|( \mathscr{E}_s  {v_y})_\tau(0,\tau)|^2\,d\tau=
\frac{d_s}{F_s^2(y)}
\int_0^\infty
\tau^{1-2s}|(\mathscr{E}_s w_y)_\tau(0,\tau)|^2\,d\tau=
\frac{\mathcal A_s(w_y)}{F_s^2(y)},\]
which is exactly~\eqref{eq:interval axis}.
\end{proof}

\medskip

\begin{Lemma}\label{lemma upper bound}
    Let $s\in[1/2,1), p>0$ and $Q\in \mathcal M_s$. Then
    \begin{equation}
        Q^{p}(0)<p+1.
    \end{equation}
\end{Lemma}
\begin{proof}
    Let  $f=Q^{p+1}$, then $f'(y)<0$ for $y>0$.
By the monotonicity and decay, we see 
\[
f(x)=f(|x|)=-\int_{|x|}^\infty f'(y)\, dy=-\int_0^\infty\mathbf1_{[-y,y]}(x)f'(y) \,dy,\qquad
\int_0^\infty -f'(y)\, dy=Q^{p+1}(0).
\]
By Tonelli's theorem in the resolvent equation, it follows that
\[
Q(x)=\int_{\R}G_s(x-y)f(y)\, dy=-\int_0^\infty w_y(x)f'(y)\,dy,\qquad
Q(0)=-\int_0^\infty F_s(y)f'(y)\, dy.
\]
We define a  probability measure $\nu$  by
\begin{equation}\label{eq:probability mixture1}
 d\nu(y)=\frac{-F_s(y)f'(y)}{Q(0)}\, dy,\qquad
\frac {Q(x)}{Q(0)}=\int_0^\infty v_y(x)\, d\nu(y),\qquad
\int_0^\infty\frac1{F_s(y)}\,d\nu(y)=Q^p(0).
\end{equation}
Thus, Tonelli's theorem  implies
\[
(\mathscr{E}_s [Q/Q(0)])(x,\tau)=
\int_{\R}P_s(x-z,\tau)
\left[
\int_0^\infty v_y(z)\,d\nu(y)
\right]dz=
\int_0^\infty (\mathscr{E}_s v_y)(x,\tau)\,d\nu(y).
\]

Let $\mathcal A_s$ be defined in~\eqref{A_s}.  We claim that 
\begin{equation}\label{vnkdsnjvckds}
    \mathcal A_s(Q/Q(0))
\leq
\int_0^\infty\mathcal A_s(v_y)\,d\nu(y)<Q^p(0).
\end{equation}
Indeed, from Lemma~\ref{lemma Poisson extension}, recall that $0<v_y\leq1$, we know that
\[
\int_0^\infty\int_{\R}
|\partial_\tau P_s(x-z,\tau)|\,|v_y(z)|\,dz\,d\nu(y)\leq
\int_0^\infty
\|\partial_\tau P_s(\cdot,\tau)\|_{L^1(\R)}
\,d\nu(y)=
\|\partial_\tau P_s(\cdot,\tau)\|_{L^1(\R)}
\leq\frac{1+2s}{\tau}<\infty.
\]
Thus, Fubini's theorem now yields
\begin{equation}\label{eq:derivative mixture}
\partial_\tau (\mathscr{E}_s [Q/Q(0)])(x,\tau)=
\int_{\R}\partial_\tau P_s(x-z,\tau)
\left[
\int_0^\infty v_y(z)\,d\nu(y)
\right]dz=
\int_0^\infty
\partial_\tau (\mathscr{E}_s v_y)(x,\tau)\,d\nu(y).
\end{equation}

We observe that,  fix $\tau>0$, $|\partial_\tau (\mathscr{E}_s v_y)(0,\tau)|\leq \frac{1+2s}{\tau}$  for every $y>0$. 
Thus, \[
\int_0^\infty
|\partial_\tau (\mathscr{E}_s v_y)(0,\tau)|^2\,d\nu(y)\leqslant \frac{(1+2s)^2}{\tau^2}.\]
As a consequence of this and Cauchy-Schwarz inequality, one has
\begin{equation}\label{eq:pointwise Jensen}
    |\partial_\tau (\mathscr{E}_s [Q/Q(0)])(0,\tau)|^2=
\left|\int_0^\infty \partial_\tau (\mathscr{E}_s v_y)(0,\tau)\cdot1\,d\nu(y)\right|^2
\leq
\int_0^\infty|\partial_\tau (\mathscr{E}_s v_y)(0,\tau)|^2\,d\nu(y).
\end{equation}
Multiplying~\eqref{eq:pointwise Jensen} by the positive weight
\(d_s\tau^{1-2s}\) and integrating over
\(0<\epsilon<T<\infty\), we obtain 
\[
\begin{aligned}
d_s\int_\epsilon^T
\tau^{1-2s}
|\partial_\tau (\mathscr{E}_s [Q/Q(0)])(0,\tau)|^2\, d\tau &\leq
d_s\int_\epsilon^T
\int_0^\infty
\tau^{1-2s}
|\partial_\tau (\mathscr{E}_s v_y)(0,\tau)|^2\,d\nu(y)\,d\tau\\
&=
\int_0^\infty
\left[
d_s\int_\epsilon^T
\tau^{1-2s}
|\partial_\tau (\mathscr{E}_s v_y)(0,\tau)|^2\,d\tau
\right]d\nu(y).
\end{aligned}
\]
By the monotone convergence theorem, we conclude that 
\[ \mathcal A_s(Q/Q(0))
\leq
\int_0^\infty\mathcal A_s(v_y)\,d\nu(y).\]
Moreover,  in the light of Lemmata~\ref{lemma F_s} and~\ref{lemma: A_s (v_y)}, we see that
\begin{equation}\label{dsvsdfg}
    0\leq\mathcal A_s(v_y)
<
\frac{1-F_s(y)}{1+F_s(y)}
<1
\qquad\text{for every }y>0.
\end{equation}
Set  $z(y)=F_s(y)^{-1}>1$  and $\Psi(z)=\frac{z-1}{z+1}$.  By~\eqref{eq:probability mixture1}, it follows that
\[
\int_0^\infty\mathcal A_s(v_y)\,d\nu(y)<\int_0^\infty z(y)\,d\nu(y)=Q^p(0).
\]
This proves~\eqref{vnkdsnjvckds}.

Since $\Psi''(z)=-4/(z+1)^3<0$, using~\eqref{vsdfvdsfcds} and ~\eqref{vnkdsnjvckds},  combining~\eqref{dsvsdfg}  with Jensen's inequality,  one has
\[
\frac{2Q^p(0)}{p+2}-1=\mathcal A_s(Q/Q(0))
\leq\int_0^\infty\mathcal A_s(v_y)\,d\nu(y)<\int_0^\infty\Psi(z(y))\,d\nu(y)\leq\Psi(Q^p(0))=\frac{Q^p(0)-1}{Q^p(0)+1}.
\]
Since $Q^p(0)>0$, this yields $Q^p(0)<p+1$ and completes  the proof.
\end{proof}

\subsection{Uniqueness}\label{sec:Nondegeneracy and uniqueness}

In this section, we complete the proof of Theorem~\ref{th:nondegeneracy s>1/2} by the sharp estimate developed in  Lemma~\ref{lemma upper bound}.

\subsubsection{Uniqueness}

\begin{Lemma}\label{lemma uniquness}
Let $0<s<1$ and $p>0$. 
Define \begin{equation}\label{n_s}
    \mathcal N_s:=\left\{Q\in\mathcal M_s:\, \frac{p+2}2<Q^p(0)\leq p+1 \right\}.
\end{equation}
Then any two positive centered solutions in $\mathcal N_s$ are identical.
\end{Lemma}
\begin{proof}
Let $u,v\in \mathcal N_s$, and denote their $s$-extensions by
$U$ and $V$, respectively. Assume that
\[
a:=u(0)\geq b:=v(0).
\]
For $t>\bigl((p+2)/2\bigr)^{1/p}$, define
\[
\mathcal R(t)=t\sqrt{\frac{2t^p}{p+2}-1}.
\]
A direct calculation yields 
\begin{equation}\label{eq:R derivatives}
\mathcal R(t)\mathcal R'(t)=t(t^p-1),\qquad
\mathcal R''(t)=\frac{pt^{p-1}(t^p-p-1)}
{(p+2)(2t^p/(p+2)-1)^{3/2}}.
\end{equation}

Recall that $\mathcal B_s(u,u)=\mathcal A_s(u)$ as defined~\eqref{A_s}. From~\eqref{vsdfvdsfcds}, we know that
\begin{equation}\label{hdfgbfdgvds}
    \mathcal B_s(u,u)=\frac{2a^{p+2}}{p+2}-a^2=\mathcal R^2(a)>0, \qquad \mathcal  B_s(v,v)=\frac{2b^{p+2}}{p+2}-b^2=\mathcal R^2(b)>0.
\end{equation}
Consider the following representation
\[\mathcal C_{u,v}=u'v'-uv+d_s\int_0^\infty \tau^{1-2s}\left(U_xV_x-U_\tau V_\tau\right)\, d\tau.\]

We claim that \begin{equation}\label{kdsncmklsd}
    \mathcal C'_{u,v}=-u'v^{p+1}-v'u^{p+1}\in L^1(\R).
\end{equation}
Indeed, since $u,v\in C^2(\R)\cap C^1_b(\R)$, Lemma~\ref{lemma Poisson extension} implies that
for every compact interval $I\subset\R$,
\begin{equation}\label{eq:weighted majorant3}
\tau^{1-2s}\bigl(|U_x V_x|+|U_\tau V_\tau|\bigr)(x,\tau)
\leq
\begin{cases}
C_I\bigl(\tau^{1-2s}+\tau^{2s-1}\bigr),
&0<\tau\leq1,\\[2mm]
C_I\tau^{-1-2s},
&\tau\geq1,
\end{cases}
\qquad x\in I.
\end{equation}
Define
\[
\mathcal J(x)=d_s\int_0^\infty \tau^{1-2s}\left(U_xV_x-U_\tau V_\tau\right)\, d\tau
\]
and, for \(0<\epsilon<1<T\),
\[
\mathcal J_{\epsilon,T}(x)=
d_s\int_\epsilon^T \tau^{1-2s}\left(U_xV_x-U_\tau V_\tau\right)\, d\tau.
\]
By \eqref{eq:weighted majorant3}, it follows that $\mathcal J_{\epsilon,T}$ converges to $\mathcal J$ uniformly on $I$ as $\epsilon\to 0$, $T\to \infty$.

Moreover, it is immediate that
\begin{equation}\label{dscvksdnkl}
    -d_s\epsilon^{1-2s}U_\tau(\cdot,\epsilon)V_x(\cdot,\epsilon)\to (-\Delta)^su\, v', \qquad -d_s\epsilon^{1-2s}V_\tau(\cdot,\epsilon)U_x(\cdot,\epsilon)\to (-\Delta)^sv\,u'
\end{equation}
locally uniformly as  $\epsilon\to0$. By the extension equation~\eqref{eq:extension energy},
it follows that 
\begin{equation}
    \begin{split}\label{sdvksbvdsv}
        \mathcal J_{\epsilon,T}'(x)
&=
-d_sT^{1-2s}\left(U_\tau(x,T)V_x(x,T)+V_\tau(x,T)U_x(x,T)\right)\\
&\qquad \qquad+
d_s\epsilon^{1-2s}\left(
U_\tau(x,\epsilon)V_x(x,\epsilon)+V_\tau(x,\epsilon)U_x(x,\epsilon)\right).
    \end{split}
\end{equation}
By Lemma~\ref{lemma Poisson extension}, for $T\geq1$,
\[
T^{1-2s}|U_\tau(x,T)V_x(x,T)+V_\tau(x,T)U_x(x,T)|
\leq
C_s\|u \|_{L^\infty(\R)}\|v \|_{L^\infty(\R)}T^{-1-2s}.
\]
Hence, integrating identity~\eqref{sdvksbvdsv} over an arbitrary interval $[x_1,x_2]$
with $x_1<x_2$, using~\eqref{dscvksdnkl}, passing to the limit, we conclude that $\mathcal J$ is locally Lipschitz continuous  and
\begin{equation}
\label{eq:mathcal J integral identity}
\mathcal J(x_2)-\mathcal J(x_1)
=
-\int_{x_1}^{x_2}
\left[
(-\Delta)^su\,v'+(-\Delta)^sv\,u'
\right]dx.
\end{equation}
The fundamental theorem of calculus yields
\[
\mathcal J'(x)=-[(-\Delta)^su(x)v'(x)+(-\Delta)^sv(x)u'(x)].
\]
As a result, since both $u$ and $v$ solve equation~\eqref{main equation},  we obtain the desired result~\eqref{kdsncmklsd}.

Since $u,v\in H^1(\R)$,
the extension energy identity gives
\[
d_s\int_{\R}\int_0^\infty
\tau^{1-2s}
\left(
|U_x|^2+|U_\tau|^2+ |V_x|^2+|V_\tau|^2
\right)\,d\tau\,dx
=
\|(-\Delta)^{s/2}u\|_{L^2(\R)}^2+\|(-\Delta)^{s/2}v\|_{L^2(\R)}^2
<\infty.
\]
Consequently, using the decay of $u,v,u',v'$ at infinity, there exists a sequence $R_j\to \infty$ such that $ \mathcal C_{u,v}(R_j)\to 0$,  as $j\to \infty$. 

Moreover, we observe that
$ \mathcal C_{u,v}(0)=-ab -\mathcal B_s(u,v)$. Thus,  
\[\mathcal B_s(u,v)=-ab+\int_{0}^\infty \mathcal C_{u,v}'(x)\, dx=-ab+\int_{0}^\infty[-u'v^{p+1}-v'u^{p+1}]\, dx.\]

Suppose that $v\neq u$. Then 
integrating  by parts, it follows that
\[\int_{0}^\infty[-u'v^{p+1}-v'u^{p+1}]\, dx=ab^{p+1}+\int_{0}^\infty[u^{p+1}-(p+1)u v^{p}](-v')\, dx.\]
Owing to strict convexity of $t^{p+1}$ on $(0,+\infty)$, one has 
\[u^{p+1}-(p+1)u v^{p}\geq -pv^{p+1}\]
where the equality holds if and only if $u=v$. Recall that $-v'>0$, we conclude that
\[\mathcal B_s(u,v)=-ab+ab^{p+1}+\int_{0}^\infty[u^{p+1}-(p+1)u v^{p}](-v')\, dx>-ab+ab^{p+1}-\frac p{p+2}b^{p+2}.\]
Recall that $\mathcal R^2(b)=\frac{2b^{p+2}}{p+2}-b^2$ and $\mathcal R(b)\mathcal R'(b)=b(b^p-1)$, thus we see that
\[\mathcal B_s(u,v)>-ab+ab^{p+1}-\frac p{p+2}b^{p+2}=\mathcal R^2(b)+(a-b)\mathcal R(b)\mathcal R'(b).\]

On the other hand, by combining~\eqref{hdfgbfdgvds} and Cauchy-Schwarz inequality, we derive 
\[\mathcal B_s(u,v)\leqslant \sqrt{\mathcal B_s(u,u)\mathcal B_s(v,v)}=\mathcal R(a)\mathcal R(b).\]
Since $\mathcal R(b)>0$, we obtain
\begin{equation}\label{fsdfsdsd}
    \mathcal R(a)> \mathcal R(b)+(a-b)\mathcal R'(b)
\end{equation}

If \(a=b\), then $\mathcal R(a)=\mathcal R(b)$, which contradicts~\eqref{fsdfsdsd}. 

If $a>b$ and $b^p<a^p\le p+1$, recalling~\eqref{eq:R derivatives}, one has $\mathcal R$ is concave on $[b,a]$, thus
\[
\mathcal R(a)\le \mathcal R(b)+(a-b)\mathcal R'(b).
\]
This is impossible, thanks to~\eqref{fsdfsdsd}. Therefore, we conclude that $u=v$.
\end{proof}

\bigskip

\appendix

\section{Properties of the Stieltjes density}
This section devoted to the identities~\eqref{rho_s identity} and~\eqref{rho_s estimate} by Cauchy's integral formula on a
truncated keyhole contour. 

\begin{Proposition}\label{prop: rho_s}
Fix $s\in(0,1)$. Define
\[Z_s(\zeta)=\frac{1}{1+\zeta+\zeta^s}\qquad \zeta\in \mathbb C.\]
Let $r>0$. We denote  by $\rho_s$ the positive density in the Stieltjes representation of $Z_s$:
\[
\rho_s(r)=\frac{\sin(\pi s)}{\pi}
\frac{r^s}
{\left(1-r+r^s\cos(\pi s)\right)^2
+r^{2s}\sin^2(\pi s)}.
\]
Then $\rho_s(r)>0$, and
\begin{equation}\label{eq:representation}
Z_s(x)
=\int_0^\infty\frac{\rho_s(r)}{x+r}\,dr,
\qquad x\ge0.
\end{equation}
The integral is absolutely convergent, including at $x=0$. 
In particular,
\begin{equation}\label{fdbfdb df}
    \int_0^\infty\frac{\rho_s(r)}r\, dr=1,
\qquad
\int_0^\infty\rho_s(r)\, dr=1.
\end{equation}
\end{Proposition}

\begin{proof}
Let
\[\Omega=\mathbb C\setminus(-\infty,0],\qquad  \zeta=|\zeta|e^{i\theta} \quad \text{with }
\quad
\theta:=\operatorname{Arg}\zeta\in (-\pi,\pi).
\]
Define $\zeta^s=|\zeta|^se^{is\theta}$.
We claim that  $F_s$ is holomorphic on $\Omega$. 

Indeed,  we observe that
\begin{equation}
    \operatorname{Im}(1+\zeta+\zeta^s)=
    \begin{cases}
|\zeta|\sin\theta+|\zeta|^s\sin(s\theta)>0\qquad\text{if } 0<\theta<\pi\\
|\zeta|\sin\theta+|\zeta|^s\sin(s\theta)>0\qquad\text{if } -\pi<\theta<0\\
    \end{cases}
\end{equation}
As for $\theta=0$, we see that $1+\zeta+\zeta^s=1+|\zeta|+|\zeta|^s>0.$
Thus $ 1+\zeta+\zeta^s\ne0 $ for all $
\zeta\in\Omega.$
This means that $Z_s$ is holomorphic throughout the slit plane. Since $x>0$ and $Z_s(x)>0$, the factor
$1/(\zeta-x)$ yields a simple pole at $\zeta=x$.
This is the only singularity of $Z_s(\zeta)/(\zeta-x)$ inside the contour $\Omega_{\epsilon,R,\delta}$ defined in~\eqref{sdfvsn}.

\medskip

 Denote $A_s(r)=1-r+r^s\cos(\pi s),$ and $ B_s(r)=r^s\sin(\pi s)>0. $
Since $B_s(r)>0$, a direct calculation yields that 
\begin{equation}\label{bdvsjvds}
    Z_{s,+}(-r):=\lim_{\theta\uparrow\pi}Z_s(r e^{i\theta})=\frac{1}{A_s(r)+iB_s(r)}
=\frac{A_s(r)-iB_s(r)}{A_s(r)^2+B_s(r)^2},
\end{equation}
and 
\begin{equation}\label{dsvsdjv}
Z_{s,-}(-r):=\lim_{\theta\downarrow-\pi}Z_s(r e^{i\theta})=\frac{1}{A_s(r)-iB_s(r)}=\frac{A_s(r)+iB_s(r)}{A_s(r)^2+B_s(r)^2}.
\end{equation}
Hence, we obtain 
\begin{equation}\label{sdvsdvsdvdds}
    Z_{s,-}(-r)-Z_{s,+}(-r)
=\frac{2iB_s(r)}{A_s(r)^2+B_s(r)^2}
=2\pi i\,\rho_s(r).
\end{equation}
\medskip

We next prove that the integrability of the density $\rho_s$.

Indeed, we notice that 
\[
\rho_s(r)\sim\frac{\sin(\pi s)}{\pi}r^s
\qquad\text{as } r\to 0\quad  \text{and}\quad \rho_s(r)\sim\frac{\sin(\pi s)}{\pi}r^{s-2}
\qquad\text{as } r\to\infty.
\] 
Since the density $\rho_s$ is continuous on $(0,\infty)$,  we derive that 
\begin{equation}\label{sdfvsdvs}
    \int_0^\infty \rho_s(r)\,\mathrm dr<\infty,
\qquad
\int_0^\infty \frac{\rho_s(r)}r\,\mathrm dr<\infty.
\end{equation}
In particular, for every $x>0$,  one has $0\le\frac{\rho_s(r)}{x+r}
\le\frac{\rho_s(r)}r.$
This proves the integrability of the desired identity~\eqref{eq:representation}.
\bigskip

We now turn to prove  the desired identity~\eqref{eq:representation} by  Cauchy’s integral formula.

For fixed $x>0$, choose $0<\epsilon<x/2$, $R>2x$, and
$0<\delta<\pi/2$.
Set $\theta_\delta=\pi-\delta$, and consider the truncated sector
\begin{equation}\label{sdfvsn}
    \Omega_{\epsilon,R,\delta}=\left\{
\zeta\in\mathbb C:
\epsilon<|\zeta|<R,\,
-\theta_\delta<\operatorname{Arg}\zeta<\theta_\delta
\right\}.
\end{equation}
We notice that, $\overline{\Omega_{\epsilon,R,\delta}}\subset\Omega$, and $x$ is an interior point. The positively oriented boundary consists of four pieces, as follows:
\begin{itemize}
\item
The outer arc $\Gamma_1:=\left\{\zeta\in \Omega:\, \zeta=R e^{i\theta} \text{ with } \theta \text{ increasing from }   -\theta_\delta\text{ to }\theta_\delta\right\}$. 
\item
The upper radial segment $\Gamma_2:=\left\{\zeta\in \Omega:\, \zeta=r e^{i\theta_\delta} \text{ with } r \text{ decreasing from }   R\text{ to }\epsilon \right\}$.  
\item
The inner arc $\Gamma_3:=\left\{\zeta\in \Omega:\, \zeta=\epsilon e^{i\theta} \text{ with } \theta \text{ decreasing from }   \theta_\delta \text{ to }-\theta_\delta \right\}$.   
\item
The lower radial segment $\Gamma_4:=\left\{\zeta\in \Omega:\, \zeta=r e^{-i\theta_\delta} \text{ with } r \text{ increasing from }  \epsilon \text{ to }R \right\}$. 
\end{itemize}
Since $Z_s$ is holomorphic in a neighborhood of $\overline{\Omega_{\epsilon,R,\delta}}$, Cauchy's integral formula yields 
\begin{equation}\label{vskjdfnvk}
   \begin{split}
        2\pi i\,Z_s(x)
&=\int_{\Gamma_1}\frac{Z_s(\zeta)}{\zeta-x}\, d\zeta+\int_{\Gamma_2}\frac{Z_s(\zeta)}{\zeta-x}\, d\zeta+\int_{\Gamma_3}\frac{Z_s(\zeta)}{\zeta-x}\, d\zeta+\int_{\Gamma_4}\frac{Z_s(\zeta)}{\zeta-x}\, d\zeta\\
&=I_{1,R,\epsilon}+I_{2,\delta}+I_{3,R,\epsilon}+I_{4,\delta}.
   \end{split}
\end{equation}

As for $I_{1,R,\epsilon}$, we choose $R$  sufficiently large such that $1+R^s\le\frac R2.$  Thus, one has 
\begin{equation}
    \begin{split}
      |I_{1,R,\epsilon}| \le \int_{\Gamma_1}\frac{1}{|\zeta-x||1+\zeta+\zeta^s|}\, d\zeta\le \int_{\Gamma_1}\frac{2}{(R-x)R}\, d\zeta
\le\frac{8\pi R}{R^2}
=\frac{8\pi}{R}\to 0 \quad \text{as } R\to \infty.  
    \end{split}
\end{equation}

\medskip

As for $I_{3,R,\epsilon}$, we  take $\epsilon$ sufficiently small such that $\epsilon+\epsilon^s\le\frac12.$ It follows that 
\begin{equation}
    \begin{split}
      |I_{3,R,\epsilon}| \le \int_{\Gamma_3}\frac{1}{|\zeta-x||1+\zeta+\zeta^s|}\, d\zeta\le \int_{\Gamma_3}\frac{2}{(x-\epsilon)}\, d\zeta
\le\frac{8\pi \epsilon}{x}\to 0 \quad \text{as } \epsilon\to 0.  
    \end{split}
\end{equation}

\medskip

We now turn to consider $I_{2,\delta}$ and $I_{4,\delta}$. 
Let $\epsilon$ and $R$ be chosen, and  $\delta\to0$, we see that
\[
I_{2,\delta}=\int_{\Gamma_2}\frac{Z_s(\zeta)}{\zeta-x}\, d\zeta=\int_R^\epsilon
\frac{Z_s(r e^{i\theta_\delta})}
{r e^{i\theta_\delta}-x}
\,e^{i\theta_\delta}\,dr.
\]
Let $m_s:=\sup_{\theta\in[\pi/2,\pi]}\sin(s\theta)>0$. Since $\theta_\delta\in (\pi/2,\pi)$, it follows that $|1+re^{i\theta_\delta}+r^se^{is\theta_\delta}|\geqslant m_sr^s$, and $|r e^{i\theta_\delta}-x|^2\geq \frac{(r+x)^2}2$ for every $x>0$. Hence
\[\left|\frac{Z_s(r e^{i\theta_\delta})e^{i\theta_\delta}}
{r e^{i\theta_\delta}-x}\right|\leqslant \frac{\sqrt{2}}{(x+r)r^sm_s}\in L^1(0,\infty). \]
Using the dominated convergence theorem, we derive that
\begin{equation}\label{dsvbsdkjnk}
\lim_{\delta\to0}I_{2,\delta}=
-\int_R^\epsilon
\frac{Z_{s,+}(-r)}{-r-x}\,dr=
-\int_\epsilon^R
\frac{Z_{s,+}(-r)}{r+x}\, dr.  
\end{equation}

Similarly, the lower radial integral satisfies
\begin{equation}\label{sdvnsdkivns}
\lim_{\delta\downarrow0}I_{4,\delta}=\int_\epsilon^R
\frac{Z_{s,-}(-r)}{r+x}\,dr.
\end{equation}
Recalling~\eqref{sdvsdvsdvdds}, combining~\eqref{dsvbsdkjnk} with~\eqref{sdvnsdkivns}, we obtain
\begin{equation}
    \lim_{\delta\downarrow0}
\left(I_{2,\delta}+I_{4,\delta}\right)=
\int_\epsilon^R
\frac{Z_{s,-}(-r)-Z_{s,+}(-r)}{r+x}\, dr=
2\pi i\int_\epsilon^R
\frac{\rho_s(r)}{r+x}\, dr.
\end{equation}

Let $\epsilon$ and $R$ be taken, and  $\delta\to0$.
Using dominated convergence theorem, together with the
bounds on the two circular arcs, we obtain 
\[
\left|
Z_s(x)-\int_\epsilon^R
\frac{\rho_s(r)}{r+x}\,dr
\right|
\le \frac4R+\frac{4\epsilon}{x}.
\]
Owing to absolute integrability of
$\rho_s(r)/(r+x)$. Taking $R\to\infty$ and $\epsilon\to0$, we conclude that
\[
Z_s(x)=
\int_0^\infty\frac{\rho_s(r)}{x+r}\,dr
\qquad \text{for every } x>0.
\]
This proves~\eqref{eq:representation} for positive real arguments.
\medskip

Finally, utilizing~\eqref{sdfvsdvs}, the dominated convergence theorem implies 
\[
1=\int_0^\infty\frac{\rho_s(r)}r\, dr.
\] 
Moreover, by multiplying ~\eqref{eq:representation} by $x>0$:
\[
\frac{x}{1+x+x^s}
=
\int_0^\infty
\frac{x}{x+r}\rho_s(r)\, dr.
\]
We observe that the left-hand side tends to $1$ as $x\to\infty$. Since  $\rho_s\in L^1(0,\infty)$, the
Monotone convergence theorem therefore yields the desired identity~\eqref{fdbfdb df}.
This completes the proof.
\end{proof}

\section{A polynomial calculation}\label{app:polynomial}

This parts proves the algebraic sign statement used to determine
the logarithmic derivative shape of the mixing density. All signs follow from explicit polynomial identities and coefficient estimates.

\begin{Lemma}\label{lem:polynomial}
Let $1\le\alpha\le2$, $0\le b\le1$, and $t,u>0$, with
\begin{equation}\label{dsdfsdsd}
\min\{t,t^2\}\le u^2\le\max\{t,t^2\}.
\end{equation}
Assume that
\[
D=(t-1)^2+2bu(t-1)+u^2>0,
\]
and define $ E=2t\partial_t+\alpha u\partial_u, $
\[
H=(\alpha-1)D-ED,\qquad
N=(ED)^2-D(E^2D)+D\,(ED)-(\alpha-1)D^2,\qquad
M=EN-5N.
\]
If $H\le0$, then $M\ge0$, with strict inequality whenever
$H<0$ or $\alpha>1$.
\end{Lemma}

\begin{proof}
Let $v=bu\ge0$, $w=u^2>0$, and
\[
C=-4\alpha^2+12\alpha-5=4-(2\alpha-3)^2\in[3,4].
\]
We observe that,
$Et=2t$, $Ev=\alpha v$, and $Ew=2\alpha w$.
A direct calculation gives
\begin{equation}\label{eq:expanded H}
H=(\alpha-5)t^2+(6-2\alpha)t+\alpha-1
+2(1-3t)v-(\alpha+1)w.
\end{equation}

\bigskip

We claim that $H>0$ if $0<t\le1/3$. 

 Indeed, recalling~\eqref{dsdfsdsd}, one has  $w\le t$, and $1-3t\ge0$. Thus
\[
H\ge(\alpha-5)t^2+(5-3\alpha)t+\alpha-1.
\]
If $\alpha=1$, then $H\geq 2t(1-2t)>0$. If $\alpha=2$, then $H\geq1-t-3t^2\ge1/3>0$.
This implies that $H$ is positive for every $1\le\alpha\le2$. 

\medskip

We consider the case of $1/3<t\le1$.
By directly applying the definition of $M$, we have
\begin{equation} \label{eq:M expansion}
M=D\mathcal F-H\mathcal J,
\end{equation}
where
\begin{align*}
\mathcal F&=E^2H-6EH+(C+5)H\\
&=4(\alpha-1)
\left[(\alpha-2)(5-\alpha)t^2
-2(\alpha-2)(3-\alpha)t+\alpha(3-\alpha)\right]+6(2-\alpha)\left[\alpha+(4-3\alpha)t\right]v,\\
\mathcal J&=E^2D-4ED+CD\\
&=C(1+t^2)+2(2\alpha-3)^2t+2(\alpha-1)\left[3\alpha-5+3(3-\alpha)t\right]v
+(4\alpha-5)w.
\end{align*}
As for $\mathcal F$, we notice that  the function $I_{1,\alpha}(t):=(\alpha-2)(5-\alpha)t^2
-2(\alpha-2)(3-\alpha)t+\alpha(3-\alpha)$ is a concave quadratic in $t$. Thus $I_{1,\alpha}(t)\geqslant \min\{I_{1,\alpha}(0),I_{1,\alpha}(1)\}=2$  for all $t\in[0,1]$.
Furthermore, $I_{2,\alpha}(t)=\alpha+(4-3\alpha)t\geqslant \min\{I_{2,\alpha}(0),I_{2,\alpha}(1)\}\ge 0$  for all $t\in[0,1]$.  Therefore,
\[
\mathcal F\ge8(\alpha-1)\ge0.
\]
Since $1\ge t>1/3$ and $4\alpha-5\ge-1$, one has $3\alpha-5+3(3-\alpha)t\ge2(\alpha-1)\ge0$ and $w\le t$. Thus, we see that
$(4\alpha-5)w\ge-w\ge-t.$ Recall that $C\ge3$,  we obtain
\[
\mathcal J\ge3(1+t^2)-t>0.
\]
From~\eqref{eq:M expansion}, we know that if $H\le0$ then  $M\ge0$ since $D>0$.

If $H<0$, then $-H\mathcal J>0$. If $\alpha>1$, then $D\mathcal F>0$.
This proves the required strict inequality in either case.
\medskip

Finally, for every $t>1$, set $x=\alpha-1\in[0,1]$ and $z=t-1\ge0$.
The exact expansion is
\begin{equation}
\label{eq:M-expansion}
M=A_0+B_1v+B_2w+B_{vv}v^2+B_{vw}vw+B_{ww}w^2,
\end{equation}
with coefficients
\begin{align*}
A_0&=(12-3x)z^4+(44-8x)z^3+64z^2+32z,\\
B_1&=(24+16x-2x^2-2x^3)z^3
+(72+40x-4x^2)z^2+(80+48x)z+32,\\
B_2&=(2+12x+12x^2-8x^3)z^2
+(12+8x+16x^2)z+16+16x,\\
B_{vv}&=(12+24x)z^2+(24+48x)z+32x,\\
B_{vw}&=(24x+18x^2-6x^3)z+8+8x+20x^2,\\
B_{ww}&=-2+7x+4x^2.
\end{align*}
Note that every coefficient of $z$ in $A_0,B_1,B_{vv},B_{vw}$ is
nonnegative for all $0\le x\le1$.
Furthermore, we have $B_{ww}\ge-2$, and
\[B_2-2t^2=4x(3+3x-2x^2)z^2
+(8+8x+16x^2)z+14+16x>0.\]
From~\eqref{dsdfsdsd}, one has $0<w\le t^2$, thus
\[
M\geq B_2w+B_{ww}w^2
\ge w(B_2-2w)\ge w(B_2-2t^2)>0.
\]
This yields that $M>0$ in this case without the hypothesis $H\le0$.
We thus complete the proof.
\end{proof}

\section{Properties of $s$-Poisson extension}

Let $0<s<1$, and define the $s$-Poisson kernel by
\[
P_s(z,\tau)
=c_s\frac{\tau^{2s}}{(z^2+\tau^2)^{(1+2s)/2}},  \quad c_s=\frac{\Gamma(s+\tfrac12)}
{\sqrt{\pi}\,\Gamma(s)}
\qquad \tau>0.
\]
For a bounded measurable function $f$, its $s$-Poisson extension  is defined by
\[(\mathscr{E}_s f)(x,\tau)
=
\int_{\R}P_s(z,\tau)f(x-z)\,dz.
\]
\begin{Lemma}\label{lemma Poisson extension}
  Let $f\in  C_b^1(\R)\footnote{$C_b^1(\mathbb R):=
\{f\in C^1(\mathbb R):f,f'\in L^\infty(\mathbb R)\}$.}
\cap C_{\mathrm{loc}}^{1,1}(\R)$.  Then
\begin{equation}\label{dvsdvdsd}
    \left|\partial_x (\mathscr{E}_s f)(x,\tau)\right|
\leq \|f'\|_{L^\infty(\R)}
\qquad
\forall x\in\R,\ \tau>0.
\end{equation}
Moreover, for every compact interval $I\subset\R$, 
there exists a constant $C_{s,I,f}>0$ such that
\begin{equation}\label{dcdsvhidsjjv}
    \left|\partial_\tau (\mathscr{E}_s f)(x,\tau)\right|
\leq C_{s,I,f}\tau^{2s-1}
\qquad \forall x\in I,\ 0<\tau\leq 1.
\end{equation}
For $\tau\geq 1$, we also have
\[
\left|\partial_x (\mathscr{E}_s f)(x,\tau)\right|
+
\left|\partial_\tau (\mathscr{E}_s f)(x,\tau)\right|
\leq
\frac{C_s\|f\|_{L^\infty(\R)}}{\tau},
\qquad \forall x\in\R.
\] 

In particular,  \begin{equation}\label{vsdknvmlds}
    \sup_{x\in I}
\int_0^\infty
\tau^{1-2s}
\left(
|\partial_x (\mathscr{E}_s f)(x,\tau)|^2
+
|\partial_\tau (\mathscr{E}_s f)(x,\tau)|^2
\right)\,d\tau
<\infty.
\end{equation}
\end{Lemma}
\begin{proof}
A direct differentiation yields
\[
\partial_z P_s(z,\tau)
= -c_s(1+2s)
\frac{\tau^{2s}z}
{(z^2+\tau^2)^{(3+2s)/2}}\quad \text{and} \quad \partial_\tau P_s(z,\tau)=
c_s\tau^{2s-1}
\frac{2sz^2-\tau^2}
{(z^2+\tau^2)^{(3+2s)/2}}.
\]
Since $|2sz^2-\tau^2|
\leq (1+2s)(z^2+\tau^2),$ we obtain
\begin{equation} \label{eq:kernel vertical bound}
|\partial_\tau P_s(z,\tau)|
\leq
c_s(1+2s)
\frac{\tau^{2s-1}}
{(z^2+\tau^2)^{(1+2s)/2}}.
\end{equation}
Fix $\tau_0>0$. For $\tau\in[\tau_0/2,2\tau_0]$, both
kernel derivatives are bounded by
\[
C_{s,\tau_0}(1+|z|)^{-1-2s},
\]
where the constant is independent of $z$ and $\tau$.
Indeed, the derivatives are uniformly bounded for $|z|\leq1$,
since $\tau$ stays away from zero, and decay at least as fast as
$|z|^{-1-2s}$ for $|z|>1$.
Since this bound is integrable over $\R$, the dominated
convergence theorem allows us to differentiate under the integral sign.

Since $\int_{\R}P_s(z,\tau)\,dz=1$, 
differentiating the normalization identity yields
\begin{equation}\label{eq:constant cancellation}
\int_{\R}\partial_\tau P_s(z,\tau)\,dz=0.
\end{equation}
Furthermore, we observe that $\partial_\tau P_s(\cdot,\tau)$ is even.
Hence, for every \(r>0\),
\begin{equation}\label{eq:linear cancellation}
\int_{-r}^{r}z\,\partial_\tau P_s(z,\tau)\,dz=0.
\end{equation}
The latter integral is over a bounded interval and therefore is absolutely convergent.

\medskip

We first prove~\eqref{dvsdvdsd}. 
Using the representation in which the derivative falls on \(f\), we obtain
\[
\partial_x (\mathscr{E}_s f)(x,\tau)
=
\int_{\R}P_s(z,\tau)f'(x-z)\,dz.
\]
By combining the positivity with normalization of the kernel, one has 
\[ |\partial_x (\mathscr{E}_s f)(x,\tau)|\leq
\int_{\R}
P_s(z,\tau)|f'(x-z)|\,dz
\leq\|f'\|_{L^\infty(\R)}
\int_{\R}P_s(z,\tau)\,dz
=\|f'\|_{L^\infty(\R)}.
\]

\medskip

We now turn to show~\eqref{dcdsvhidsjjv}.
Fix a compact interval $I\subset\R$, and set
\[
I_1=\{x+z:x\in I,\ |z|\leq1\},
\qquad  L_I \text{ is the Lipschitz constant  of } f'  \text{ on } I_1.
\]
Since $f\in C_{\mathrm{loc}}^{1,1}(\mathbb R)$,
the constant $L_I$ is finite.
By~\eqref{eq:constant cancellation} and~\eqref{eq:linear cancellation},
\begin{equation}
   \begin{split}
       \partial_\tau (\mathscr{E}_s f)(x,\tau)&=
\int_{\R}
\partial_\tau P_s(z,\tau)
\left(f(x-z)-f(x)\right)\,dz\\
&=\int_{|z|\leq1}
\partial_\tau P_s(z,\tau)
\bigl(f(x-z)-f(x)+zf'(x)\bigr)\,dz+
\int_{|z|>1}
\partial_\tau P_s(z,\tau)
\bigl(f(x-z)-f(x)\bigr)\,dz.
   \end{split} 
\end{equation}
For $x\in I$ and $|z|\leq1$, Taylor's formula with integral remainder implies
\[
|f(x-z)-f(x)+zf'(x)|
=\left|
-z\int_0^1\bigl(f'(x-\theta z)-f'(x)\bigr)\,d\theta\right|\leq \frac{L_I}{2}|z|^2.
\]
Using \eqref{eq:kernel vertical bound}, we find
\begin{equation}
    \begin{split}\label{eq:local contribution}
        \left|
\int_{|z|\leq1}
\partial_\tau P_s(z,\tau)
\bigl(f(x-z)-f(x)+zf'(x)\bigr)\,dz
\right|&\leq
\frac{c_s(1+2s)L_I}{2}\tau^{2s-1}
\int_{|z|\leq1}
\frac{z^2}{(z^2+\tau^2)^{(1+2s)/2}}\,dz\\
&\leq \frac{c_s(1+2s)L_I}{2(1-s)}\tau^{2s-1}.
    \end{split}
\end{equation}
For the far-field contribution, one has
\begin{equation}\label{eq:far contribution}
 \left|
\int_{|z|>1}
\partial_\tau P_s(z,\tau)
\bigl(f(x-z)-f(x)\bigr)\,dz
\right|\leq \frac{2c_s(1+2s)}{s}
\|f\|_{L^\infty(\R)}\tau^{2s-1}.
\end{equation}
Combining~\eqref{eq:local contribution} with~\eqref{eq:far contribution},
we conclude that
\[
|\partial_\tau (\mathscr{E}_s f)(x,\tau)|
\leq C_{s,I,f}\tau^{2s-1},
\]
where 
\[
C_{s,I,f}=c_s(1+2s)
\left(
\frac{L_I}{2(1-s)}
+
\frac{2\|f\|_{L^\infty(\R)}}{s}
\right).
\]
In particular, this estimate holds uniformly for
$x\in I$ and $0<\tau\leq1$.

\medskip

We next estimate the derivatives at large heights by
differentiating the convolution kernel.
Since $P_s(\cdot,\tau)$ is even and decreasing on $(0,\infty)$, we obtain
\begin{equation}\label{bsxvcsd}
    |\partial_x (\mathscr{E}_s f)(x,\tau)|\leq
\int_{\R}
|\partial_z P_s(x-y,\tau)||f(y)|\,dy\leq\frac{2c_s}{\tau}\|f\|_{L^\infty(\R)}.
\end{equation}
Similarly, \eqref{eq:kernel vertical bound} and the
normalization of $P_s$ yield
\begin{equation}\label{vdsvsdvdsvs}
    |\partial_\tau (\mathscr{E}_s f)(x,\tau)|
\leq
\int_{\R}|\partial_\tau P_s(x-y,\tau)||f(y)|\,dy\leq
c_s(1+2s)\tau^{2s-1}
\int_{\R}
\frac{\|f\|_{L^\infty(\R)}dz}{(z^2+\tau^2)^{(1+2s)/2}}=
\frac{1+2s}{\tau}\|f\|_{L^\infty(\R)}.
\end{equation}
Both~\eqref{bsxvcsd} and~\eqref{vdsvsdvdsvs}  hold for every $\tau>0$;
we use them in particular when $\tau\geq1$.

\medskip

Finally, for $x\in I$ and $0<\tau\leq1$, the estimates~\eqref{dvsdvdsd} and~\eqref{dcdsvhidsjjv} imply
\[ \tau^{1-2s}
\left(
|\partial_x (\mathscr{E}_s f)(x,\tau)|^2
+
|\partial_\tau (\mathscr{E}_s f)(x,\tau)|^2
\right)\leq
\|f'\|_{L^\infty(\R)}^2\tau^{1-2s}
+C_{s,I,f}^2\tau^{2s-1}.
\]
For $\tau\geq1$, we have
\[
\tau^{1-2s}
\left(
|\partial_x (\mathscr{E}_s f)(x,\tau)|^2
+
|\partial_\tau (\mathscr{E}_s f)(x,\tau)|^2
\right)\leq
\left(4c_s^2+(1+2s)^2\right)
\|f\|_{L^\infty(\R)}^2\tau^{-1-2s}.
\]
Both of these are integrable at zero and infinity, respectively, independent of $x\in I$. This implies the desired result~\eqref{vsdknvmlds}.
\end{proof}

\medskip

\section{Positive semi-definiteness on the odd subspace}

Our goal is to prove that the quadratic form of $L_{s,Q}$ is nonnegative on odd functions for every $Q\in\mathcal M_s$. 

We first determine the essential spectrum of the odd restriction. Moreover, an energy comparison and the strict negativity of $Q'$ on the positive half-line then exclude negative eigenvalues through orthogonality.

\subsection{Essential spectrum on the odd subspace}

We first determine the spectrum of $H=-\Delta+(-\Delta)^s+1$ on the odd subspace and show that the corresponding restriction of $L_{s,Q}$ has essential spectrum $[1,\infty)$. Hence, any negative spectral value must  be an isolated eigenvalue. This observation will be used in the proof of nonnegativity.

\begin{Lemma}\label{lemma specturm of H}
Let $s\in (0,1)$ and $p>0$.   Let $H_{\rm odd}$ be the restriction of $H=-\Delta+(-\Delta)^s+1$ to $L^2_{\rm odd}(\R) $. 
    Then \begin{equation}\label{specturm of H }
      \sigma_{\rm ess}(H_{\rm odd})=  \sigma(H_{\rm odd})
=
[1,\infty).
    \end{equation}
Moreover,  let $Q\in \mathcal M_s$, and  $L_{s,Q}=H-(p+1)Q^p$ be its corresponding linearized operator.  Let $L_{s,Q}|_{L^2_{\mathrm{odd}}(\mathbb R)}$
denote the self-adjoint realization of $L_{s,Q}$ on the odd subspace. Then
\begin{equation}\label{3.9}
    \sigma_{\mathrm{ess}}(L_{s,Q}|_{L^2_{\mathrm{odd}}(\mathbb R)})=
[1,\infty).
\end{equation}
\end{Lemma}
\begin{proof} 
Recall that the unitary Fourier transform is defined by
\[
\widehat f(\xi)
=
\frac{1}{\sqrt{2\pi}}
\int_{\mathbb R}e^{-ix\xi}f(x)\,dx.
\]
It follows that $\widehat{Hf}(\xi) =m_s(\xi)\widehat f(\xi)$ with $m_s(\xi)=1+\xi^2+|\xi|^{2s}.$  Since $1+\xi^2\le m_s(\xi)\le 2(1+\xi^2),$ the operator domain of $H$ is $H^2(\mathbb R)$. Furthermore, since $m_s$ is even and $H$  preserves the odd subspace, we see  that $L^2_{\rm odd}(\mathbb R)$   is a reducing subspace and
$H_{\rm odd}$ is self-adjoint on
\[
D(H_{\rm odd})
=
H^2(\mathbb R)\cap L^2_{\rm odd}(\mathbb R).
\]
Owing to the Fourier sine transform, we know that $H_{\mathrm{odd}}$
is unitarily equivalent to multiplication by $m_s(\xi) $ with $ \xi>0.$  From Weyl's criterion, it follows that the essential range of $m_s$ is $[1,\infty)$.
Hence
\[
\sigma_{\mathrm{ess}}(H_{\mathrm{odd}})
=\sigma(H_{\mathrm{odd}})
=[1,\infty).
\]

We next verify that $\sigma_{\mathrm{ess}}(L_{s,Q}|_{L^2_{\mathrm{odd}}(\mathbb R)})
=[1,\infty).$

Since $L_{s,Q}|_{L^2_{\mathrm{odd}}(\mathbb R)}=H_{\rm odd}-(p+1)Q^p$ is self-adjoint on
$D(H_{\mathrm{odd}})$, and  $Q^p$ is bounded, even and $Q^p(x)\to 0$ as $|x|\to \infty$, the Weyl's theorem guarantees that it suffices to show that
\begin{equation}\label{dsvsdvcsc}
    -(p+1)Q^p(H_{\rm odd}+1)^{-1} \text{ is compact }.
\end{equation}
Let $h_k $ be a bounded sequence in $L_{\rm odd}^2(\R)$, and $u_k=(H_{\rm odd}+1)^{-1}h_k$. Hence $u_k$ is  uniformly bounded in $H^2(\R)\cap L_{\rm odd}^2(\R)$,
 passing to a subsequence,
we may therefore assume that $u_k\rightharpoonup u
$ in $H^2(\mathbb R), $ and $
u_k\to u $ in $
L^2_{\mathrm{loc}}(\mathbb R).$ Furthermore, 
 for every $R>0$, we have
\[
\|Q^p(u_k-u)\|_{L^2(\mathbb R)}^2
\le
\|Q^p\|_{L^\infty}^2
\|u_k-u\|_{L^2(-R,R)}^2+
\sup_{|x|>R}|Q^p(x)|^2
\|u_k-u\|_{L^2(\mathbb R)}^2.
\]
For fixed $R$, the first term tends to zero as
$k\to\infty$. We then take $R\to\infty$ and obtain $Q^pu_k\to Q^pu $ in $ L^2(\mathbb R).$
This proves~\eqref{dsvsdvcsc}, and we complete the proof.
\end{proof}

\subsection{Nonnegativity of the quadratic form}

Based on the essential spectrum identified, we use an energy
comparison and the strict negativity of $Q'$ on the positive half-line to rule out negative eigenvalues in the odd subspace.

\begin{Lemma}\label{lemma: nonnegative}
Let  $0<s<1$, $p>0$, and  $Q\in \mathcal M_s$, then the corresponding linearized operator $L_{s,Q}$ is nonnegative on $H^1(\R)\cap L^2_{\rm{odd}}(\R)$:
\begin{equation}\label{sdvkidsjmf}
    \langle f,L_{s,Q}f\rangle\geq 0
\qquad
\text{for every }
f\in H^1(\R)\cap L^2_{\rm{odd}}(\R).
\end{equation}
\end{Lemma}
\begin{proof}
From~\cite[Lemma A.1]{SZZ2026}, we know that  $Q\in H^2(\R)\subset C^{1,1/2}(\R)$.
This allows us to differentiate the equation  in the distribution sense, 
\begin{equation}\label{vdsvdsdsfds}
    -(Q'')'+(-\Delta)^sQ'+Q'
=(p+1)Q^pQ'
\end{equation}
and hence $L_{s,Q}Q'=0.$ Since $Q$ is even, $Q'$ is odd. Moreover,
$Q'(x)<0$ for every $x>0$ by Proposition~\ref{prop: Q'>0}.

We observe that  $Q'\in H^1(\mathbb R)$ and $Q\in L^\infty(\R)$, applying the same  argument developed in~\cite[Lemma A.1]{SZZ2026} to the equation~\eqref{vdsvdsdsfds}, obtaining that  $Q'\in H^2(\mathbb R).$
Thus $Q'$ is an eigenfunction  of   operator $L_{s,Q}$ corresponding to the
eigenvalue zero.

\medskip

First of all, we rewrite the representation of $\langle f,L_{s,Q}f\rangle$.

Indeed,  for simplify, set $V(x)=1-(p+1)Q(x)^p$, let $f\in H^1(\R)\cap L^2_{\rm{odd}}(\R),$ one has 
\[
\langle f,L_{s,Q}f\rangle=
\int_{\mathbb R}
\left(
|f'(x)|^2+V(x)|f(x)|^2
\right)\,dx+
\frac{c_{1,s}}{2}
\int_{\R}\int_{\mathbb R}
\frac{|f(x)-f(y)|^2}
{|x-y|^{1+2s}}
\,dx\,dy.
\]
We denote the restriction of $f$ to the positive half-line by
\[
w=f|_{(0,\infty)}\in H^1(0,+\infty).
\]
Since $H^1(\mathbb R)\subset C^{1/2}(\mathbb R)$ and
$f$ is odd, it follows that $f(0)=0,$
and hence $w(0)=0.$
Thus $w\in H^1_0(0,+\infty)$ in the trace sense. 

Hence  the local part is immediately reduced to the half-line:
\[
\int_{\mathbb R}
\left(|f'|^2+Vf^2
\right)\,dx
=2\int_0^\infty
\left(
|w'|^2+Vw^2
\right)\,dr.
\]
As for  the fractional part, since $H^1(\R)\subset H^s(\R)$, Tonelli's theorem allows us to
split the integral over the four quadrants of $\mathbb R^2$, we obtain the half-line representation below
\begin{equation}
  \begin{split}
       \frac{c_{1,s}}2\iint_{\mathbb R^2}
\frac{(f(x)-f(y))^2}
{|x-y|^{1+2s}}
\,dx\,dy
&=  \frac{c_{1,s}}2\left(2\int_{0}^\infty\int_{0}^\infty
\frac{(w(r)-w(t))^2}
{|r-t|^{1+2s}}
\,dr\,dt+2\int_{0}^\infty\int_{0}^\infty\frac{(w(r)+w(t))^2}
{(r+t)^{1+2s}}
\,dr\,dt\right).
  \end{split}
\end{equation}
Since $(w(r)+w(t))^2+(w(r)-w(t))^2
=2w^2(r)+2w^2(t),$
we obtain
\[
\int_{0}^\infty\int_{0}^\infty\frac{(w(r)+w(t))^2+(w(r)-w(t))^2}
{(r+t)^{1+2s}}
\,dr\,dt=
4\int_0^\infty
w^2(r)
\left(
\int_0^\infty
\frac{dt}{(r+t)^{1+2s}}
\right)dr=\frac{2}{s}
\int_0^\infty
r^{-2s}w^2(r)\,dr.
\]
It follows that
\begin{equation}\label{sdcfdsvds}
    \begin{split}
        \frac12\langle f,L_{s,Q}f\rangle
=&
\int_0^\infty
\left(
|w'(r)|^2
+
V(r)w^2(r)
+\frac{c_{1,s}}{s}r^{-2s}w^2(r)
\right)\,dr
\\
&+
\frac{c_{1,s}}{2}
\int_{0}^\infty\int_{0}^\infty
(w(r)-w(t))^2
\left(
\frac{1}{|r-t|^{1+2s}}
-
\frac{1}{(r+t)^{1+2s}}
\right)
\,dr\,dt.
    \end{split}
\end{equation}
Note  that $|r-t|<r+t$ if $r,t>0$ whenever $r\neq t$.

\medskip

Moreover,  for every $w\in H^1(0,\infty)$ with $w(0)=0$,  using the fundamental theorem of calculus, we observe that
\begin{equation}\label{vdsvsdvsdvs}
    \int_0^\infty r^{-2s}w^2(r)\,dr<\infty.
\end{equation}

Next, we consider a absolute value function. Given $w=f|_{(0,\infty)}$, and let $f^*$ be the odd extension of $|w|$:
\[f^*(x) = |w(x)| \quad \text{for } x>0, \qquad f^*(x) = 0 \quad \text{for } x=0,\qquad f^*(x) = -|w(-x)| \quad \text{for } x<0.\]
Then $f^*\in H^1(\R)\cap L^2_{\rm{odd}}(\R)$, and  $\|f^*\|_{L^2(\mathbb R)}
=\|f\|_{L^2(\mathbb R)}.$

We claim  that
\begin{equation}\label{3.7}
    \langle f^*,L_{s,Q}f^*\rangle
\leq
\langle f,L_{s,Q}f\rangle.
\end{equation}
Indeed, since the absolute value map $t\mapsto |t|$ is Lipschitz,
the Sobolev chain rule implies that $|w|\in H^1(0,\infty)$ and $(|w|)'=\operatorname{sgn}(w)w' $ a.e. on $(0,\infty).$
Here we set $\operatorname{sgn}(0)=0$ and use the
standard fact that $w'=0 $ a.e. on $
\{w=0\}.$ Hence,
\[
\left|(|w|)'\right|=|w'|
\qquad\text{a.e. on }(0,\infty),
\]
and all the terms in the local part of \eqref{sdcfdsvds} are unchanged
when $w$ is replaced by $|w|$.

For the double integral term, we use the elementary inequality $\bigl(|a|-|b|\bigr)^2\leq(a-b)^2.$ Therefore, we obtain~\eqref{3.7}. 
\medskip

Furthermore, it is immediate that 
\begin{equation}\label{ksdjsdf}
    \text{if } w \text{ changes sign on }  (0,\infty), \text{ then }
\langle f^*,L_{s,Q}f^*\rangle
<\langle f,L_{s,Q}f\rangle.
\end{equation}
Indeed, if $ab<0$, then  $(a-b)^2-\bigl(|a|-|b|\bigr)^2
=4|ab|>0.$  Also, we  see that    
\[
\frac{1}{|r-t|^{1+2s}}
-\frac{1}{(r+t)^{1+2s}}
>0\qquad \text{for every } r,t>0,
\]
thus, the decrease in the double integral is strictly positive whenever $w(r)$ and $w(t)$ have opposite signs on sets of positive measure.

\bigskip

With the preliminary work at hand, we now show the nonnegativity of the quadratic form on the odd subspace.

Let $\mathcal L:=L_{s,Q}|_{L^2_{\mathrm{odd}}(\mathbb R)}$
denote the self-adjoint realization of $L_{s,Q}$ on the odd subspace. Its
operator domain is $D(\mathcal L)
=H^2(\R)\cap L^2_{\mathrm{odd}}(\mathbb R)$, 
while the domain of the associated quadratic form $\mathfrak B_Q[f]=
\langle f, \mathcal Lf\rangle$   is $D(\mathfrak B_Q) =H^1(\R)\cap L^2_{\mathrm{odd}}(\mathbb R).$

Suppose, for contradiction, that
\[
\lambda_1^{\mathrm{odd}}
:=\inf\sigma(\mathcal L)<0.
\]
Since $\sigma_{\mathrm{ess}}(\mathcal L)=[1,\infty)$,
the lowest point of the spectrum is therefore an isolated eigenvalue.
We choose a corresponding real valued eigenfunction $\phi\in H^2(\R)\cap L^2_{\mathrm{odd}}(\mathbb R)$ with $\|\phi\|_2=1$.

If $\phi$ changed sign on $(0,\infty)$, let $\phi^*$
be the odd extension of $|\phi|$ from the positive
half-line. Utilizing the strict absolute value comparison~\eqref{ksdjsdf}, we see that
\[
\mathfrak B_Q[\phi^*]
<
\mathfrak B_Q[\phi]
=\lambda_1^{\mathrm{odd}},
\qquad
\|\phi^*\|_2=1,
\]
contradicting the variational characterization of
$\lambda_1^{\mathrm{odd}}$.

Hence, we may assume that $\phi\geq0$ on $(0,\infty)$,  changing the sign of $\phi$ if necessary. 
Since $\phi$ is nontrivial and $Q'<0$ on this half-line, one derives that
$\langle\phi,Q'\rangle_{L^2(\mathbb R)}
=
2\int_0^\infty\phi(x)Q'(x)\,dx<0.$
This is impossible, since $\phi$ and $Q'$ correspond to the distinct eigenvalues $\lambda_1^{\mathrm{odd}}<0$ and $0$
of the self-adjoint operator $\mathcal L$, they must  be orthogonal.

Hence $\mathfrak B_Q$ is nonnegative on
$H^1(\R)\cap L^2_{\rm{odd}}(\R)$. This ends the proof.
\end{proof}

\section*{AI usage statement}

Here, ChatGPT was used to discuss possible approaches to the uniqueness proof and to assist with some mathematical calculations and intermediate arguments. Some ideas developed through these discussions were incorporated into the proof. All mathematical arguments and calculations presented in this paper were carefully examined and independently verified by the authors, who take full responsibility for the results.

\section*{Conflict-of-Interest Statement}
The authors declare that they have no known competing financial interests or personal
relationships that could have appeared to influence the work reported in this paper.

\bibliographystyle{is-abbrv}

\bibliography{manuscript}

	\end{document}